\documentclass[hidelinks, 11pt, oneside, reqno]{amsart}
\usepackage{dsliheader}
\usepackage{cite}
\usepackage{caption}
\usepackage{enumitem}

\newtheoremstyle{named}{}{}{\itshape}{}{\bfseries}{.}{.5em}{\thmnote{#3}}
\theoremstyle{named}

\renewcommand{\thefootnote}{\fnsymbol{footnote}}

\newcommand{\sqg}{\sqrt{1+|\nabla \varphi(x)|^2}}
\newcommand{\wGk}{\widetilde{\Gamma}_j}

\DeclareMathOperator{\Id}{Id}
\DeclareMathOperator{\divg}{div}
\DeclareMathOperator{\I}{I}
\DeclareMathOperator{\II}{II}
\DeclareMathOperator{\III}{III}

\DeclareMathOperator{\Graph}{Graph}
\DeclareMathOperator{\dist}{dist}
\DeclareMathOperator{\Tan}{Tan}

\title[Generic boundary behavior of Dirichlet eigenfunctions]{Generic boundary behavior of Dirichlet eigenfunctions}

\author[Payne]{Alec Payne}
\address{Williams College, Department of Mathematics, 18 Hoxsey Street, Williamstown, MA 01267}
\email{ap30@williams.edu}
\author[Zhao]{Zihui Zhao}
\address{Johns Hopkins University, Department of Mathematics, 404 Krieger Hall, Baltimore, MD 21218}
\email{zhaozh@jhu.edu}

\begin{document}

\begin{abstract}
We extend Uhlenbeck's classical results on the generic properties of Laplacian eigenfunctions. In particular, we find a residual set of $C^k$ domains in $\mathbb{R}^n$, $k > n$, such that the nodal set of any Dirichlet eigenfunction intersects the boundary orthogonally along a $C^{k-1}$ embedded submanifold, whenever the intersection is non-empty. This is achieved by analyzing the boundary asymptotics of Dirichlet eigenfunctions via homogeneous harmonic polynomials and comparing their respective nodal sets. Along the way, we give a complete proof that generic Euclidean domains have simple Dirichlet eigenfunctions with zero as a regular value in the interior and zero as a regular value of the boundary normal derivative. In the special case of planar domains, a different argument provides one-parameter families of perturbed domains with the desired properties.
\end{abstract}

\thanks{The authors started working together while in residence at the Simons Laufer Mathematical Sciences Institute (formerly MSRI), partially supported by the National Science Foundation Grant No. DMS-1928930. The second-named author is also supported by the NSF Grant No. DMS-2350351.}

\maketitle

\tableofcontents

\section{Introduction}

In this paper, we establish the generic boundary behavior of Dirichlet Laplacian eigenfunctions on Euclidean domains. We say that $u$ is a Dirichlet Laplacian eigenfunction of a domain $\Om \subset \bR^n$ if
\begin{equation}
\begin{aligned}\label{equation Dirichlet eigenfunction problem}
\left\{\begin{array}{ll}
		-\Delta u = \lambda u, & \text{ in } \Omega, \\
		\quad\,\,\,\, u=0, & \text{ on } \partial\Omega,
	\end{array} \right.
\end{aligned}
\end{equation}
for some $\la \in \bR$.
It is a general principle that the analytic properties of Dirichlet Laplacian eigenfunctions improve when the domain or the metric is generically perturbed~\cite{UhlenbeckGeneric, Albert, Micheletti-1, Micheletti-2, HenryPerturbationbook05}. Here, generic means that these properties hold for a residual subset of the moduli space of domains or metrics. 
 
What previous work has not addressed is the generic boundary behavior of eigenfunctions (see for example the last paragraph of \S 3 in ~\cite{UhlenbeckGeneric}). The main result of our paper is to show that transversal intersection, in fact orthogonal intersection, of the interior nodal set with the boundary is generic for Dirichlet Laplacian eigenfunctions.

\begin{thm}\label{thm:main}
Let $n\geq 2$ and $k> n$. Let $M$ be a smooth compact connected $n$-manifold with boundary which admits a $C^k$ embedding into $\bR^n$, and let $\mathrm{Emb}_k(M, \bR^n)$ be the Banach manifold of $C^k$ embeddings.

There is a residual set $\mathfrak{X} \subset \mathrm{Emb}_k(M, \bR^n)$ such that if $F \in \mathfrak{X}$ and if $u$ is a Dirichlet Laplacian eigenfunction of the domain $\Omega = \inte F(M)$, then the following properties hold:
\begin{quote}
    The closure of the interior nodal set $N_u := \overline{\{u=0\} \cap \Omega}$ is a compact, possibly disconnected, embedded hypersurface which intersects $\partial\Omega$ orthogonally, and its manifold boundary $\partial N_u$, if not empty, is a $C^{k-1}$-embedded submanifold of $\partial\Omega$.
\end{quote}
\end{thm}

Theorem \ref{thm:main} is also true for generic domains in a complete ambient Riemannian manifold, but we work with Euclidean space for simplicity. We may also prove the analogous statement for generic $C^k$ metrics on a fixed Riemannian manifold $M$ with boundary. After our discussion of the proof techniques, we explain the necessary modifications to our argument in this case.

Our proof of Theorem \ref{thm:main} builds on previously established generic properties of Dirichlet Laplacian eigenfunctions, summarized in Theorem \ref{thm:generic}. These properties were stated in Uhlenbeck's seminal paper, where an outline of the proof was provided (see \cite[Example 3]{UhlenbeckGeneric}). We give a complete and slightly different proof of Theorem \ref{thm:generic} in Section \ref{sec:generic}, taking the Lagrangian rather than the Eulerian point of view to analyze PDEs. We use the same general framework of infinite-dimensional Sard--Smale and transversality theorems as Uhlenbeck, but we use a different PDE argument to verify the assumptions of the transversality theorems, inspired by Dan Henry's work on boundary value problems of PDEs under perturbations of domains~\cite{HenryPerturbationbook05}. We believe this approach to studying eigenfunctions under domain perturbations makes the key ingredients clearer and easier to adapt to other problems. See Section \ref{sec:generic} for further discussion and details.

The next step is to study the boundary behavior of Dirichlet eigenfunctions. By \cite{Bers}, Laplacian eigenfunctions are well approximated near their interior nodal set by homogeneous harmonic polynomials. In Section \ref{sec:bd}, we prove a boundary version of this statement (see Proposition \ref{prop:bdexp}), which is of independent interest for general boundary analysis of Dirichlet eigenfunctions. In the classical interior approximations, one gets $C^{k+2}$ error estimates if the coefficients of the relevant operator have $C^k$ regularity. However, in the boundary case, no matter how smooth the domain is, the approximations are at best $C^2$. In fact, even the $C^2$ estimate of the approximation is borderline and can only be understood as one-sided derivatives (see \eqref{est:grad2} and Remark \ref{rmk:reg}). Fortunately, $C^2$ regularity is exactly enough to characterize the generic boundary behavior.

The results of Sections \ref{sec:generic} and \ref{sec:bd} imply that generically, eigenfunctions are approximated by degree-$2$ homogeneous harmonic polynomials near the boundary (see Claim \ref{cl:P2}). The $C^2$ estimates from Proposition \ref{prop:bdexp} are enough to show that the nodal set of the eigenfunction and the harmonic polynomial are close (see Lemma \ref{lm:HS} for details).
Despite the freedom of degree-$2$ homogeneous harmonic polynomials in high dimensions (the space of such polynomials in $\bR^n$ has dimension $\frac{n(n+1)}{2} - 1$), the Dirichlet boundary condition implies that the relevant approximating polynomial must be the simplest $(n-2)$-invariant one. Using such a polynomial as the model for the nodal sets of Dirichlet eigenfunctions, we prove the geometric properties of nodal sets in Proposition \ref{prop:main}, from which Theorem \ref{thm:main} follows.

Theorem \ref{thm:main} can be easily adapted to hold for generic $C^k$ metrics on a fixed Riemannian manifold $M$ with boundary. That is, rather than working with generic domains in a fixed ambient manifold or Euclidean space, we may work with generic metrics on a fixed Riemannian manifold with boundary. The necessary adaptations to our argument are the following: prove the existence and requisite properties of the reduced Green's function and Poisson kernel for the Laplace--Beltrami operator $\Delta_g $, as in Section \ref{sec:prelim}, which are in turn used to prove the analogue of Theorem \ref{thm:generic} in this setting, or \cite[Theorem 8]{UhlenbeckGeneric}; apply the boundary unique continuation property for $\Delta_g + \lambda$ (or, for elliptic operators with Lipschitz coefficients) with vanishing boundary condition; prove the local boundary asymptotics for Dirichlet eigenfunctions of the operator $\Delta_g$ as in Proposition \ref{prop:bdexp}.

The case for planar domains (i.e.\ $n=2$) is very special, thanks to conformal maps and the full classification of homogeneous harmonic polynomials in $\bR^2$. In particular, nodal critical points are isolated, and at such a point, finitely many nodal curves meet at equal angles, assuming the planar domain has $C^{1,\alpha}$ regularity or is convex (see \cite{HartmanWintner, Albert, Alessandrini, SpectralPartition} etc.). We can take advantage of these special features to prove a stronger result in dimension $2$: The properties in Theorem \ref{thm:main} hold for one-parameter families of domains along a generic perturbation direction.

\begin{prop}\label{prop:2D}
    Suppose $\lambda$ is a simple Dirichlet eigenvalue of a $C^{2,\alpha}$ bounded planar domain $\Omega$. There is an open dense set $\mathcal{A}_{\lambda} \subset C^{2,\alpha}(\overline{\Omega}, \mathbb{R}^2)$ such that for any $s\in \mathcal{A}_\lambda$, the one-parameter family of domains $\Omega_t := h_t(\Omega)$ given by the diffeomorphism $h_t = \id + ts: \overline{\Omega} \to \overline{\Omega_t}$ satisfies the following: 
    \begin{quote}
        For any $0<|t|<t_0$, the Dirichlet Laplacian $-\Delta_{\Omega_t}$ has a one-dimensional eigenspace with eigenvalue $\lambda_t$, $|\lambda_t - \lambda| \ll 1$; and the nodal curve of the corresponding eigenfunction $u_t$ consists of a finite union of closed curves in $\Omega_t$ and simple curves that start and terminate at $\partial\Omega_t$, with orthogonal intersection. 
    \end{quote}
    Here the admissible set $\mathcal{A}_\lambda$ is given explicitly in \eqref{def:good-perturbation}, and the parameter $t_0$ depends on $s\in \mathcal{A}_\lambda$ and the eigenspace corresponding to $\lambda$.
\end{prop}

An analogous statement for finitely many eigenfunctions immediately follows, but the argument cannot be extended to include all eigenfunctions, since the small parameter $t_0$ may degenerate.
\begin{cor}
    Suppose the first $N$ Dirichlet eigenvalues of a $C^{2,\alpha}$ bounded planar domain $\Omega$ are simple. There is an open dense set $\mathcal{A}_N \subset C^{2,\alpha}(\overline{\Omega}, \mathbb{R}^2)$ such that for any $s\in \mathcal{A}_N$, the properties in Proposition \ref{prop:2D} hold for the first $N$ Dirichlet eigenfunctions of a one-parameter family of domains $\Omega_t = h_t(\Omega)$ given by the diffeomorphism $h_t = \id + ts$ with $0<|t|< t_0$.
\end{cor}

We point out that we use a different argument  to obtain a stronger result in Proposition \ref{prop:2D} than what is possible using the transversality theorem: The transversality argument deployed in Theorem \ref{thm:main} only yields a \emph{residual set of domains} near $\Omega$ with desirable properties; moreover, applying that argument to a one-parameter family of domains $\Omega_t$ yields the desirable properties only \emph{for almost every $t$} (see for example \cite[Section 3]{UhlenbeckGeneric}). Nonetheless, it is worth mentioning that even though we obtain the regularity of eigenfunctions for $\Omega_t$ with all $0<|t|<t_0$, it is possible that the nodal curve has bifurcation behaviors for $t<0$ and $t>0$, and we give a concrete example to illustrate this phenomenon at the end of Section \ref{sec:2d}. Our proof idea is closer in spirit to \cite{Albert}, where the author studies the eigenfunctions of $L+\rho$ for a fixed second-order elliptic operator $L$ on a closed manifold $M$ and generic functions $\rho \in C^\infty(M)$. 

\medskip

\textbf{Disclosure of AI use}. Large language models were used to assist with proofreading. The manuscript was entirely written by the authors, who take full responsibility for its content.  

\section{Preliminaries}\label{sec:prelim}
Let $\Omega \subset \bR^n$ be a bounded $C^{2,\alpha}$ domain\footnote{For the existence and properties of the reduced Green's function, it suffices to consider domains with $C^{1,\alpha}$ regularity. The $C^{2,\alpha}$ assumption is only needed for the reduced Poisson kernel in Lemma \ref{lem:reduced-poisson-kernel}.}, and let $L= \Delta +\lambda$, which is a self-adjoint Fredholm operator of index zero. We prove the existence and properties of the \emph{reduced Green's function} for the operator $L$, which morally speaking satisfies $LG(x,\cdot) = -\delta_x$, just like the classical Green's function for the Laplacian.
However, since the operator $L$ is not invertible (when $\lambda$ is an eigenvalue), this is impossible without some modifications.
To that end, let $\{\psi_i\}_{i=1}^m$ be an arbitrary $L^2(\Omega)$-orthonormal basis of $\ker L$, and let $\Pi_{(\ker L)^\perp}$ denote the $L^2(\Omega)$-orthogonal projection onto $(\ker L)^\perp$.

\begin{lem}[Reduced Green's function]\label{lm:Green}
    For any $x\in \Omega$ there exists a unique weak solution $G(x,\cdot) \in W_0^{1,q}(\Omega) \cap (\ker L)^\perp$ for any $1<q<\frac{n}{n-1}$ such that\footnote{Clearly any kernel element $u\in \ker L$ annihilates the left hand side, so we need the projection operator $\Pi_{(\ker L)^\perp}$ on the right hand side to make sure they are compatible.}
    \begin{equation}\label{eq:Green}
        \int_\Omega G(x,\cdot) Lu = -\delta_x \circ \Pi_{(\ker L)^\perp}(u), \qquad \text{ for any } u\in W_0^{1,p}(\Omega).
    \end{equation}
    
    Moreover, we have
    \begin{enumerate}[label=\roman*)]
        \item $\|G(x,\cdot)\|_{W^{1,q}(\Omega)} \leq C$ with a uniform constant $C$ independent of $x$,
        \item $LG(x,\cdot) = \sum_{i=1}^m \psi_i(x) \psi_i$ in the sense of distributions\footnote{We may also write this PDE as $LG(x,\cdot) = \delta_x \circ \Pi_{\ker L}$, just to show that it is independent of the specific choice of the orthonormal basis of $\ker L$.} in $\Omega \setminus\{x\}$, and $G(x,\cdot) \in C^\infty(\Omega \setminus\{x\})$,
        \item $\|G(x,\cdot) - G^0(x,\cdot)\|_{W^{1,q}(\Omega)} \leq C$ with a uniform constant $C$ independent of $x$, where $G^0(x,\cdot)$ denotes the classical Green's function for $\De$ in $\Omega$  with pole at $x$,
        \item $G(x,y) = G(y,x)$.
    \end{enumerate}
\end{lem}
\begin{proof}
    By the Fredholm alternative, for any $f\in H^{-1}(\Omega) \cap (\ker L)^\perp$, there is a unique solution $v\in H_0^1(\Omega) \cap (\ker L)^\perp$ such that $Lv = f$; moreover, when $f$ is assumed to be in the space $W^{-1,p}(\Omega)$ for some $p>n$, the corresponding solution $v\in C(\Omega)$. 

    For any $f\in W^{-1,p}(\Omega) \subset H^{-1}(\Omega)$, let
    \[ f^\perp := f- \sum_{i=1}^m \langle f, \psi_i \rangle \psi_i \in W^{-1,p}(\Omega), \]
    and let $u_f \in W_0^{1,p}(\Omega)\cap (\ker L)^\perp$ denote the unique weak solution of $Lu_f = f^\perp$. We have that
    \[ \|u_f\|_{C^{0,\alpha'}(\Omega)} \lesssim \|u_f\|_{W^{1,p}(\Omega)} \lesssim \|f^\perp\|_{W^{-1,p}(\Omega)} \lesssim \|f\|_{W^{-1,p}(\Omega)}, \]
    where $\alpha' = 1-\frac{n}{p} \in (0,1)$.
    Then the map
    \[ \ell: f\in W^{-1,p}(\Omega) \mapsto -\delta_x \circ \Pi_{(\ker L)^\perp}(u_f) \in \bR \]
    is a bounded linear functional. Hence there exists a unique $G(x,\cdot) \in W_0^{1,q}(\Omega)$ such that
    \[ \int_\Omega G(x,\cdot) f = \ell(f), \]
    where $q$ is the conjugate of $p$; moreover $\|G(x,\cdot)\|_{W^{1,q}(\Omega)}$ is uniformly bounded independent of $x$. Note that when $f=\psi_i$, the corresponding solution $u_f$ is trivial and thus $\int_\Omega G(x,\cdot) \psi_i = 0$ for all $i\in \{1, \cdots, m\}$. 
    Equivalently, for any $u\in W_0^{1,p}(\Omega)$, we have
    \[ f:= Lu \in W^{-1,p}(\Omega) \cap (\ker L)^\perp \rightsquigarrow u_f = u-\sum_{i=1}^m \left(\int_\Omega u\psi_i \right) \psi_i \rightsquigarrow \Pi_{(\ker L)^\perp} (u_f) = \Pi_{(\ker L)^\perp} u. \]
    Hence
    \begin{equation}\label{tmp:Green}
        \int_\Omega G(x,\cdot) Lu = \ell(Lu) = -\delta_x \cdot \Pi_{(\ker L)^\perp} (u),
    \end{equation}
    which proves the Green's identity \eqref{eq:Green}.

    In particular, by taking $\varphi \in C_c^\infty(\Omega\setminus\{x\})$ (by which we mean a smooth function with compact support in the interior of $\Omega\setminus\{x\}$) in \eqref{eq:Green}, it follows by the divergence theorem and \eqref{eq:Green} that
    \begin{equation}\label{tmp:Green-awayfrompole}
        \int_\Omega LG(x,\cdot) \varphi = \int_\Omega G(x,\cdot) L\varphi = -\delta_x \cdot \Pi_{(\ker L)^\perp} (\varphi) = \sum_{i=1}^m \left( \int_\Omega \varphi \psi_i \right) \psi_i(x) = \delta_x \circ \Pi_{\ker L} (\varphi).
    \end{equation}
    Hence $LG(x,\cdot) = \sum_{i=1}^m \psi_i(x) \psi_i$ in the distributional sense in $\Omega \setminus\{x\}$. By the smoothness of $\psi_i$, it follows that $G(x,\cdot)$ is also smooth away from $\{x\}$. More precisely, for any ball $B$ such that $2B \subset\subset \Omega\setminus\{x\}$, we have that 
    \[ \|G(x,\cdot)\|_{W^{2,q}(B)} \lesssim \|G(x, \cdot)\|_{L^q(2B)} + \|\sum_{i=1}^m \psi_i(x) \psi_i\|_{L^{q}(2B)} \leq \|G(x, \cdot)\|_{L^q(\Omega)} + \|\sum_{i=1}^m \psi_i(x) \psi_i\|_{L^{q}(\Omega)}. \]
    By iterating this (and making $B$ smaller), we get a uniform $C^1$ bound on $G(x,\cdot)$, which depends only on the distance to $\partial(\Omega \setminus\{x\})$.

    Let $G^0(x,\cdot) \in W_0^{1,q}(\Omega)$ denote the classical Green's function for $\Omega$ with pole at $x$. Then for any $u\in W_0^{1,p}(\Omega)$,
    \[ \int_{\Omega} G^0(x,\cdot) \Delta u = -u(x). \]
    Combined with \eqref{tmp:Green}, we have that
    \[ \int_{\Omega} (G(x,\cdot) - G^0(x,\cdot) ) \Delta u = \sum_{i=1}^m \left( \int_\Omega u\psi_i \right) \psi_i(x) - \lambda \int_\Omega G(x,\cdot) u. \]    
    Combined with the invertibility of the operator $\Delta$ on $H_0^1(\Omega)$, this implies that $G(x,\cdot) - G^0(x,\cdot)$ coincides with a bounded linear functional defined on $f\in W^{-1,p}(\Omega)$, with a uniform upper bound that depends on the norms of $\psi_i$, $\lambda$, and the uniform upper bound of $\|G(x,\cdot)\|_{L^q(\Omega)}$. Therefore $\|G(x,\cdot) - G^0(x,\cdot)\|_{W^{1,q}(\Omega)}$ is uniformly bounded. We will use this to show that the singular behavior of $G(x,\cdot)$ at the pole is the same as the classical Green's function $G^0(x,\cdot)$.

    The proof of the symmetry $G(x,y)=G(y,x)$ is similar to the analogous proof for the classical Green's function. Let $u := G(x, \cdot)$ and $v:= G(y,\cdot)$. Let $\epsilon$ be arbitrarily small and let $\Omega_\epsilon:= \Omega \setminus (B_\epsilon(x) \cup B_\epsilon(y))$. Then
    \begin{equation}\label{tmp:divergence-polexy}
        \int_{\Omega_\epsilon} \left( uLv - vLu \right) = -\int_{\partial B_\epsilon(x) \sqcup \partial B_\epsilon(y)} \left( u\partial_r v - v \partial_r u \right), 
    \end{equation} 
    where we use the fact that the unit normal derivative of $B_\epsilon(x)$ (resp., $B_\epsilon(y)$) is the radial derivative with vertex at $x$ (resp., $y$). By \eqref{tmp:Green-awayfrompole} and a density argument, we obtain that
    \begin{equation}\label{tmp:interior}
        \int_{\Omega_\epsilon} uLv = \sum_{i=1}^m \left( \int_{\Omega_\epsilon } u\psi_i \right)\psi_i(y) \rightarrow 0,
    \end{equation} 
    as $\epsilon \to 0$, since $u=G(x,\cdot)$ and $ G(x, \cdot) \perp \psi_i$ in $\Omega$. Similarly $\int_{\Omega_\epsilon} vLu \to 0$ as $\epsilon \to 0$. On the other hand, we claim that
    \begin{equation}\label{tmp:polex}
        -\int_{\partial B_\epsilon(x)} \left( u\partial_r v - v \partial_r u \right) \to -c_n v(x)= -c_n G(y,x),
    \end{equation} 
    and
    \begin{equation}\label{tmp:poley}
        -\int_{\partial B_\epsilon(y)} \left( u\partial_r v - v \partial_r u \right) \to  c_n u(y) = c_n G(x,y), 
    \end{equation}
    for some dimensional constant $c_n$. It suffices to prove the estimate near $x$, and the estimate near $y$ is analogous. Recall that away from the pole $y$, $\nabla v = \nabla G(y,\cdot)$ and $\Delta v$ are uniformly bounded, and the function $h: = G(x, \cdot) - G^0(x,\cdot) $ satisfies the equation $\Delta h
  =\sum_{i=1}^m\psi_i(x)\psi_i-\lambda G(x,\cdot) \in L^q_{\mathrm{loc}}$. By the divergence theorem, 
    \begin{align*}
        \int_{\partial B_\epsilon(x)} u\partial_r v & = \int_{\partial B_\epsilon(x)} h\, \partial_r v + \int_{\partial B_\epsilon(x)} G^0(x,\cdot) \partial_r v \\
        & = \int_{B_\epsilon(x)}
  \bigl(\nabla h\cdot\nabla v+h\Delta v\bigr) + \int_{\partial B_\epsilon(x)} G^0(x,\cdot) \partial_r v .
    \end{align*}
    Combined with the uniform bound of $\|h\|_{W^{1, q}}$ and
    \[ \int_{\partial B_\epsilon(x)} G^0(x,\cdot) \approx \int_{\partial B_\epsilon(x)} \frac{1}{|z - x|^{n-2}} d\mathcal{H}^{n-1}(z) \text{ or } \int_{\partial B_\epsilon(x)} \left| \log|z-x| \right| d\mathcal{H}^{1}(z) \to 0, \]
    we conclude that $\int_{\partial B_\epsilon(x)} u\partial_r v \to 0$ as $\epsilon \to 0$. Analogously, by the divergence theorem
    \begin{align*}
        \int_{\partial B_\epsilon(x)} v\partial_r u & = \int_{\partial B_\epsilon(x)} v\partial_r h + \int_{\partial B_\epsilon(x)} v\partial_r G^0(x,\cdot) \\
        & = \int_{B_\epsilon(x)}
  \bigl(\nabla v\cdot\nabla h+v\Delta h\bigr)  + \int_{\partial B_\epsilon(x)} v\partial_r G^0(x,\cdot).
    \end{align*}
    By uniform norm bound on $h$, the continuity of $v$ at $x$, and $\int_{\partial B_\epsilon(x)} \partial_r G^0(x,\cdot) \to -c_n$
    for some dimensional constant $c_n>0$ independent of $x$, we conclude that 
    \[ \int_{\partial B_\epsilon(x)} v\partial_r u \to -c_n v(x) = -c_n G(y,x).\]
    We finish the proof of the claim \eqref{tmp:polex}. Combining \eqref{tmp:divergence-polexy}, \eqref{tmp:interior}, \eqref{tmp:polex} and \eqref{tmp:poley}, we show that the reduced Green's function is symmetric, i.e.\ $G(x,y) = G(y,x)$.
\end{proof}
 
Define
\[
    V:=\mathrm{span}\left\{\frac{\partial\psi_1}{\partial\nu},\dots,\frac{\partial\psi_m}{\partial\nu}\right\}\subset C^{1,\alpha}(\partial \Omega)\subset L^2(\partial \Omega),
\]
and let $\Pi_{V^\perp}:L^2(\partial \Omega)\to V^\perp$ be the $L^2(\partial \Omega)$-orthogonal projection.\footnote{We caution that the functions $\frac{\partial\psi_i}{\partial\nu}$ need not be $L^2(\partial\Omega)$-orthogonal, so the projection $\Pi_{V^\perp}$ does not have the usual nice expression.} By the continuity of elements in $V$, when restricted to continuous functions $g\in C(\partial\Omega)$ we have $\Pi_{V^\perp} g \in C(\partial\Omega)$. 
\begin{lem}[Reduced Poisson kernel]\label{lem:reduced-poisson-kernel}
For any $z\in \partial\Omega$ and $1<q<\frac{n}{n-1}$ there exists a unique function $P_z\in L^q(\Omega)\cap (\ker L)^\perp$ such that\footnote{Clearly any kernel element $u\in \ker L$ annihilates the left hand side, so we need the projection operator $\Pi_{V^\perp}$ on the right hand side to make sure they are compatible.}
\begin{equation}\label{eq:Poisson-kernel}
    \int_{\Omega} P_z Lu = \delta_z \circ \Pi_{V^\perp}\left(\frac{\partial u}{\partial\nu}\right), \qquad \text{ for any } u \in W^{2,p}\cap W^{1,p}_0(\Omega).
\end{equation}
In particular $LP_z=0$ in $\Omega$ in the sense of distributions.

Assume that 
\begin{equation}\label{asmp:vanish-Poisson-kernel}
	\frac{\partial\psi_i}{\partial\nu}(z) = 0, \qquad \text{ for all } i\in \{1, \cdots, m\}.
\end{equation} 
Then
\begin{enumerate}[label=\roman*)]
    \item $P_z = -\frac{\partial G}{\partial \nu}(z,\cdot)$ in $L^q(\Omega)$, where the outer normal derivative $\frac{\partial}{\partial\nu}$ is taken with respect to the first variable,
    \item $P_z = 0$ on $\partial\Omega \setminus\{z\}$ in the classical sense, and $\frac{\partial P_z}{\partial \nu}$ is well-defined away from $z$,
    \item\label{item:Pk-near-pole} for any $w\in C^2(\Omega)$ such that $w|_{\partial \Omega}$ is  compactly supported in $\partial\Omega \setminus\{z\}$, we have
\begin{equation}\label{eq:Poisson-kernel-boundary}
    \int_\Omega P_z Lw = \frac{\partial w}{\partial \nu}(z) - \int_{\partial \Omega }w\,\frac{\partial P_z}{\partial \nu}.
\end{equation}
\end{enumerate}
\end{lem}
\begin{rmk}
    The assumption that $w|_{\partial\Om}$ is supported away from $z$ in \ref{item:Pk-near-pole} can be replaced by a superlinear decay assumption on $w$ near $z$; see the proof of Lemma \ref{lm:unique-triple}.
\end{rmk}
\begin{proof}
Fix $p>n$ and let $q$ be its dual exponent. By the Fredholm alternative of $L$, for each $f\in L^p(\Omega)$ there is a unique solution $u_f\in W^{2,p}\cap W^{1,p}_0(\Omega)\cap (\ker L)^\perp$ of
\[ L u_f=f^\perp,   \]
where $f^{\perp} := f- \sum_{i=1}^m \left( \int_\Omega f \psi_i \right) \psi_i$ is the $L^2$-projection of $f$ onto $(\ker L)^\perp$; moreover,
\[
    \|u_f\|_{W^{2,p}(\Omega)}\lesssim \|f^{\perp}\|_{L^p(\Omega)} \lesssim \|f\|_{L^p(\Omega)}.
\]
Since $p>n$, Sobolev inequality implies $W^{2,p}(\Omega)\hookrightarrow C^{1,\alpha'}(\Omega)$, where $\alpha'=1-\frac np>0$. In particular
\begin{equation}
	\left| \frac{\partial u_f}{\partial\nu}\right|
    \leq \|u_f\|_{C^1(\Omega)} \lesssim \|u_f\|_{W^{2,p}(\Omega)}\lesssim \|f\|_{L^p(\Omega)}.
\end{equation}
Thus
\[ \ell: f \in L^p(\Omega) \mapsto \delta_z \circ \Pi_{V^\perp}\left( \frac{\partial u_f}{\partial \nu} \right) \in \bR  \] 
is a bounded linear functional. By duality, there exists a unique $P_z\in L^q(\Omega)$ such that $\int_{\Omega} P_z f = \ell(f)$.
In particular, for any $f= \psi_i$ the unique solution $u_f \equiv 0$ and thus $\int_\Omega P_z \psi_i = 0$.

Equivalently, for any $u\in W^{2,p}\cap W_0^{1,p}(\Omega)$, we have
\[ f:= Lu \rightsquigarrow u_f = u-\sum_{i=1}^m \left(\int_\Omega u\psi_i \right) \psi_i \rightsquigarrow \Pi_{V^\perp} \left( \frac{\partial u_f}{\partial \nu} \right) = \Pi_{V^\perp} \left( \frac{\partial u}{\partial \nu} \right). \]
Therefore
\[ \int_\Omega P_z Lu = \int_\Omega P_z f = \delta_z \circ \Pi_{V^\perp}\left( \frac{\partial u_f}{\partial \nu} \right) = \delta_z \circ \Pi_{V^\perp}\left( \frac{\partial u}{\partial \nu} \right), \]
which proves the identity in \eqref{eq:Poisson-kernel}.
In particular, applying any $u\in C_c^\infty(\inte \Omega)$ to \eqref{eq:Poisson-kernel} we obtain that $LP_z = 0$ in the sense of distributions.

Under the assumption \eqref{asmp:vanish-Poisson-kernel}, we have $\delta_z \circ \Pi_{V^\perp} = \delta_z$. We claim that 
\begin{equation}\label{claim:Poisson-kernel-Green}
    P_z = -\frac{\partial G}{\partial \nu}(z,\cdot) := \lim_{\tau \to 0+} \frac{G(z-\tau \nu_z,\cdot )}{\tau}, \qquad \text{ in } L^q(\Omega), 
\end{equation} 
where $G(\cdot, \cdot)$ is the reduced Green's function in Lemma \ref{lm:Green}, and $\nu_z$ denotes the unit outer normal vector of $\Omega$ at $z$. Firstly, we note that for fixed $y\in \Omega$, by the symmetry of the reduced Green's function, the function $G(\cdot, y) = G(y,\cdot)$ also satisfies an (inhomogeneous) elliptic PDE with vanishing boundary data in $\Omega \setminus\{y\}$. Hence 
\[ -\frac{\partial G}{\partial \nu}(z,y) := \lim_{\tau \to 0+} \frac{G(z-\tau \nu_z,y)}{\tau} \]
is well-defined pointwise. Moreover, for any $y\in \Omega \setminus B_{2r}(z)$ and $\tau$ sufficiently small,
\begin{equation}\label{tmp:Taylor-Green-boundary}
    \left|\frac{G(z-\tau \nu_z,y)}{\tau} - \left(-\frac{\partial G}{\partial \nu}(z,y) \right) \right| \leq \tau \sup_{x\in \Omega \cap B_r(z) }| \nabla_x^2 G(x,y)| \leq \tau \|G(\cdot, y)\|_{C^2(\Omega \cap B_r(z))}.
\end{equation}
By elliptic boundary estimates and the uniform bound on $\|G(\cdot, y)\|_{L^q(\Omega)}$, we have that $\|G(\cdot, y)\|_{C^2(\Omega \cap B_r(z))}$ is uniformly bounded for $y\in \Omega \setminus B_{2r}(z)$, with a constant depending on $r$. Thus
\begin{equation}\label{eq:CV-bdder-Green}
	\frac{G(z-\tau \nu_z,\cdot )}{\tau} \to -\frac{\partial G}{\partial \nu}(z,\cdot) \qquad \text{ as } \tau \to 0+ \text{ uniformly on compact subsets of } \Omega \setminus \{z\}.   
\end{equation} 

For any $f\in L^p(\Omega) \cap (\ker L)^\perp$, let $u_f \in  W^{2,p} \cap W_0^{1,p}(\Omega) \cap (\ker L)^\perp$ be the unique weak solution of $Lu=f$. Then by \eqref{eq:Green}, we have
\[ \int_\Omega G(z-\tau \nu_z, \cdot) f = -u_f(z-\tau \nu_z), \]
and thus
\begin{equation}\label{tmp:Poisson-kernel-Green1}
	\int_\Omega \frac{G(z-\tau \nu_z, \cdot)}{\tau} f = -\frac{u_f(z-\tau \nu_z)}{\tau} \to \frac{\partial u_f}{\partial \nu}(z), \qquad \text{ as } \tau \to 0+. 
\end{equation} 
Moreover, by the $C^{1,\alpha}$ regularity of $u_f$ up to the boundary, the convergence above is uniform in $\|u_f\|_{W^{2,p}} \lesssim \|f\|_{L^p}$.
On the other hand, by \eqref{eq:Poisson-kernel} and the assumption \eqref{asmp:vanish-Poisson-kernel}, we have
\begin{equation}\label{tmp:Poisson-kernel-Green2}
	\int_\Omega P_z f = \delta_z\left( \frac{\partial u_f}{\partial \nu} \right) = \frac{\partial u_f}{\partial \nu}(z). 
\end{equation} 
Combining \eqref{tmp:Poisson-kernel-Green1} and \eqref{tmp:Poisson-kernel-Green2}, we know that 
\[ \frac{G(z-\tau \nu_z, \cdot)}{\tau} \to P_z, \qquad \text{ in }L^q(\Omega). \]
Combined with \eqref{eq:CV-bdder-Green}, this proves the claim \eqref{claim:Poisson-kernel-Green}.
In particular, by \eqref{eq:CV-bdder-Green} and the fact that $G(z-\tau\nu_z,\cdot) = 0$ on $\partial\Omega$, it follows that $P_z = 0$ on $\partial\Omega \setminus\{z\}$ in the classical sense.

To show that $P_z$ has better regularity (than merely $L^q$) away from $z$, by the claim \eqref{claim:Poisson-kernel-Green} we need to study the mixed derivatives of the reduced Green's function on both variables. To that end, let $\beta$ be the multi-index $(\beta_1, \cdots, \beta_n)\in (\mathbb{N} \cup\{0\})^n$ with $|\beta|:= \beta_1 + \cdots + \beta_n \leq 2$, and denote
\[ \partial_y^\beta G(\cdot, y) = \partial_{y_1}^{\beta_1} \cdots \partial_{y_n}^{\beta_n} G(\cdot, y). \]
As discussed above, for $y\in \Omega \setminus B_{2r}(z)$, we have that $L G(\cdot, y) = \sum_{i=1}^m \psi_i(\cdot) \psi_i(y)$ in $\Omega \cap B_r(z)$, where the differential operator $L$ is applied to the first variable. Applying the $y$-derivative and using the regularity of $G$ on $\Omega\times \Omega$ away from the diagonal, we get that
\[ L \left( \partial_y^{\beta} G(\cdot, y) \right) = \partial_y^{\beta} \left( LG(\cdot, y) \right) = \sum_{i=1}^m \psi_i(\cdot) \partial_y^{\beta} \psi_i(y), \qquad \text{ on } \Omega \cap B_r(z). \]
By boundary elliptic estimates, we can iterate from the uniform upper bound on $\|G(\cdot, y)\|_{W^{1,q}}$, $1<q<\frac{n}{n-1}$ and obtain uniform pointwise bounds on $\partial_y^{\beta} G(\cdot, y)$ (depending on $r$) in $B_{r/2}(z)$; moreover, the $C^{2,\alpha}$ regularity of the domain guarantees that $\partial_y^{\beta} G(\cdot, y)$ also has $C^{2,\alpha}$ regularity, and in particular $\frac{\partial }{\partial\nu} \partial_y^\beta G(\cdot,y)$ is well-defined at the boundary. By an analogous argument as in \eqref{tmp:Taylor-Green-boundary}, this implies 
\[ \partial_y^\beta P_z = -\partial_y^\beta \left(\frac{\partial G}{\partial \nu}(z,\cdot) \right) = -\frac{\partial }{\partial\nu} \partial_y^\beta G(z,\cdot) \]
in $\Omega \setminus B_{2r}(z)$. Thus $P_z \in C^\infty(\inte\Omega\setminus B_{2r}(z))$ (and $C^{1,\alpha}$ up to the boundary $\partial\Omega \setminus B_{2r}(z)$).

Let $w\in C^2(\Omega)$ be such that $w|_{\partial \Omega}$ is supported in a relatively open subset $\Gamma \subset \partial\Omega$ with $z\notin \overline{\Gamma}$. Let $\tilde{w}$ be a $C^2$ extension of $w|_{\partial\Omega}$ supported in a neighborhood $U$ of $\Gamma$ with $z\notin \overline{U}$. Then $w-\tilde{w} \in W^{2,p}\cap W_0^{1,p}(\Omega)$, and hence \eqref{eq:Poisson-kernel} implies
\begin{equation}\label{tmp:Poisson-kernel-bd1}
	\int_\Omega P_z L(w-\tilde{w}) = \frac{\partial (w-\tilde{w})}{\partial \nu}(z) = \frac{\partial w}{\partial \nu}(z). 
\end{equation} 
On the other hand, since $P_z, \tilde{w}$ are both $C^2$ in $\Omega \cap U$, we may apply the classical divergence theorem, $LP_z = 0$, and the boundary data of $P_z, \tilde{w}$ to get
\begin{equation}\label{tmp:Poisson-kernel-bd2}
	\int_\Omega P_z L\tilde{w} =\int_{\Omega \cap U} P_z L\tilde{w} = \int_{\partial \Omega } P_z \frac{\partial \tilde{w}}{\partial \nu} - \int_{\partial\Omega} \tilde{w} \frac{\partial P_z}{\partial \nu} = -\int_{\partial\Omega} w  \frac{\partial P_z}{\partial \nu}. 
\end{equation} 
Summing up \eqref{tmp:Poisson-kernel-bd1} and \eqref{tmp:Poisson-kernel-bd2} gives the identity \eqref{eq:Poisson-kernel-boundary}.
\end{proof}

\section{Generic Properties of Dirichlet Eigenfunctions}\label{sec:generic}

The following result is stated in ~\cite[Theorem 9]{UhlenbeckGeneric} with an outline of the proof. We give the details in this section.

\begin{thm}[\!\!{\cite[Theorem 9]{UhlenbeckGeneric}}]\label{thm:generic}
    Let $k > n$. Let $M$ be a smooth compact connected $n$-manifold with boundary which admits a $C^k$ embedding into $\bR^n$, and let $\mathrm{Emb}_k(M, \bR^n)$ be the Banach manifold of $C^k$ embeddings.

Then, there is a residual set $\mathfrak{X} \subset \mathrm{Emb}_k(M, \bR^n)$ with the following property: If $F \in \mathfrak{X}$ and $u: \Omega = F(M) \to \bR$ is a Dirichlet eigenfunction for the Laplacian on $\Omega$, then
\begin{enumerate}[label=(\alph*)]
    \item\label{item:1d-eigenspace} The eigenspace of $u$ is one-dimensional,
    \item\label{item:interior-regular-value} Zero is a regular value of $u$ on $\inte\Omega$,
    \item\label{item:boundary-regular-value} Zero is a regular value of $\dd_{\nu} u: \dd \Omega \to \bR$, the normal derivative of $u$ on the boundary.
\end{enumerate}
\end{thm}

Note that the properties above are invariant with respect to $C^k$ reparametrizations of the domain $\Omega$, so a residual set among embeddings gives a residual set among the images, namely among $C^k$ domains diffeomorphic to $M$.
We also note that, to prove Theorem \ref{thm:generic}, it is more convenient to work with closures of bounded open domains rather than the open domains themselves: In this setting, the $C^k$ diffeomorphisms of a fixed closed domain form a Banach manifold, and vector fields on the domain naturally induce perturbations of its boundary, among other things. Hence, despite the slight abuse of notation, we use $\Omega$ to denote the \emph{closure} of an open domain throughout this section. Accordingly, $C_c^\infty(\Omega)$ denotes the class of smooth functions with compact support in \emph{the interior of} $\Omega$, and $H_0^1(\Omega)$ denotes the closure of $C_c^\infty(\Omega)$ under the $H^1$ norm.

Let $k\geq 2$. Fix an embedding $F_0 \in \mathrm{Emb}_{k}(M, \bR^n)$ with image $\Omega := F_0(M) \subset \bR^n$.  
We define 
\begin{align*}
    B_{k} &= B_{k}(\Omega) := \{h \in C^{k}(\Omega, \bR^n)\,|\,h: \Omega \to h(\Omega) \text{ is a } C^{k} \text{ diffeomorphism}\}\\
    &\cong \mathrm{Emb}_{k}(M, \bR^n).
\end{align*}
The identification in the last line follows since $\Omega$ can be identified with $M$ via a $C^{k}$ embedding. 
In particular, each $\Omega_h:= h(\Omega)$ is a bounded connected \emph{closed} domain with $C^k$ boundary.
With this definition, $B_{k}$ is an open subset of the separable Banach space $C^{k}(\Omega, \bR^n)$, so $B_{k}$ is a smooth Banach manifold modeled on $C^k(\Omega,\bR^n)$. If $h \in B_k$, then 
\begin{align}\label{equation:tangent-space-B}
T_h B_k \cong C^{k}(\Omega, \bR^n).
\end{align}

When we work with solutions on a class of domains, there is a subtle difference between taking the Eulerian point of view and the Lagrangian point of view. We need to work with Banach manifolds and function spaces on a fixed domain $\Omega$, in order to apply the transversality theorems; but the statements we are interested in and the proof arguments are about solutions to PDEs on the variable domains $\Omega_h$. To that end, we consider the operator $L_h$, which is the pullback of the Dirichlet Laplacian on $\Omega_h$ to the fixed domain $\Omega$. That is,
\begin{equation}\label{def:Lh}
    L_h u := h^{\star} \Delta_{\Omega_h} (h^{\star})^{-1} u = \Delta_{\Omega_h}(u \circ h^{-1}) \circ h, \qquad \text{ for any }u \in H^2(\Omega).
\end{equation}
After considering the composition with $h$ we have that
\[ L_h: H^2 \cap H_0^1(\Omega) \to L^2(\Omega) \text{ is an isomorphism}. \]
It is well known that in a $C^k$ domain $\Omega_h$, any eigenfunction in $H^2 \cap H_0^1(\Omega_h)$ automatically belongs to $H^k(\Omega_h)$ (see \cite[Theorems 8.9 and 8.13]{GT}). It also automatically belongs to $C^{k-1,\beta}(\Omega_h)$ for any $\beta \in (0,1)$ by Schauder theory (see \cite[Theorem 6.6 and Problem 6.2]{GT}).
So, even though we choose $H^2(\Omega)$ as the domain of $L_h$ for simplicity, the functions involved often have improved regularity thanks to the elliptic regularity theory. Lastly, we remark that the operator $L_h$ is invariant with respect to $C^k$ reparametrizations of the domain $\Omega_h$.

Now, define
\begin{equation}\label{eq:S}
    S  = S(\Omega) :=\left\{u \in H^2 \cap H_0^1(\Omega) \,\Big|\,\int_\Omega u^2 = 1\right\}. 
\end{equation}
This is the function space we will work with; or rather, we work with functions $u\circ h^{-1} \in H^2 \cap H_0^1(\Omega_h)$, for $u\in S(\Omega)$. We remark that this is where our approach departs from the sketch Uhlenbeck outlined. In our setup the variations of eigenfunctions always vanish on the boundary, and the boundary variations are entirely encoded in the domain variations. 
Then we use PDE arguments to verify the assumptions about regular values (or equivalently, the necessary density properties, which Uhlenbeck described as ``almost impossible to verify''~\cite[Example 3]{UhlenbeckGeneric}), in order to apply the transversality theorem. Our approach is greatly inspired by the work of Henry in \cite{HenryPerturbationbook05}.

We will usually suppress the dependence on $\Om$ from our notation for $B_k$ and $S$.
Regarding $S$ as a smooth codimension-one Banach submanifold of $H^2 \cap H_0^1(\Omega)$, its tangent space $T_u S$ at $u \in S$ is given by
\begin{align}\label{eqn:TuS}
T_u S = \left\{v \in H^2 \cap H_0^1(\Omega) \,\Big|\,\int_{\Omega} uv = 0\right\}.
\end{align}
Finally, define $\varphi: S \times \bR \times B_k \to L^2(\Omega)$ by
\begin{align}\label{eqn:definition-varphi}
    \varphi(u,\lambda, h) := (L_h + \lambda) u.
\end{align}
Note that $\varphi$ depends linearly on $u$, when extended to the Banach space $H^2 \cap H_0^1(\Omega)$, and $\lambda \in \bR$.
By definition \eqref{def:Lh}, $L_h$ is a linear elliptic operator with second-order coefficients given by the product of $\nabla (h^{-1})(h(x))$, and first-order coefficients given by $\nabla^2 (h^{-1})(h(x))$. Hence the map $h\in B_k \mapsto L_h \in \mathcal{L}(H^2 \cap H_0^1(\Omega), L^2(\Omega))$ is smooth. Therefore, $\varphi$ is a smooth map of Banach manifolds.

\subsection{Generic simplicity of eigenvalues}\label{section:generic-simplicity}

In this subsection, we prove \ref{item:1d-eigenspace} from Theorem \ref{thm:generic}. 
This result has been detailed in \cite{Micheletti-1}, \cite{Micheletti-2} and ~\cite[Example 6.1]{HenryPerturbationbook05}, but we provide our own treatment for completeness since the setup is used in the proofs of \ref{item:interior-regular-value} and~\ref{item:boundary-regular-value} from Theorem \ref{thm:generic}.

Using the Sard-Smale or transversality theorem, Uhlenbeck proved a general theorem relating regular values of certain maps on Banach manifolds to the generic properties of eigenfunctions for associated operators~\cite[Theorem 1]{UhlenbeckGeneric}. We formulate this theorem below, applied to our setting.

\begin{thm}[\!{\cite[Theorem 1]{UhlenbeckGeneric}}]\label{thm:Uhlenbeck-regular-value-phi}
    If $k\geq 2$ and if zero is a regular value of $\varphi: S \times \bR \times B_k \to L^2(\Om)$, then the set 
    $$\{h \in B_k \,|\,L_h \text{ has one-dimensional eigenspaces}\}$$ 
is residual in $B_k$. 

    In particular, if zero is a regular value of $\varphi$, then there is a residual set $\mathfrak{X} \subset \mathrm{Emb}_k(M, \bR^n)$ such that for each $F \in \mathfrak{X}$, every Dirichlet eigenfunction for the Laplacian on $F(M)$ has a one-dimensional eigenspace, i.e.\ \ref{item:1d-eigenspace} from Theorem \ref{thm:generic} holds.
\end{thm}
\begin{proof}
The first statement is a reformulation of \cite[Theorem 1]{UhlenbeckGeneric}, which follows by applying the Thom‑ Abraham transversality theorem \cite[Transversality Theorem 1]{UhlenbeckGeneric} to the map $\varphi$, together with the observation \cite[Lemma 2.3]{UhlenbeckGeneric}. For the second statement, let 
\[ \mathfrak{X} = \{h \circ F_0 \,|\,h \in B_k \text{ and }L_h \text{ has one-dimensional eigenspaces}\}. \] If  $h \circ F_0 \in \mathfrak{X}$ and $\phi$ is a Dirichlet eigenfunction for $h \circ F_0(M)$, then $\phi \circ h$ is a Dirichlet eigenfunction of $L_h$ and thus spans its eigenspace. This implies that $\phi$ spans its eigenspace for the Dirichlet Laplacian on $h \circ F_0(M)$. The second statement then follows from the first.
\end{proof}

Our goal now is to show that $0$ is a regular value of $\varphi: S \times \bR \times B_k \to L^2(\Om)$, and Theorem \ref{thm:Uhlenbeck-regular-value-phi} will allow us to conclude ~\ref{item:1d-eigenspace} of Theorem \ref{thm:generic}. 

If $E$ is a closed domain in $\bR^n$ with $C^k$ boundary and $E$ is $C^k$-diffeomorphic to $M$, we define 
$$\widetilde{\varphi}_E: H^2 \cap H_0^1(E) \times \bR \times B_k(E) \to L^2(E),$$
by
\begin{equation}\label{eq:defvarp-anydomain}
    \widetilde{\varphi}_E(u,\lambda, h) := \Delta_{E_h}(u \circ h^{-1}) \circ h + \lambda u,
\end{equation}
where $E_h = h(E)$. Our definition of $\widetilde{\varphi}_E$ is a generalization of $\varphi$ which allows the domain on which $\varphi$ is defined to vary and also extends $\varphi$ from $S$ to the unnormalized space $H^2 \cap H_0^1$. 

\noindent Let 
\begin{align*}
    (u,\la, h) &\in H^2 \cap H_0^1(E) \times \bR \times B_k(E) \text{ be arbitrary},\\
    (v, \eta, s) &\in T_{(u, \la, h)}\big(H^2 \cap H_0^1(E) \times \bR \times B_k(E)\big) \cong H^2 \cap H_0^1(E) \times \bR \times C^k(E, \bR^n),
\end{align*}
where we use ~\eqref{equation:tangent-space-B} and think of $s \in T_h B_k(E)$ as a $C^k$ vector field on $E$.

Consider the differential of $\widetilde{\varphi}_E$ at $(u,\la, h)$, denoted by
$$d\tilde{\varphi}_E|_{(u,\la, h)}: H^2 \cap H_0^1(E) \times \bR \times C^k(E, \bR^n) \to T_{\tilde{\varphi}_E(u,\la, h)}L^2(E) \cong L^2(E).$$
To compute $d\widetilde{\varphi}_E$ at $h = \Id\!|_{E}$, we consider a $C^1$ path of embeddings $h(t,\cdot) \in B_k(E)$, a $C^1$ path of functions $u(t,\cdot) \in H^2 \cap H_0^1(E)$ and a $C^1$ path of real numbers $\lambda_t$, and take the $t$-derivative of \eqref{eq:defvarp-anydomain}. We obtain
\begin{align}\label{equation:dphi-at-identity}
d\widetilde{\varphi}_E|_{(u, \lambda, \Id\!|_{E})}(v, \eta, s) = (\Delta_{E} + \la) v + \eta u + s \cdot \nabla\left((\Delta_{E} + \la) u\right) - (\Delta_{E} + \la)(s \cdot \nabla u). 
\end{align}
See also \cite[Theorem 2.2 and remark]{HenryPerturbationbook05} and \cite[Theorem 4]{PereiraHenry22}. Even though it appears that the last two terms contain $\nabla^3 u \in H^{-1}$, they actually cancel each other out as a commutator $[s\cdot \nabla, \Delta_E + \lambda]$ and take values in $L^2$.

Let $h_0 = h \in B_k = B_k(\Om)$, and let $h_t$ be a $C^1$ path in $B_k$ with $\frac{d}{dt}h_t|_{t=0} = s \in T_h B_k$. Define the $C^1$ path $\tilde{h}_t := h_t \circ h^{-1}$ in $B_k(\Om_h)$. Note that $\tilde{h}_0 = \mathrm{Id}|_{\Omega_h}$. Likewise, let $u_0 = u \in S$, and let $u_t$ be a $C^1$ path with $\frac{d}{dt} u_t |_{t=0} = v \in T_u S$. Define the $C^1$ path $\tilde{u}_t := u_t \circ h^{-1} \in H^2 \cap H_0^1(\Omega_h)$. Also, let $\la_t$ be a $C^1$ path in $\bR$ with $\la_0 = \la$ and $\frac{d}{dt}\la_t|_{t=0} = \eta$. We then have that $\frac{d}{dt}\tilde{h}_t|_{t=0} = s \circ h^{-1}$ and $\frac{d}{dt} \tilde{u}_t|_{t=0} = v \circ h^{-1}$. Due to the $L^2$ normalization in the definition of $S$, it may be that $\tilde{u}_t \notin S(\Omega_h)$. This is why we extended $\varphi$ to  $\widetilde{\varphi}_E$, in order to accommodate $\tilde{u}_t$. 

Now, 
\begin{align}\label{equation:change-of-origin}
\widetilde{\varphi}_{\Omega_h}(\tilde{u}_t, \la_t, \tilde{h}_t) \circ h = \Delta_{\tilde{h}_t(\Omega_h)}(\tilde{u}_t \circ \tilde{h}_t^{-1}) \circ \tilde{h}_t \circ h + \la_t \tilde{u}_t \circ h &= \Delta_{h_t(\Omega)}(u_t\circ h_t^{-1})\circ h_t + \la_t u_t\nonumber\\
&= \varphi(u_t, \la_t, h_t).
\end{align}
Since $\varphi: S \times \bR \times B_k \to L^2(\Om)$, we have that 
\begin{align}
d\varphi|_{(u,\la, h)}: T_u S \times \bR \times C^k(\Om, \bR^n) \to T_{\varphi(u,\la, h)}L^2(\Om) \cong L^2(\Om).\end{align}
Differentiating ~\eqref{equation:change-of-origin} at $t=0$ and applying~\eqref{equation:dphi-at-identity}, we find that
\begin{align}
    d\varphi|_{(u, \la, h)}(v, \eta, s) &= d\widetilde{\varphi}_{\Om_h}|_{(u \circ h^{-1}, \la, \mathrm{Id}|_{\Om_h})}(v \circ h^{-1}, \eta, s \circ h^{-1}) \circ h\nonumber\\
    &=\left((\Delta_{\Om_h} + \la) (v\circ h^{-1}) + \eta \,u \circ h^{-1}  + (s\circ h^{-1}) \cdot \nabla\big((\Delta_{\Om_h} + \la) (u\circ h^{-1})\big) \right.\nonumber\\
    &\hspace{2.5in} \left. - (\Delta_{\Om_h} + \la)((s \circ h^{-1}) \cdot \nabla (u \circ h^{-1}))\right) \circ h.
    \nonumber\\
    &= (L_h + \la) v + \eta u + \Big((s\circ h^{-1}) \cdot \nabla \big(((L_h + \la) u)\circ h^{-1}\big)\Big) \circ h  \nonumber\\
    &\hspace{2.5in} - (\Delta_{\Om_h} + \la)((s \circ h^{-1}) \cdot \nabla (u \circ h^{-1})) \circ h.\nonumber
\end{align}
From now on, we consider $(u, \la, h) \in S \times \bR \times B_k$ where $\varphi(u,\la, h) = 0$. Then, we have that 
\begin{align}\label{eqn:differential-phi-0}
    d\varphi|_{(u, \la, h)}(v, \eta, s) = (L_h + \la) v + \eta u - (\Delta_{\Om_h} + \la)((s \circ h^{-1}) \cdot \nabla (u \circ h^{-1})) \circ h.
\end{align}
In particular, if $\eta = 0$ and $s=0$, then $d\varphi|_{(u, \la, h)}(v, 0, 0) = (L_h + \la) v$. 

\begin{lem}\label{lemma:appropriate-surjectivity}
Let $k\geq 2$, and let $(u, \la, h) \in S \times \bR \times B_k$ such that $\varphi(u, \la, h)=0$. Then, for each $g \in L^2(\Om)$, there exists $(v, 0, s) \in T_u S \times \bR \times C^k(\Om, \bR^n)$ such that 
\begin{align}
    d\varphi|_{(u, \la, h)}(v, 0, s) = g.        
\end{align}
In other words, zero is a regular value of $\varphi: S \times \bR \times B_k \to L^2(\Om)$.
\end{lem}
\begin{proof}
Let $V \subset L^2(\Om)$ be the subspace defined by
\begin{align*}
    V &:= d\varphi|_{(u, \la, h)}(T_u S \times \{0\} \times C^k(\Om, \bR^n))\\
    &= \{g \in L^2(\Omega) \,\big|\,d\varphi|_{(u, \la, h)}(v, 0, s) = g \text{ for some }(v, s) \in T_u S \times C^k(\Om, \bR^n)\}
\end{align*}
Suppose that the conclusion of the lemma does not hold. Then, $V$ is a proper subspace of $L^2(\Om)$.
Since $L_h+ \la$ is a Fredholm operator, we find that $\mathrm{im}(d\varphi|_{(u, \la, h)}(\cdot, 0, 0))\subset L^2(\Om)$ has finite codimension and is closed. So, $V\supset \mathrm{im}(d\varphi|_{(u, \la, h)}(\cdot, 0, 0))$ is closed.

Now we switch gears to the domain $\Omega_h$ in order to use the PDE.
For each $h \in B_k$, the map $F_h: L^2(\Om) \to L^2(\Om_h)$, defined by $F_h(u) = u \circ h^{-1}$, is a continuous invertible linear operator. Then, $F_h(V) \subset L^2(\Om_h)$ is a closed, proper subspace. Hence there exists a non-trivial $\widetilde{\psi} \in L^2(\Om_h)$ such that 
\begin{align}\label{eqn:psi-tilde-orthogonal-0}
\widetilde{\psi} \perp F_h(V).\end{align}

For $v \in T_u S$, we have $(L_h + \la)v \in \mathrm{im}(d\varphi|_{(u, \la, h)}(\cdot, 0, 0))\subset V$, and so $(\Delta_{\Om_h}+\la)(v \circ h^{-1}) \in F_h(V)$. Thus, ~\eqref{eqn:psi-tilde-orthogonal-0} implies that for each $v \in T_u S$,
\begin{align}\label{eqn:psi-tilde-orthogonal-1}
\int_{\Om_h} \widetilde{\psi}(\Delta_{\Om_h} +\la)(v \circ h^{-1})  = 0. 
\end{align}
Note that $\varphi(u, \la, h)=0$ means that $u \circ h^{-1}$ is a Dirichlet eigenfunction on $\Om_h$ with eigenvalue $\la$. 
So, $(\Delta_{\Om_h} + \la)(u \circ h^{-1})=0$ and~\eqref{eqn:psi-tilde-orthogonal-1} holds for $v=u$ as well.
Together with the definition of $T_u S$ in \eqref{eqn:TuS}, this implies \eqref{eqn:psi-tilde-orthogonal-1} holds for all $v \in H^2 \cap H_0^1(\Om)$. 
Since $F_h: H^2 \cap H_0^1(\Om)\subset L^2(\Om) \to H^2 \cap H_0^1(\Om_h) \subset L^2(\Om_h)$ is a surjective map, it follows that for each $w \in H^2 \cap H_0^1(\Om_h)$,
\begin{equation}\label{equation:}
\int_{\Om_h} \widetilde{\psi}(\Delta_{\Om_h} +\la)w  = 0. 
\end{equation}
This says that $\widetilde{\psi} \in D((\De_{\Om_h} + \la)^\star)$ and $(\De_{\Om_h} + \la)^\star \widetilde{\psi} = 0$. By self-adjointness of the operator $\De_{\Om_h} + \la$ on $H^2 \cap H_0^1$, $\widetilde{\psi} \in D(\De_{\Om_h} +\la) = H^2 \cap H_0^1(\Om_h)$ and 
\begin{align}\label{eqn:psi-tilde-eigenfunction}
(\De_{\Om_h} + \la)\widetilde{\psi} = 0. \tag{$\dagger$}
\end{align}
By assumption, $\widetilde{\psi}$ is non-trivial, so $\widetilde{\psi}$ does not vanish identically in $\Om_h$. 

On the other hand for $s \in C^k(\Om, \bR^n)$, we have
\begin{align}\label{eqn:orthog-scirch}
    -(\Delta_{\Om_h} + \la)\left((s \circ h^{-1}) \cdot \nabla (u \circ h^{-1})\right) &= d\varphi|_{(u, \la, h)}(0, 0, s) \circ h^{-1} \in F_h(V
    ),
\end{align}
Combining ~\eqref{eqn:psi-tilde-orthogonal-0} with ~\eqref{eqn:orthog-scirch}, we have that for each $s \in C^k(\Om, \bR^n)$,
\begin{align}
    0 &= \int_{\Om_h} \widetilde{\psi} (\Delta_{\Om_h} + \la)\left((s \circ h^{-1}) \cdot \nabla (u \circ h^{-1})\right),\nonumber\\
\intertext{Applying ~\eqref{eqn:psi-tilde-eigenfunction},}
    &= \int_{\Om_h} \widetilde{\psi} (\Delta_{\Om_h} + \la)\left((s \circ h^{-1}) \cdot \nabla (u \circ h^{-1})\right) - \left((s \circ h^{-1}) \cdot \nabla (u \circ h^{-1})\right) (\Delta_{\Om_h} + \la)\widetilde{\psi} \nonumber\\
    &= \int_{\Om_h} \widetilde{\psi} \Delta_{\Om_h}\left((s \circ h^{-1}) \cdot \nabla (u \circ h^{-1})\right) - \left((s \circ h^{-1}) \cdot \nabla (u \circ h^{-1})\right)\Delta_{\Om_h}\widetilde{\psi}\nonumber\\
    &= \int_{\dd \Om_h} \widetilde{\psi}\frac{\partial}{\partial \nu} \left((s \circ h^{-1}) \cdot \nabla (u \circ h^{-1})\right) - \left((s \circ h^{-1}) \cdot \nabla (u \circ h^{-1})\right) \frac{\partial \widetilde{\psi}}{\partial \nu},\nonumber\\
\intertext{Using that $\widetilde{\psi}=0$ on $\dd \Om_h$,}
    &= -\int_{\dd \Om_h} \left((s \circ h^{-1}) \cdot \nabla (u \circ h^{-1})\right) \frac{\partial \widetilde{\psi}}{\partial \nu},\nonumber\\
\intertext{Since $u\circ h^{-1} = 0$ on $\dd \Om_h$, $\nabla (u \circ h^{-1}) = \frac{\partial}{\partial \nu}(u \circ h^{-1}) \nu$.}\label{eqn:integral-normal-derivs}
    &= -\int_{\dd \Om_h} (s \circ h^{-1}) \cdot \nu \frac{\partial(u \circ h^{-1})}{\partial \nu} \frac{\partial \widetilde{\psi}}{\partial \nu}.
\end{align}
Note that we have used the assumption $k\geq 2$ to justify that the eigenfunction $u\circ h^{-1} \in H^2(\Omega_h)$, and thus $(s \circ h^{-1}) \cdot \nabla (u \circ h^{-1}) \in H^1(\Omega_h)$, in order to apply the divergence theorem.

Note that $(s \circ h^{-1}) \cdot \nu$ corresponds to a normal variation of the domain $\Omega$ applied to $\Omega_h$ after a change of variables. We now regard $L^2(\dd\Om_h)$ as the closure of the continuous real-valued functions on $\dd\Om_h$ with respect to $L^2(\cH^{n-1}\lfloor \dd\Om_h)$. We note that $\{(s \circ h^{-1}) \cdot \nu \in L^2(\dd\Om_h)\,|\, s \in C^k(\Om, \bR^n)\}$ is dense in $L^2(\dd\Om_h)$. Since ~\eqref{eqn:integral-normal-derivs} holds for all $ s\in C^k(\Om, \bR^n)$, it implies that $\frac{\partial (u \circ h^{-1})}{\partial \nu} \frac{\partial \widetilde{\psi}}{\partial \nu}$ is orthogonal to a dense set in $L^2(\dd\Om_h)$. The only vector orthogonal to a dense subspace is zero, so
\begin{align}\label{eqn:one-step-to-contradiction1}
    \frac{\partial (u \circ h^{-1})}{\partial \nu} \frac{\partial \widetilde{\psi}}{\partial \nu} = 0, \qquad \text{ on } \partial\Omega_h.\tag{$\dagger\dagger$}
\end{align}
 
Since $u \circ h^{-1}$ is a non-trivial Dirichlet eigenfunction on the $C^k$ domain $\Om_h$, by uniqueness of the Cauchy problem (see for example ~\cite[pg.\ 83]{HenryPerturbationbook05}) $\frac{\partial (u \circ h^{-1})}{\partial \nu}$ cannot vanish almost everywhere on $\partial\Omega_h$; in fact by the unique continuation property, it can at most vanish on a codimension-one set in $\partial\Omega_h$.
Since $\widetilde{\psi}$ is also a Dirichlet eigenfunction by \eqref{eqn:psi-tilde-eigenfunction}, the same property holds for $\frac{\partial \widetilde{\psi}}{\partial \nu}$. Hence \eqref{eqn:one-step-to-contradiction1} is impossible.

Thus, $V = L^2(\Om)$, and this concludes the lemma. 
\end{proof}

Since zero is a regular value of $\varphi$, Theorem \ref{thm:Uhlenbeck-regular-value-phi} applies. We conclude that ~\ref{item:1d-eigenspace} from Theorem \ref{thm:generic} holds.

\subsection{Zero is a regular value in the interior}\label{sec:interior-reg}

In this subsection, we prove \ref{item:interior-regular-value} from Theorem \ref{thm:generic}. 

Consider $\varphi: S \times \bR \times B_k \to L^2(\Om)$ as in the previous subsection. We define
\begin{align}
    Q_k := \varphi^{-1}(0) = \{(u, \lambda, h) \in S \times \bR \times B_k \,|\,(L_h + \lambda)u = 0\},
\end{align}
to incorporate all (normalized) eigenfunctions and eigenvalues on $\Omega_h$.
In Section \ref{section:generic-simplicity}, we showed that $0$ is a regular value of $\varphi$. Hence as discussed in ~\cite[Lemma 2.5]{UhlenbeckGeneric}, $Q_k$ is a smooth Banach submanifold of $S \times \bR \times B_k$; and moreover by ~\eqref{eqn:differential-phi-0}, the tangent space of $Q_k$ at $(u, \la, h) \in Q_k$ is given by
\begin{align}
    T_{(u, \la, h)} Q_k &= \ker d\varphi|_{(u, \la, h)}\nonumber\\\label{eqn:TQ}
    &= \{(v, \eta, s) \in T_u S \times \bR \times C^k(\Omega, \bR^n) \,|\, (L_h + \la) v + \eta u\\
    &\hspace{2.75in}  - (\Delta_{\Om_h} + \la)((s \circ h^{-1}) \cdot \nabla (u \circ h^{-1})) \circ h=0\}. \nonumber
\end{align}
We remark that $Q_k$ is separable as it is a subset of the separable Banach manifold $S \times \bR \times B_k$.

Furthermore, if $(v,\eta,s) \in T_{(u,\lambda,h)} Q_k$ for some $\eta\in \bR$ and $s\in C^k(\Omega,\bR^n)$, then $v$ has improved regularity, namely $v \in H^k(\Omega) \cap C^{k-1,\beta}(\Omega)$ for any $\beta\in (0,1)$. This is because $v\circ h^{-1} \in H^2 \cap H_0^1(\Omega_h)$ satisfies the inhomogeneous equation
\begin{equation}\label{eq:v-inhom}
    (\Delta_{\Omega_h} + \lambda) (v\circ h^{-1}) = -\eta \, u\circ h^{-1} + (\Delta_{\Om_h} + \la)((s \circ h^{-1}) \cdot \nabla (u \circ h^{-1})). 
\end{equation} 
Since $u\circ h^{-1}$ is a Dirichlet eigenfunction in the $C^k$ domain $\Omega_h$, we know $u\circ h^{-1} \in H^k(\Omega_h)$ and $C^{k-1,\beta}$ for any $\beta\in (0,1)$. This implies that the right hand side of the equation \eqref{eq:v-inhom} lies in $H^{k-2}(\Omega_h)$ and $C^{k-3,\beta}(\Omega_h)$. Here we also use the condition $(\Delta_{\Omega_h} + \lambda) (u\circ h^{-1})=0$ to get rid of the third-order term in $u$ (see the remark below \eqref{equation:dphi-at-identity}). The improved regularity of $v$ then follows by \cite[Theorem 8.13 and Problem 6.2]{GT}.

We now define a map $\alpha: Q_k \times \inte\Om \to \bR$ by
\begin{equation}\label{def:alpha}
	\alpha(u, \la, h, z) := u(z).
\end{equation}
The eigenfunction $u\circ h^{-1}$ is smooth in the interior of $\Omega_h$, so $u = (u\circ h^{-1}) \circ h$ has $C^k$ regularity in the interior of $\Omega$. In other words, $\alpha$ is $C^k$ regular in the last variable. To show $\alpha$ is also a regular map on the triplet $(u,\lambda, h) \in Q_k$, where $u$ is only assumed to be in the weak space $H^2 \cap H_0^1(\Omega)$, we want to express the point evaluation $u(z)$ in an integral form. Near the point $h(z) \in \inte \Omega_h$, by lifting the eigenfunction $u\circ h^{-1}$ to a harmonic function in a higher-dimensional space, we may use the mean-value property for harmonic functions to write
\[ u\circ h^{-1}(h(z))  = \int_{\Omega_h } \left( \frac{1}{2\pi} \int_0^{2\pi} e^{\sqrt{\lambda} \sqrt{r^2 - |y-h(z)|^2} \cos \theta } \, d\theta \right) u\circ h^{-1}(y) \, \chi\left( \frac{|y-h(z)|}{r} \right) \, dy. \]
Here $\chi: \mathbb{R}_+ \to [0,1]$ is a smooth cut-off function, such that it vanishes on $[1,+\infty)$ and $\int_{\mathbb{R}_+} \chi(\rho) \rho^{n-1} d\rho = \omega_n$, the volume of the unit ball in $\mathbb{R}^n$.
By the change of variables, it follows that
\begin{align*}
	u(z) & = \int_{\Omega} \left( \frac{1}{2\pi} \int_0^{2\pi} e^{\sqrt{\lambda} \sqrt{r^2 - |h(x)-h(z)|^2} \cos \theta } \, d\theta \right) u(x) \left| \det Dh(x) \right| \, \chi\left(\frac{|h(x)-h(z)|}{r} \right) \, dx \\
	& =: \mathcal{A}(u,\lambda, h, z).
\end{align*}
Let $(u_0, \lambda_0, h_0)\in Q_k$ and $z_0 \in \inte \Omega$ be fixed, and let $r>0$ be a fixed value satisfying $2r< \dist(h_0(z_0), \partial\Omega_{h_0})$. Nearby, the map $\mathcal{A}$ defined above is a smooth map with respect to $(u,\lambda, h) \in L^2(\Omega) \times \R \times B_k$, thus for eigenfunctions we have $\alpha = \mathcal{A}|_{Q_k \times \inte\Omega}$ is $C^k$ regular.

The next theorem is a reformulation of \cite[Theorem 2]{UhlenbeckGeneric}. Note that the hypothesis that zero is a regular value of $\varphi$ already follows from the argument of Section \ref{section:generic-simplicity}, as long as $k\geq 2$. 

\begin{thm}[\!{\cite[Theorem 2]{UhlenbeckGeneric}}]\label{thm:Uhlenbeck-regular-value-alpha} 
If $k \geq n$, and zero is a regular value of both $\varphi: S \times \bR \times B_k \to L^2(\Om)$ and $\alpha: Q_k \times \inte\Om \to \bR$, then the set 
    $$\{h \in B_k \,|\,\text{The eigenfunctions of } L_h \text{ have }0 \text{ as a regular value in} \inte \Om\}$$ 
is residual in $B_k$. 

In particular, if zero is a regular value of both $\varphi$ and $\al$, then there is a residual set $\mathfrak{X} \subset \mathrm{Emb}_k(M, \bR^n)$ such that for each $F \in \mathfrak{X}$, every Dirichlet eigenfunction for the Laplacian on $F(M)$ has zero as a regular value on $\inte F(M)$, i.e.\ \ref{item:interior-regular-value} from Theorem \ref{thm:generic} holds.
\end{thm}
\begin{proof}
Applying \cite[Transversality Theorem 2]{UhlenbeckGeneric} to the map $\alpha$, with finite-dimensional manifolds $X=\inte\Omega, Y= \bR, Y'=\{0\}$, we need the regularity of $\alpha$ to satisfy
\[ k > \max\{0, n-1\} = n-1. \]

For the second statement, let 
$$\mathfrak{X} = \{h \circ F_0\,|\,h \in B_k \text{ and the eigenfunctions of }L_h \text{ have }0 \text{ as a regular value in}\inte\Om\}.$$
If  $h \circ F_0 \in \mathfrak{X}$ and $\phi$ is a Dirichlet eigenfunction for $h \circ F_0(M)$, then $\phi \circ h$ is an eigenfunction of $L_h$ and thus has $0$ as a regular value in $\inte \Om$. Let $p \in \inte (h \circ F_0(M))$ such that $\phi(p) = 0$. Letting $q  = h^{-1}(p) \in \inte F_0(M)$, $\phi\circ h(q) = 0$ and so $d(\phi \circ h)|_q$ is surjective. Since $d(\phi \circ h)|_q = d\phi|_p\circ Dh|_q$ and $Dh|_q$ is a linear isomorphism, we conclude that $d\phi|_p$ is surjective. This implies that $\phi$ has $0$ as a regular value in $\inte (h \circ F_0(M))$. The first statement implies that $\mathfrak{X}$ is residual, so the second statement follows.
\end{proof}

Our goal now is to show that zero is a regular value of $\al$, as long as $k\geq 3$ (so that we have the existence and properties of the reduced Green's function as in Lemma \ref{lm:Green}). Afterwards, for any $k\geq \max\{n,3\}$, we can apply Theorem \ref{thm:Uhlenbeck-regular-value-alpha} to conclude ~\ref{item:interior-regular-value} from Theorem \ref{thm:generic}. 

First, we will show that elements of $T_{(u, \la, h)}Q_k$ are abundant.

\begin{lem}\label{lemma:solvability-for-TQ}
Let $k \geq 2$, let $\Om$ be the image of a $C^k$ embedding of $M$, and let $(u, \la, h) \in S \times \bR \times B_k$ such that $\varphi(u, \la, h)=0$. Let $\{\psi_i\}_{i=1}^m$ be an $L^2(\Om_h)$-orthonormal basis of Dirichlet eigenfunctions with eigenvalue $\lambda$, i.e.\ $\psi_i \in \ker(\De_{\Om_h}+\la) \cap \left(H^2 \cap H_0^1 \right) (\Om_h)$, and let $W \subset L^2(\dd \Om_h)$ be the finite-dimensional subspace defined by
\begin{align}\label{eqn:W-subspace-definition}
    W &:= \mathrm{span}\left\{\frac{\dd(u \circ h^{-1})}{\dd \nu} \frac{\dd \psi_1}{\dd \nu}, \cdots, \frac{\dd(u \circ h^{-1})}{\dd \nu} \frac{\dd \psi_m}{\dd \nu}\right\}.
\end{align}
Then, for any $\rho \in W^\perp$ and $s \in C^k(\Om, \bR^n)$ with boundary data $(s \circ h^{-1})\cdot \nu = \rho$ on $\dd\Om_h$, there exists $v \in H^2 \cap H_0^1(\Om)$ such that $(v, 0, s) \in T_{(u, \la, h)}Q_k$. 

In particular, there is a dense set $Z \subset W^\perp$ such that for each $\rho \in Z$, there are $v \in H^2 \cap H_0^1(\Om)$ and $s \in C^k(\Om, \bR^n)$ such that $(v, 0, s) \in T_{(u, \la, h)}Q_k$ and $(s \circ h^{-1})\cdot \nu = \rho$ on $\dd\Om_h$.

\end{lem}
\begin{rmk}\label{rmk:simple}
    Loosely speaking, what we prove is the following: A dense set of $L^2$ boundary data $\rho$ on $\partial\Om_h$ that is orthogonal to $W$ can be realized as a domain variation $s$ of $\Om_h$; moreover, after the domain variation, there is a corresponding variation $v$ of the Laplacian eigenfunction $u$ with the same eigenvalue.

    In fact, if $\lambda$ is a simple eigenvalue for $-\Delta_{\Omega_h}$, we can strengthen the statement to hold for all $\rho \in L^2(\partial\Omega_h)$ agreeing with a normal variation of the domain, i.e.\ $\rho = (s\circ h^{-1}) \cdot \nu$ on $\partial \Omega_h$ for some $s\in C^k(\Omega, \bR^n)$. 
    In fact, by \eqref{eqn:TQ} we know $(v,\eta,s) \in T_{(u,\lambda, h)} Q_k$ if and only if
    \[ (\Delta_{\Omega_h} + \lambda) (v\circ h^{-1}) + \eta u \circ h^{-1} - (\Delta_{\Omega_h} + \la)((s \circ h^{-1}) \cdot \nabla (u \circ h^{-1})) =0, \qquad \text{ on } \Om_h. \]
    Under the assumption of simplicity, it suffices to take
    \[ \eta := -\dfrac{\int_{\Omega_h} (u\circ h^{-1}) \left( \Delta_{\Omega_h} + \lambda \right) \left( (s\circ h^{-1}) \cdot \nabla (u\circ h^{-1} ) \right)}{\int_{\Omega_h} |u\circ h^{-1}|^2 } = -\dfrac{\int_{\partial\Omega_h} (s\circ h^{-1}) \cdot \nu \left| \frac{\partial (u\circ h^{-1})}{\partial \nu} \right|^2}{{\int_{\Omega_h} |u\circ h^{-1}|^2 }}. \]
    This choice guarantees that $g:= (\De_{\Om_h} + \la)((s \circ h^{-1})\cdot \nabla (u \circ h^{-1})) - \eta \, u\circ h^{-1}$ is orthogonal to $\ker(\Delta_{\Omega_h} + \lambda) = \mathrm{span} \{u\circ h^{-1} \}$, and thus by the Fredholm alternative we can solve for a unique $w \perp u\circ h^{-1}$ such that $(\Delta_{\Omega_h} + \lambda) w = g$. Therefore, letting $v:= w\circ h - \langle w\circ h, u \rangle_{L^2(\Omega)} u$, we have $(v,\eta, s) \in T_{(u,\lambda,h)} Q_k$. Unfortunately this stronger statement only holds when $\lambda$ is a simple eigenvalue; otherwise, even with the flexibility of choosing $\eta \in \bR$ we still cannot recover the $m$-dimensional space $W$ with $m\geq 2$.
\end{rmk}

\begin{proof}

Composing the defining equation of $T_{(u, \la, h)} Q_k$ with $h^{-1}$ and letting $\eta =0$, we find
\begin{align}\label{eqn:defining-equation-TQ-mod}
    (\De_{\Om_h} + \la)(v \circ h^{-1}) - (\De_{\Om_h} + \la)((s \circ h^{-1})\cdot \nabla (u \circ h^{-1})) = 0, \qquad \text{ on } \Om_h.
\end{align}
We next study the solvability of ~\eqref{eqn:defining-equation-TQ-mod} for $v$, given $s$. By the Fredholm alternative (see ~\cite[Theorem 8.6]{gilbarg_trudinger_1977}), given $s \in C^k(\Om, \bR^n)$, there is $v \in H^2 \cap H_0^1(\Om)$ solving ~\eqref{eqn:defining-equation-TQ-mod} if
\begin{align}\label{eqn:solvability-condition-TQ}
    \int_{\Om_h} f (\De_{\Om_h} + \la)((s \circ h^{-1})\cdot \nabla (u \circ h^{-1})) = 0 \tag{$\star$}
\end{align}
for each $f \in \ker(\De_{\Om_h}+ \la)^* = \ker(\De_{\Omega_h}+\lambda)$ with $f=0$ on $\dd\Om_h$.

After integration by parts and using that $f$ is a Dirichlet eigenfunction of $\Om_h$ with eigenvalue $\la$, ~\eqref{eqn:solvability-condition-TQ} becomes
\begin{align}\label{eqn:solvability-condition-TQ-2}
    0 &= \int_{\dd\Om_h} (s \circ h^{-1})\cdot \nabla (u \circ h^{-1}) \frac{\dd f}{\dd \nu} =  \int_{\dd\Om_h} (s \circ h^{-1})\cdot\nu \frac{\dd(u \circ h^{-1})}{\dd \nu} \frac{\dd f}{\dd \nu},
\end{align}
where we used that $u \circ h^{-1}=0$ on $\dd\Om_h$ in the second equality. 
In other words, \eqref{eqn:solvability-condition-TQ-2}, and thus the condition \eqref{eqn:solvability-condition-TQ}, is equivalent to $(s\circ h^{-1})\cdot \nu \in W^\perp$.

Consider the dense subspace $Y \subset L^2(\dd \Om_h)$ given by
\begin{align}
    Y &:= \{(s \circ h^{-1}) \cdot \nu \in L^2(\dd\Om_h)\,|\, s \in C^k(\Om, \bR^n)\}\subset L^2(\dd \Om_h).
\end{align}
Since $Y$ is a dense subspace of $L^2(\dd\Om_h)$ and $W^\perp$ is a finite-codimensional subspace, 
\begin{align}
Z:=Y \cap W^\perp
\end{align}
is dense in $W^\perp$. 

Indeed, let $\pi := \mathrm{proj}_W: L^2(\dd\Om_h) \to W \cong \bR^m$ be the projection map, which exists by the finite-dimensionality of $W$. Since $Y$ is a dense subspace, $\pi(Y) \subset \bR^m$ is a dense subspace and hence $\pi|_Y: Y \to W$ is surjective. By the surjectivity of $\pi|_Y$ and the finite-dimensionality of $W$, there exists a continuous linear right inverse of $\pi|_Y$ which we denote by $R: W \to Y$. Let $b \in W^\perp$ be arbitrary, and let $x_i \in Y$ be a sequence such that $x_i \to b$ in $L^2(\dd\Om_h)$; since $Y$ is a subspace, $x_i - R(\pi|_Y(x_i)) \in Y$. Since $\pi(x_i - R(\pi|_Y(x_i))) = 0$, $x_i - R(\pi|_Y(x_i)) \in W^\perp$. So, $x_i - R(\pi|_Y(x_i)) \in Z$, and $x_i - R(\pi|_Y(x_i)) \to b - R(0) = b$ in $L^2(\dd\Om_h)$. This proves that $Z$ is dense in $W^\perp$.

Let $\rho \in Z$. Since $\rho \in Y$, there exists $s \in C^k(\Om, \bR^n)$ by definition such that on $\dd\Om_h$,
\begin{align}\label{eqn:rho-on-boundary}
    \rho = (s \circ h^{-1}) \cdot \nu.
\end{align}
We also have $\rho \in W^\perp$, so the condition \eqref{eqn:solvability-condition-TQ} is verified for such $s$. 
Hence, the previous discussion implies there exists $\tilde{v} \in H^2 \cap H_0^1(\Om)$ solving ~\eqref{eqn:defining-equation-TQ-mod}, and equivalently,
\begin{align}\label{eqn:solvability-TQ-4}
(L_h + \la) \tilde{v} - (\Delta_{\Om_h} + \la)((s \circ h^{-1}) \cdot \nabla (u \circ h^{-1})) \circ h=0, \qquad \text{ on }\Om.
\end{align}
We have $||u||_{L^2(\Om)}=1$ since $u \in S$, and we normalize $\tilde{v}$ to
$$v := \tilde{v} - \langle \tilde{v}, u\rangle_{L^2(\Om)} u \in T_u S. $$

Since $\varphi(u, \la, h) = (L_h + \la)u=0$, $(v, 0, s) \in T_u S \times \bR \times C^k(\Omega, \bR^n)$ still satisfies the defining equation ~\eqref{eqn:TQ} for $T_{(u, \la, h)}Q_k$ with $\eta =0$:
\begin{align}\label{eqn:defining-equation-TQ-mod-2}
    (L_h + \la) v - (\Delta_{\Om_h} + \la)((s \circ h^{-1}) \cdot \nabla (u \circ h^{-1})) \circ h=0,
\end{align}
and so
\begin{align}\label{eqn:(v,0,s)}
(v, 0, s) \in T_{(u, \la, h)}Q_k.
\end{align}

Overall, we have shown that $Z$ is dense in $W^\perp$ and for each $\rho \in Z$ and any $s \in C^k(\Om, \bR^n)$ such that ~\eqref{eqn:rho-on-boundary} holds, we have found $v \in H^2 \cap H_0^1(\Om)$ such that~\eqref{eqn:(v,0,s)} holds. This concludes the lemma. 
\end{proof}

The differential $d\alpha|_{(u, \la, h, z)}: T_{(u, \la, h)}Q_k \times T_z \inte \Om \to T_{\alpha(u, \la, h, z)}\bR \cong \bR$ is given by
\begin{align}\label{eqn:differential-alpha}
    d\alpha|_{(u, \la, h, z)}(v, \eta, s, y) = v(z) + y \cdot \nabla u(z).
\end{align}
We remark that the domain perturbation parameter $s$ does not enter into the differential above, but it certainly determines the allowable functions $v$. Recall that the function $v$ such that $(v,\eta,s)\in T_{(u,\lambda,h)}Q_k$ has improved regularity $v\in C^{k-1,\alpha}(\Omega)$, and $k-1 > 0$, so the pointwise evaluation of $v$ makes sense.

If $\nabla u (z) \neq 0$, then we may choose $y = \nabla u(z)$ so that
\begin{align}
    d\alpha|_{(u, \la, h, z)} (0,0,0, y) = |\nabla u(z)|^2 \neq 0.
\end{align}
In this case, $d\alpha|_{(u, \la, h, z)}$ is surjective. So, without loss of generality, we may assume $\nabla u(z) = 0$. 
Suppose that $0$ is not a regular value of $\alpha$. Then, there is $(u, \la, h, z) \in Q_k \times \inte \Om$ such that 
\begin{align}\label{eqn:u(z)=0}
\alpha(u, \la, h, z) = u(z) &= 0,\\
    \nabla u(z) &= 0,
\end{align}
and
\begin{align}\label{eqn:dalpha-v(z)=0}
   d\alpha|_{(u, \la, h, z)}(v, \eta, s, y) = v(z) = 0, \qquad \text{  for all } (v, \eta, s, y) \in T_{(u, \la, h)}Q_k \times T_z \inte \Om.
\end{align}

Let $\{\psi_i\}_{i=1}^m$ be an $L^2(\Om_h)$-orthonormal basis of Dirichlet eigenfunctions with eigenvalue $\lambda$. Define $v_i := \psi_i \circ h - \langle\psi_i \circ h, u\rangle_{L^2(\Om)}u$, so that $v_i \in T_u S$, and
$(L_h + \la)v_i=0$. By ~\eqref{eqn:TQ}, $(v_i, 0, 0) \in T_{(u, \la, h)}Q_k$. And by ~\eqref{eqn:dalpha-v(z)=0}, $v_i(z)=0$. Since $u(z)=0$ also,
\begin{align}\label{eqn:psii(h(z))=0}
\psi_i(h(z))=0
\end{align}
for $i = 1, \dots, m$. As a side remark, notice that the same argument (of using the reduced Green's function as a test function) cannot work to show any nonzero real number $c$ is a regular value of $\alpha$. This points to the well-known fact that the zero level set of eigenfunctions is different from other level sets.

Let $G^h$ be the reduced Green's function on $\Om_h$ for the operator $\De_{\Om_h} + \la$. In light of the conditions above, we consider the reduced Green's function with pole at $h(z) \in \Omega_h$, denoted as $G^h_{h(z)} := G^h(h(z), \cdot)$.

The next claim helps us make a suitable choice of $\rho$ to which we apply Lemma \ref{lemma:solvability-for-TQ}.
\begin{claim}\label{Claim:notorthogonalGh} Under the assumption that zero is not a regular value of $\alpha$, with given $(u, \la, h, z) \in Q_k \times \inte \Om$, 
\begin{align}
    \frac{\dd(u \circ h^{-1})}{\dd \nu}\frac{\dd G^h_{h(z)}}{\dd \nu} \notin W,
\end{align}
where $W \subset L^2(\dd \Om_h)$ is the finite-dimensional linear subspace as defined in \eqref{eqn:W-subspace-definition}.
\end{claim}
\begin{proof}
Suppose the claim is not true. Then, there exist coefficients $d_1, \dots, d_m \in \bR$ such that
\begin{align}
    \frac{\dd(u \circ h^{-1})}{\dd \nu} \frac{\partial}{\partial \nu}\left(G^h_{h(z)} - \sum_{i=1}^m d_i\psi_i\right)  \equiv 0.
\end{align}
Since $u \circ h^{-1}$ is a non-trivial Dirichlet eigenfunction on $\Om_h$, $\frac{\partial (u \circ h^{-1})}{\partial \nu}$ cannot vanish identically on any relatively open subset of $\dd \Om_h$. So
\begin{align}\label{eqn:boundary-zero-Gh}
    \frac{\partial}{\partial \nu}\left(G^h_{h(z)} - \sum_{i=1}^m d_i \psi_i\right) &= 0, \qquad \text{ on }\dd\Om_h.
\end{align}

By definition of $G^h_{h(z)}$, $(\De_{\Om_h}+\la)G^h_{h(z)} = -\de_{h(z)} \circ \Pi_{(\ker(\De_{\Om_h}+\la))^\perp}$. 
By assumption we have $\psi_i(h(z))=0$ in ~\eqref{eqn:psii(h(z))=0}, and thus $(\De_{\Om_h}+\la)G^h_{h(z)} = -\de_{h(z)}$. On $\Om_h \setminus \{h(z)\}$,
$$(\De_{\Om_h}+\la)\left(G^h_{h(z)} - \sum_{i=1}^m d_i \psi_i\right) = 0.$$
Together with ~\eqref{eqn:boundary-zero-Gh}, the uniqueness of the Cauchy problem implies that $G^h_{h(z)} - \sum_{i=1}^m d_i \psi_i = 0$ on $\Om_h \setminus \{h(z)\}$. In particular, it follows that $h(z)$ is a removable singularity for $G^h_{h(z)}$, which contradicts the definition of the reduced Green's function. We conclude the claim.
\end{proof}

We may consider the orthogonal projection onto $W^\perp$ since $W^\perp$ is closed. By Claim \ref{Claim:notorthogonalGh}, 
\begin{align}
    a:= \mathrm{proj}_{W^\perp}\left( \frac{\dd(u \circ h^{-1})}{\dd \nu}\frac{\dd G^h_{h(z)}}{\dd \nu}\right) \text{ is a non-trivial function in } L^2(\partial\Omega_h),
\end{align}
and in particular
\begin{align}
    \int_{\dd\Om_h} a \frac{\dd(u \circ h^{-1})}{\dd \nu}\frac{\dd G^h_{h(z)}}{\dd \nu} \neq 0. 
\end{align}

By Lemma \ref{lemma:solvability-for-TQ}, the set $Z\subset W^{\perp}$ is dense, and thus there is a sequence $\rho_i \in Z$  such that $\rho_i \to a$ in $L^2(\dd\Om_h)$. 
So, there exists $\rho \in Z$ such that 
\begin{align}\label{eqn:nonzero-integral-with-rho}
\int_{\dd \Om_h} \rho \frac{\partial}{\partial \nu} (u \circ h^{-1}) \frac{\partial}{\partial \nu}G^h_{h(z)} \neq 0.\end{align}
Applying Lemma \ref{lemma:solvability-for-TQ} to this choice of $\rho$, there is $v \in H^2 \cap H_0^1(\Om)$ and $s \in C^k(\Om, \bR^n)$ such that $(v, 0, s) \in T_{(u, \la, h)}Q_k$ and $\rho = (s\circ h^{-1})\cdot \nu$ on $\dd\Om_h$. By ~\eqref{eqn:nonzero-integral-with-rho} and $u \circ h^{-1}=0$ on $\dd\Om_h$,
\begin{equation}\label{eqn:nonzero-integral-rho2}
    \int_{\dd \Om_h} (s \circ h^{-1}) \cdot \nabla (u \circ h^{-1}) \frac{\partial}{\partial \nu}G^h_{h(z)} = \int_{\dd \Om_h} (s \circ h^{-1}) \cdot \nu \frac{\partial}{\partial \nu} (u \circ h^{-1}) \frac{\partial}{\partial \nu}G^h_{h(z)} \neq 0.
\end{equation}
Now, we apply the defining equation of ~\eqref{eqn:TQ} to $(v, 0, s)$, compose with $h^{-1}$, and integrate against $G^h_{h(z)}$,
\begin{align}\label{eqn:integral-alpha-reg-1}
    \int_{\Om_h} G^h_{h(z)} (\Delta_{\Om_h} + \la) (v\circ h^{-1})- \int_{\Om_h} G^h_{h(z)} (\Delta_{\Om_h} + \la)((s \circ h^{-1}) \cdot \nabla (u \circ h^{-1})) = 0.
\end{align}

Analyzing the first term of ~\eqref{eqn:integral-alpha-reg-1}, using the identity for the reduced Green's function \eqref{eq:Green} and the improved regularity $v\circ h^{-1} \in H^2 \cap H_0^1 \cap C^{1,\alpha}(\Omega_h) \subset W_0^{1,p}(\Omega_h)$ for any $p>n$, 
\begin{equation}
   \int_{\Om_h} G^h_{h(z)} (\Delta_{\Om_h} + \la) (v \circ h^{-1})  =  - v(z) + \sum_{i=1}^m \psi_i(h(z))\int_{\Om_h} (v \circ h^{-1})\psi_i = 0,
   \label{eqn:integral-alpha-reg-3}
\end{equation}
where we used~\eqref{eqn:dalpha-v(z)=0} and ~\eqref{eqn:psii(h(z))=0} in the last line.
We cannot directly use the identity \eqref{eq:Green} to analyze the second term of ~\eqref{eqn:integral-alpha-reg-1}, since the function $w:= (s \circ h^{-1}) \cdot \nabla (u \circ h^{-1})$ does not vanish at the boundary. Instead, let $\tilde{w}$ be a $C^1$ extension of the function $w|_{\partial\Omega_h} \in C^1(\partial\Omega_h)$ such that $\tilde{w}$ is supported in a neighborhood $U$ of $\partial\Omega_h$ with $h(z) \notin \overline{U}$. By applying the identity \eqref{eq:Green} to $w-\tilde{w} \in W_0^{1,p}(\Omega_h)$, applying the divergence theorem to $\tilde{w}$ and $G^h_{h(z)}$ away from the pole in $\Omega_h \cap U$, and using the assumption \eqref{eqn:psii(h(z))=0}, we obtain
\begin{align}\label{eqn:integral-alpha-reg-4}
    \int_{\Om_h} &G^h_{h(z)} (\Delta_{\Om_h} + \la)w =   \int_{\Om_h} G^h_{h(z)} (\Delta_{\Om_h} + \la)(w-\tilde{w}) + \int_{\Om_h \cap U } G^h_{h(z)} (\Delta_{\Om_h} + \la)\tilde{w}  \nonumber\\
    & = - w(h(z)) + \sum_{i=1}^m \left( \int_{\Omega_h} (w-\tilde{w}) \psi_i \right) \psi_i(h(z)) - \int_{\partial\Omega_h} \tilde{w} \frac{\partial}{\partial \nu} G^h_{h(z)} \nonumber \\
    &= -s(z) \cdot \nabla (u \circ h^{-1})(h(z))  - \int_{\dd \Om_h} (s \circ h^{-1}) \cdot \nabla (u \circ h^{-1}) \frac{\partial}{\partial \nu}G^h_{h(z)},\nonumber\\
\intertext{We have $\nabla (u \circ h^{-1})(h(z)) = \nabla u (z) D(h^{-1})(h(z)) = 0$ since $\nabla u(z)=0$, so}
    &= - \int_{\dd \Om_h} (s \circ h^{-1}) \cdot \nabla (u \circ h^{-1}) \frac{\partial}{\partial \nu}G^h_{h(z)}.
\end{align}
Applying ~\eqref{eqn:integral-alpha-reg-3} and ~\eqref{eqn:integral-alpha-reg-4} to ~\eqref{eqn:integral-alpha-reg-1}, we obtain
\begin{align}\label{eqn:integral-alpha-reg-5}
   \int_{\dd \Om_h} (s \circ h^{-1}) \cdot \nabla (u \circ h^{-1}) \frac{\partial}{\partial \nu}G^h_{h(z)} = v(z) - \sum_{i=1}^m \psi_i(h(z))\int_{\Om_h} (v \circ h^{-1})\psi_i = 0.
\end{align}
This is a contradiction with \eqref{eqn:nonzero-integral-rho2}, so we conclude that zero is a regular value of $\alpha$. As discussed above, we find that ~\ref{item:interior-regular-value} from Theorem \ref{thm:generic} holds, using Theorem \ref{thm:Uhlenbeck-regular-value-alpha}.

\subsection{Zero is a regular value for the boundary normal derivative}\label{sec:boundary-reg}

In this subsection, we prove \ref{item:boundary-regular-value} from Theorem \ref{thm:generic}.

Using the same notation as in the previous subsections, we define a map $\ga: Q_k \times \dd\Om \to \bR$ by
\begin{align}
    \ga(u, \la, h, z) := \frac{\dd u}{\dd \nu}(z),
\end{align}
where $\frac{\dd}{\dd \nu}$ denotes the normal derivative on $\dd \Om$ pointing outwards.
Recall that in a $C^k$ domain $\Omega_h$, any Dirichlet eigenfunction $u\circ h^{-1}$ has $C^{k-1,\beta}$ regularity up to the boundary for any $\beta \in (0,1)$. Thus, combined with the regularity of $h$, 
\[ \frac{\partial u}{\partial \nu}(z) = \frac{\partial\left( u \circ h^{-1} \right)}{\dd\nu_h}(h(z)) \, Dh|_z^T \nu_h \cdot \nu \in C^{k-2,\beta}(\partial\Omega), \]
where $\frac{\partial}{\partial \nu_h}$ denotes the normal derivative on $\partial\Omega_h$. (See the analogous discussion in \eqref{eq:change-bdnormalder}.)
Therefore, $\gamma$ is a $C^{k-2}$ map on the last variable. To show that $\gamma$ is also a regular map on the triplet $(u,\lambda, h) \in Q_k$, where $u$ is only assumed to be in the weak space $H^2 \cap H_0^1(\Omega)$, we use the improved elliptic boundary regularity of eigenfunctions. Recall that for any $1<p<\infty$ and $1 \leq j\leq k$, the operator $-L_h: W^{j, p}\cap W_0^{1,p}(\Omega) \to W^{j-2,p}(\Omega)$ is an isomorphism, and thus its inverse $T_h := (-L_h)^{-1}$ is a bounded operator from $ W^{j-2,p}(\Omega)$ to $W^{j, p}\cap W_0^{1,p}(\Omega)$; moreover, the map $h\in B_k \mapsto T_h \in \mathcal{L}(W^{j-2,p}(\Omega), W^{j,p} \cap W_0^{1,p}(\Omega))$ is smooth as long as $k\geq 2$. By letting $p>n$ and choosing $N$ large enough (depending on the values of $k, n$ and $p$), the map
	\[ u\in H^2 \cap H_0^1(\Omega) \longmapsto \lambda^N T_h^N u \in W^{k,p} \cap W_0^{1,p}(\Omega) \hookrightarrow C^{k-1}(\Omega) \]
	is bounded. Combining all the above ingredients, we get the following map $\mathcal{B}$
	\[ (u,\lambda, h, z)\in \left( H^2 \cap H_0^1(\Omega)\right) \times \R \times B_k \times \partial\Omega \overset{C^\infty}{\longmapsto} (\lambda^N T_h^N u, z) \in C^{k-1}(\Omega) \times \partial\Omega \overset{C^{k-2}}{\longmapsto} \lambda^N  \frac{\partial}{\partial_\nu}\left(T_h^N u\right)(z). \]
	In particular, when we restrict the map $\mathcal{B}$ to eigenfunctions, $(L_h + \lambda) u= 0$ implies that $u = \lambda T_h u$, and thus $u = \lambda^N T_h^N u$. It follows that $\gamma = \mathcal{B}|_{Q_k \times \partial\Omega}$ is $C^{k-2}$ regular. We remark that the interior map $\alpha$ in Subsection \ref{sec:interior-reg} has the optimal regularity allowed by the $C^k$ class of domains, and here the regularity of $\gamma$ can not be improved to $C^{k-1}$, even for fixed domains.

\begin{thm}[\!{\cite[Theorem 4]{UhlenbeckGeneric}}]\label{thm:Uhlenbeck-regular-value-gamma} 
If $k > n$, and zero is a regular value of both $\varphi: S \times \bR \times B_k \to L^2(\Om)$ and $\ga: Q_k \times \dd\Om \to \bR$, then the set 
    $$\{h \in B_k \,|\,\text{The normal derivatives of the eigenfunctions of } L_h \text{ on } \dd\Om\text{ have }0 \text{ as a regular value}\}$$ 
is residual in $B_k$. 

In particular, if zero is a regular value of both $\varphi$ and $\ga$, then there is a residual set $\mathfrak{X} \subset \mathrm{Emb}_k(M, \bR^n)$ such that for each $F \in \mathfrak{X}$, every Dirichlet eigenfunction for the Laplacian on $F(M)$ has zero as a regular value of its normal derivative on $\partial F(M)$, i.e.\ \ref{item:boundary-regular-value} from Theorem \ref{thm:generic} holds.
\end{thm}
\begin{proof}
The first statement is a reformulation of \cite[Theorem 4]{UhlenbeckGeneric}. 
More precisely, we apply \cite[Transversality Theorem 2]{UhlenbeckGeneric} to the map $\gamma$, with finite-dimensional manifolds $X=\partial\Omega, Y= \bR, Y'=\{0\}$, assuming the regularity of $\gamma$ satisfies
\[ k-2 > \max\{0, n-2\} = n-2. \]

For the second statement, let $\Om=F_0(M)$ and
$$\mathfrak{X} = \{h \circ F_0\,|\,h \in B_k, \text{ the normal derivatives of the eigenfunctions of }L_h \text{ on $\partial\Om$ have }0 \text{ as a regular value}\}.$$
Suppose  $h \circ F_0 \in \mathfrak{X}$ and $\phi$ is a Dirichlet eigenfunction for the domain $\Om_h$. Then $\phi \circ h$ is an eigenfunction of $L_h$ and $\frac{\dd}{\dd \nu}(\phi \circ h)$ has $0$ as a regular value on $\dd \Om$. Let $p \in \partial\Om_h$, and let $q = h^{-1}(p) \in \dd \Om$. Let $\nu_h$ denote the unit normal to $\Om_h$ and let $\nu$ denote the unit normal to $\Om$. Since $\phi \circ h=0$ on $\partial \Om$, 
\begin{equation}\label{eq:change-bdnormalder}
    \frac{\dd}{\dd \nu}(\phi \circ h)(q) = \nabla (\phi \circ h)|_q \cdot \nu = \nabla \phi(p) Dh|_q \nu = \frac{\dd \phi}{\dd\nu_h}(p)\nu_h^T Dh|_q \nu = \frac{\dd \phi}{\dd\nu_h}(p) Dh|_q^T \nu_h \cdot \nu.
\end{equation}
Since $h:\Om \to \Om_h$ is a diffeomorphism, we have that $Dh|_q^T \nu_h \cdot w = 0$ for each $w \in T_q \partial\Om$, and there exists $c(q) \neq 0$ such that $Dh|_q^T \nu_h  = c(q)\nu$ and so
\begin{align}\label{eq:change-normalder}
    \frac{\dd}{\dd \nu}(\phi \circ h)(q) = c(q)\frac{\dd \phi}{\dd\nu_h}(p). 
\end{align}

Assuming $\frac{\dd \phi}{\dd \nu_h}(p)=0$, then $\frac{\dd}{\dd \nu}(\phi \circ h)(q) =0$. By assumption, $d(\frac{\dd}{\dd \nu}(\phi \circ h))|_q: T_q \partial\Om \to \bR$ is surjective.

\begin{align}
    d(\frac{\dd}{\dd \nu}(\phi \circ h))|_q &=  \frac{\dd \phi}{\dd\nu_h}(p) dc|_q + c(q) d\left( \frac{\dd\phi}{\dd\nu_h}\right)\Big|_p \circ Dh|_q \nonumber \\
    &= c(q) d\left( \frac{\dd\phi}{\dd\nu_h}\right)\Big|_p \circ Dh|_q.\label{eq:change-normalder-2}
\end{align}
Since $c(q) \neq 0$ and $Dh|_q$ is a linear isomorphism from $T_q \partial\Om$ onto $T_p \partial\Om_h$, we conclude that $d\left( \frac{\dd\phi}{\dd\nu_h}\right)\Big|_p: T_p \partial\Om_h \to \bR$ is surjective. Thus, $\phi$ has $0$ as a regular value of its normal derivative on $\partial\Om_h$. The first statement implies that $\mathfrak{X}$ is residual, so the second statement follows. 
\end{proof}

Our goal is to show that zero is a regular value of $\ga$ as long as $k\geq 3$, and the result of Section \ref{section:generic-simplicity} combined with Theorem \ref{thm:Uhlenbeck-regular-value-gamma} will allow us to conclude ~\ref{item:boundary-regular-value} from Theorem \ref{thm:generic}.

Suppose $(u,\lambda, h) \in Q_k$, and $z\in \partial\Om$. The differential $d\ga|_{(u, \la, h, z)}: T_{(u, \la, h)}Q_k \times T_z \dd\Om \to T_{\ga(u, \la, h, z)}\bR \cong \bR$ is given by
\begin{align}\label{eqn:differential-gamma}
    d\ga|_{(u, \la, h, z)}(v, \eta, s, y) = \frac{\dd v}{\dd \nu}(z) + y \cdot \nabla\left(\frac{\dd u}{\dd \nu}\right) (z).
\end{align}
We note that $\ga$ is defined as a normal derivative on the fixed domain $\Om$, not the normal derivative of Laplacian eigenfunctions on $\Om_h$. So the domain $\Om_h$ and its variation $s$ do not appear in the differential above, and they only affect the allowable perturbations $v$ and $\eta$. And $v$ has improved regularity, so that the pointwise evaluation of the boundary normal derivative of $v$ is well-defined.

Suppose that $0$ is not a regular value of $\ga$. Then, there is $(u, \la, h, z) \in Q_k \times \dd \Om$ such that 
\begin{align}\label{eqn:du/dnu(z)=0}
\ga(u, \la, h, z) = \frac{\dd u}{\dd \nu}(z) &= 0,
\end{align}
and
\begin{align}\label{eqn:dgamma-v(z)=0--0}
   d\ga|_{(u, \la, h, z)}(v, \eta, s, y) = \frac{\dd v}{\dd \nu}(z) + y \cdot \nabla\left(\frac{\dd u}{\dd \nu}\right) (z) = 0, \qquad \text{  for all } (v, \eta, s, y) \in T_{(u, \la, h)}Q_k \times T_z \dd \Om.
\end{align}
In particular, any $(v, \eta, s, 0) \in T_{(u, \la, h)}Q_k \times T_z \dd \Om$ satisfies
\begin{align}\label{eqn:dgamma-v(z)=0--1}
    \frac{\dd v}{\dd \nu}(z) &=0,\\\label{eqn:nablav(z)}
   \nabla v(z) &= 0,
\end{align}
where the last line follows from the first since $v=0$ on $\dd\Om$.

Let $\{\psi_i\}_{i=1}^m$ be an $L^2(\Om_h)$-orthonormal basis of Dirichlet eigenfunctions with eigenvalue $\lambda$. Let $v_i := \psi_i \circ h - \langle\psi_i \circ h, u\rangle_{L^2(\Om)}u$. We have $v_i \in T_u S$, and
$(L_h + \la)v_i=0$. By ~\eqref{eqn:TQ}, $(v_i, 0, 0) \in T_{(u, \la, h)}Q_k$. And by ~\eqref{eqn:dgamma-v(z)=0--1}, $\frac{\dd v_i}{\dd \nu}(z) = 0$. Then, ~\eqref{eqn:du/dnu(z)=0} implies that $\frac{\dd (\psi_i \circ h)}{\dd \nu}(z)=0$ for $i = 1, \dots, m$. We have that $\nabla(\psi_i \circ h)(z)=0$ since $\psi_i \circ h =0$ on $\dd \Om$. So, $\nabla \psi_i(h(z))=0$ and we find that 
\begin{align}\label{eqn:ddnupsii(h(z))=0}
\frac{\dd \psi_i}{\dd \nu}(h(z))=0,
\end{align}
where the normal derivative is taken with respect to $\dd\Om_h$. 

In light of the assumptions above, consider $P^h_{h(z)}$, the reduced Poisson kernel for the operator $\De_{\Om_h} + \la$ anchored at $h(z) \in \partial\Omega_h$, see Lemma \ref{lem:reduced-poisson-kernel}.
Let $\Gamma \subset \dd\Om_h$ be a relatively open subset of $\dd\Om_h$ such that $h(z) \notin \overline{\Gamma}$. Let $W_\Ga\subset L^2(\Ga)$ be the subspace defined by
\begin{align}
    W_\Ga &:= \mathrm{span}\left\{\frac{\dd(u \circ h^{-1})}{\dd \nu} \frac{\dd \psi_1}{\dd \nu}\Big|_{\Ga}, \cdots, \frac{\dd(u \circ h^{-1})}{\dd \nu} \frac{\dd \psi_m}{\dd \nu}\Big|_{\Ga}\right\}\subset L^2(\Ga).
\end{align}

The next claim is analogous to Claim \ref{Claim:notorthogonalGh}. 
\begin{claim}\label{Claim:notorthogonalPh} Under the assumption that zero is not a regular value of $\ga$, with given $(u, \la, h, z) \in Q_k \times \dd \Om$, 
\begin{align}
    \frac{\dd(u \circ h^{-1})}{\dd\nu} \frac{\dd}{\dd\nu}P^h_{h(z)}\Big|_{\Ga} \notin W_\Ga.
\end{align}
\end{claim}
\begin{proof}
Suppose the claim is not true. Then, there exist $d_1, \dots, d_m \in \bR$ such that on $\Gamma$, 
\begin{align}
    \frac{\dd(u \circ h^{-1})}{\dd \nu} \frac{\partial}{\partial \nu}\left(P^h_{h(z)} - \sum_{i=1}^m d_i\psi_i\right)  = 0.
\end{align}
Since $u \circ h^{-1}$ is a non-trivial Dirichlet eigenfunction on $\Om_h$, $\frac{\partial (u \circ h^{-1})}{\partial \nu}$ vanishes at most on a $\mathcal{H}^{n-1}$-null set on $\Gamma$. So, on $\Gamma$,
\begin{align}\label{eqn:deriv-boundary-zero-Ph}
    \frac{\partial}{\partial \nu}\left(P^h_{h(z)} - \sum_{i=1}^m d_i \psi_i\right) &= 0.
\end{align}

Recall that $(\De_{\Om_h}+\la)P^h_{h(z)} = 0$ in $\Om_h$, so
\begin{align}
    (\De_{\Om_h}+\la)\left(P^h_{h(z)} - \sum_{i=1}^m d_i \psi_i\right)=0
\end{align}
in $\Om_h$. Moreover, $P^h_{h(z)} - \sum_{i=1}^m d_i \psi_i = 0$ on $\Gamma$ since $h(z) \notin \overline{\Gamma}$. Thus the uniqueness of the Cauchy problem implies that $P^h_{h(z)} = \sum_{i=1}^m d_i \psi_i$ in a neighborhood of $\Gamma$, and thus by unique continuation, they must agree on all of $\Om_h$. This contradicts the definition of $P^h_{h(z)}$ in \eqref{eq:Poisson-kernel}, and we conclude the claim.
\end{proof}

Since $W^\perp_\Ga$ has finite codimension in $L^2(\Gamma)$, it is closed. Consider the orthogonal projection onto $W^\perp_\Ga$. By Claim \ref{Claim:notorthogonalPh}, 
\begin{align}
    a:= \mathrm{proj}_{W^\perp_\Ga}\left( \frac{\dd(u \circ h^{-1})}{\dd \nu}\frac{\dd}{\dd \nu}P^h_{h(z)}\right) \text{ is a non-trivial function in } L^2(\Gamma),
\end{align}
and in particular
\begin{align}
    \int_{\Ga} a \frac{\dd(u \circ h^{-1})}{\dd \nu}\frac{\dd}{\dd \nu}P^h_{h(z)} \neq 0. 
\end{align}
Let
$$Y_\Ga := \{(s \circ h^{-1}) \cdot \nu \in L^2(\Gamma)\,|\, s \in C^k(\Om, \bR^n) \text{ and } (s\circ h^{-1})|_{\partial\Omega_h} \text{ has compact support in }\Gamma \}\subset L^2(\Ga).$$
We have that $Y_\Ga$ is dense in $L^2(\Ga)$. Since $Y_\Ga$ is a dense subspace of $L^2(\Ga)$ and since $W^\perp_\Ga$ has finite codimension in $L^2(\Ga)$, we find that $Y_\Ga \cap W^\perp_\Ga$ is dense in $W^\perp_\Ga$. This follows from the same argument that $Y\cap W^\perp$ is dense in $W^\perp$ in the proof of Lemma \ref{lemma:solvability-for-TQ}. So, there exists $\rho \in Y_\Ga \cap W^\perp_\Ga$ such that 
\begin{align}\label{eqn:nonzero-integral-with-rho-for-gamma}
\int_{\Ga} \rho \frac{\partial}{\partial \nu} (u \circ h^{-1}) \frac{\partial}{\partial \nu}P^h_{h(z)} \neq 0.\end{align}
By construction of $Y_\Ga$, $\rho$ is compactly supported on $\Ga$ and can be extended by zero to $\rho: \dd\Om_h \to \bR$, which we denote with the same symbol. Thus,
\begin{align}\label{eqn:nonzero-integral-with-rho-for-gamma-2}
\int_{\dd\Om_h} \rho \frac{\partial}{\partial \nu} (u \circ h^{-1}) \frac{\partial}{\partial \nu}P^h_{h(z)} \neq 0.\end{align}
Likewise, $\rho \in W^\perp$ by extension.

Applying Lemma \ref{lemma:solvability-for-TQ} to this choice of $\rho$ and letting $s \in C^k(\Om, \bR^n)$ such that $(s\circ h^{-1})\cdot \nu = \rho$ on $\dd\Om_h$ and $s \circ h^{-1} \equiv 0$ in a small neighborhood of $h(z)$, there is $v \in H^2 \cap H_0^1(\Om)$ such that $(v, 0, s) \in T_{(u, \la, h)}Q_k$. By ~\eqref{eqn:nonzero-integral-with-rho-for-gamma-2} and $u \circ h^{-1}=0$ on $\dd\Om_h$,
\begin{equation}\label{eqn:nonzero-integral-rho2-gamma}
    \int_{\dd \Om_h} (s \circ h^{-1}) \cdot \nabla (u \circ h^{-1}) \frac{\partial}{\partial \nu}P^h_{h(z)} = \int_{\dd \Om_h} (s \circ h^{-1}) \cdot \nu \frac{\partial}{\partial \nu} (u \circ h^{-1}) \frac{\partial}{\partial \nu}P^h_{h(z)} \neq 0.
\end{equation}

Now, we apply the defining equation of ~\eqref{eqn:TQ} to $(v, 0, s)$, compose with $h^{-1}$, and integrate against $P^h_{h(z)}$,
\begin{align}\label{eqn:integral-alpha-reg-1-ga}
    \int_{\Om_h} P^h_{h(z)} (\Delta_{\Om_h} + \la) (v\circ h^{-1})- \int_{\Om_h} P^h_{h(z)} (\Delta_{\Om_h} + \la)((s \circ h^{-1}) \cdot \nabla (u \circ h^{-1})) = 0.
\end{align}
Analyzing the first term of ~\eqref{eqn:integral-alpha-reg-1-ga}, using $v \circ h^{-1} \in W^{2,p} \cap W^{1,p}_0 (\Omega_h)$ for $p>n$, the identity for the reduced Poisson kernel \eqref{eq:Poisson-kernel} and the assumption \eqref{eqn:ddnupsii(h(z))=0},
\begin{align}
   \int_{\Om_h} P^h_{h(z)} (\Delta_{\Om_h} + \la) (v \circ h^{-1})  &=  \frac{\dd(v \circ h^{-1})}{\dd\nu}(h(z)) = 0,\label{eqn:integral-alpha-reg-3-ga}
\end{align}
where the last equality follows since $\nabla v(z)=0$. 
Analyzing the second term of ~\eqref{eqn:integral-alpha-reg-1-ga}, using the assumption that $(s \circ h^{-1}) \cdot \nabla (u \circ h^{-1})\equiv 0$ in a neighborhood of $h(z)$ (denoted as $B_r(h(z))$), the identity of the reduced Poisson kernel \eqref{eq:Poisson-kernel-boundary}, $u\circ h^{-1} \in C^2(\Omega_h) \subset W^{2,p}(\Omega_h)$ and a density argument,\footnote{We may approximate the function $(s\circ h^{-1}) \cdot \nabla (u\circ h^{-1}) \in W^{1,p}(\Omega_h)$ by smooth functions $w_i \in C^2(\Omega_h)$, which vanish on $B_r(h(z))$, then pass the integral $\int_{\Omega_h} \nabla P^h_{h(z)} \cdot \nabla w_i$ to the limit.}
\begin{align}\label{eqn:integral-alpha-reg-4-ga}
    - \int_{\Om_h \setminus B_r(h(z)) } &P^h_{h(z)} (\Delta_{\Om_h} + \la)((s \circ h^{-1}) \cdot \nabla (u \circ h^{-1})) \nonumber\\
    &= - \frac{\dd}{\dd\nu}\left((s\circ h^{-1}) \cdot \nabla (u \circ h^{-1})\right)(h(z)) + \int_{\dd \Om_h} (s \circ h^{-1}) \cdot \nabla (u \circ h^{-1}) \frac{\partial}{\partial \nu}P^h_{h(z)}\nonumber\\
    &= \int_{\dd \Om_h} (s \circ h^{-1}) \cdot \nabla (u \circ h^{-1}) \frac{\partial}{\partial \nu}P^h_{h(z)}.
\end{align}
Applying ~\eqref{eqn:integral-alpha-reg-3-ga} and ~\eqref{eqn:integral-alpha-reg-4-ga} to ~\eqref{eqn:integral-alpha-reg-1-ga},
\begin{align}
   0 &= \int_{\dd \Om_h} (s \circ h^{-1}) \cdot \nabla (u \circ h^{-1}) \frac{\partial}{\partial \nu}P^h_{h(z)}.
\end{align}
However, with our choice of $(v,0,s)$, the above integral is nonzero by~\eqref{eqn:nonzero-integral-rho2-gamma}. This is a contradiction, so we conclude that zero is a regular value of $\ga$. As discussed above, we find that ~\ref{item:boundary-regular-value} from Theorem \ref{thm:generic} holds, using Theorem \ref{thm:Uhlenbeck-regular-value-gamma}. 

\section{Boundary asymptotics of Dirichlet eigenfunctions}\label{sec:bd}

\begin{prop}\label{prop:bdexp}
	Let $\Omega$ be an open $C^{2,\alpha}$ domain in $\bR^n$, with $\alpha\in (0,1)$ and $n\geq 2$. Suppose $u$ is a Dirichlet Laplacian eigenfunction in $\Omega$, namely $u$ satisfies
	\[ \left\{\begin{array}{ll}
		-\Delta u = \lambda u, & \text{ in } \Omega \\
		\quad\,\,\,\, u=0, & \text{ on } \partial\Omega
	\end{array} \right. \]
	for some $\lambda \in \bR$. Then for any $p\in \partial\Omega$ and any $0<\beta<1$, there exist some $r_0>0$ and an extension $\tilde{u}$ of $u$ across the boundary such that 
	\begin{equation}\label{eq:expansion}
		\tilde{u}(x) = P(x-p) + \psi(x-p), \quad \text{ in } B_{r_0}(p) 
	\end{equation} 
	where $P$ is a non-trivial homogeneous harmonic polynomial of degree $k\in \mathbb{N}$, and the remaining term $\psi$ satisfies
	\begin{equation}\label{est:grad1}
		|\psi(z)| = O(|z|^{k+\beta}), \quad \text{ and } \quad |\nabla \psi(z)| = O(|z|^{k+\beta-1}). 
	\end{equation}  
	Furthermore, we have 
	\begin{equation}\label{est:grad2}
		|\nabla^2 \psi(z)| = O(|z|^{k+\beta-2}), \qquad \text{ for all } z\in B_{r_0} (0) \cap (\overline{\Omega} - p) \text{ and } z\in B_{r_0}(0) \cap (\Omega^c - p ),
	\end{equation}
	where for $z\in \partial(\Omega-p)$, $\nabla^2 \psi(z)$ is understood to be either the one-sided derivatives from $\Omega$ or $\overline{\Omega}^c$, respectively, though it may have jump discontinuities across the boundary.
\end{prop}

\emph{Step one.} Firstly, we use the reflection principle to construct the extension $\tilde{u}$ of $u$, which moreover satisfies a divergence-form elliptic equation with Lipschitz coefficients. The construction is well-known to experts studying boundary regularity of elliptic PDEs, but we provide the detailed proof for the sake of completeness.
\begin{lem}\label{lm:ext}
    Let $\Omega \subset \bR^n$ be a $C^{k}$ domain with $k\geq 3$. Suppose $u$ solves 
    \begin{equation}
        \left\{ \begin{array}{ll}
             -\Delta u = \lambda u & \text{ in } \Omega \\
            u=0 & \text{ on } \partial \Omega.  
        \end{array} \right.
    \end{equation} 
     Then for any $p\in \partial\Omega$, there exist some $r>0$ and an extension $\tilde{u}$ of $u$ across the boundary satisfying
     \begin{equation}\label{eq:PDEext}
         -\divg(B\nabla \tilde{u}) = \lambda c \tilde{u}, \quad \text{ in } B_{2r}(p),
     \end{equation}
     where $B(y,t)$ is an elliptic symmetric $n\times n$ matrix, $c(y,t)\geq 0$ is a scalar function, and both have Lipschitz regularity.
\end{lem}
\begin{rmk}\label{rmk:reg}
	In fact, using a slightly more complicated diffeomorphism than the map $G$ in the proof below, the extension works if $\Omega$ is only assumed to be a $C^{1,1}$ domain. See the construction in \cite[Section 4]{KZcritical} with parameter $\theta(s) \approx s$. However, the regularity of the new coefficient matrix $B$ does not improve as we improve the regularity of the domain: $B$ is Lipschitz but not $C^1$ across $\partial\Omega$, due to the curvature of $\partial\Omega$. This lack of regularity makes our analysis more difficult. In particular, $\tilde{u}$ is not necessarily in $C^2$.
\end{rmk}

\begin{proof}
    Without loss of generality assume $p=0$, and that
    \[ \Omega \cap B_{R}(0) = \{(x,t) \in \mathbb{R}^{n-1}\times \mathbb{R}: t> \varphi(x) \} \cap B_{R}(0) \]
    for some $R>0$ and some function $\varphi \in C^k(\mathbb{R}^{n-1};\mathbb{R})$ satisfying $\varphi(0) = |\nabla \varphi(0)| = 0$. Then the unit normal vector at $(x,\varphi(x))$ pointing inwards is given by
    \[ n_x = \dfrac{(-\nabla \varphi(x), 1) }{\sqrt{1+|\nabla \varphi(x)|^2}}. \]
    Consider the following transformation $G: \mathbb{R}^n \to \mathbb{R}^n$ that flattens the domain $\Omega$, where
    \[ G: (x,s)\in \mathbb{R}^{n-1} \times \mathbb{R} \mapsto (x,\varphi(x)) + s \, n_{x}.\footnote{If the domain $\Omega$ (and thus the function $\varphi$) is only assumed to have $C^{1,1}$ regularity, we instead define $G(x,s) := (x,\varphi(x)) + s \, \rho_s \ast n_x $, where $\rho_s$ is an approximation of the identity in $\mathbb{R}^{n-1}$. After mollifying the unit normal vector $n_x$, it suffices to assume less regularity for $\Omega$. We invite interested readers to see the details in \cite[Section 4]{KZcritical} (with the parameter $\theta(s)\approx s$), but the essence of the argument is the same as is presented here for a simpler setup. } \]
    In the literature, these are sometimes called Fermi coordinates of $\Omega$.
    Note that
    \begin{equation}
        DG(x,s) = \begin{pmatrix}
			\delta_{ij} - s \partial_j \left( \frac{\partial_i \varphi(x)}{\sqrt{1+|\nabla \varphi(x)|^2} } \right) & -\frac{\left( \nabla \varphi(x) \right)^T}{\sqrt{1+|\nabla \varphi(x)|^2}} \\
			\nabla \varphi(x) + s \nabla \left(\frac{1}{\sqrt{1+|\nabla \varphi(x)|^2}} \right) & \frac{1}{\sqrt{1+|\nabla \varphi(x)|^2}}
        \end{pmatrix}.
    \end{equation}
    In particular, when $|(x,s)|< 2r$ is small enough, we have that $1/2 \leq \det DG(x,s) \leq 3/2$ and it is invertible. Also, the matrix $DG(x,s)$ has $C^{k-2}$ regularity. Thus, $G$ is a local $C^{k-1}$ diffeomorphism between $\Omega$ and $\bR^n_+$, as well as between $\overline{\Omega}^c$ and $\bR^n_-$.
    
    Let $v(x,s) = u\circ G(x,s)$ in $B_{2r}^+(0) = B_{2r}(0)\cap \mathbb{R}^n_+$. For any $\rho \in C_c^{\infty}(G(B_{2r}^+(0)))$, we have that
\begin{align*}
    \int_{B^+_{2r}(0)} \nabla(\rho \circ G) \cdot A \nabla v = \int_{B^+_{2r}(0)} \lambda \left|\det DG(x,s)\right| (\rho \circ G) v,
\end{align*}
where the matrix $A(x,s)$ is given by
    \begin{equation}\label{eq:Adt}
    	A(x,s) = \left|\det DG(x,s)\right| (DG(x,s))^{-1} ((DG(x,s))^T)^{-1}. 
    \end{equation} 
Thus, $v$ satisfies
    \[ \left\{\begin{array}{ll}
    	-\divg(A \nabla v) = \lambda \left|\det DG(x,s)\right|\, v(x,s) & \text{ in } B_{2r}^+(0), \\
    	\hspace{.7in} v =0 & \text{ on } B_{2r}(0) \cap \partial \mathbb{R}^n_+.
    \end{array} \right. \]
    
    Let $T(x)$ denote the $(n-1)$-dimensional square matrix given by
    \begin{align}
        T(x) = (\Id_{n-1} + (\nabla \varphi(x))^T\nabla \varphi(x))^{-1}.
    \end{align}
    Linear algebra calculations show that
    \begin{equation}\label{tmp:laOs}
    	(DG(x,s))^{-1} ((DG(x,s))^T)^{-1} = \begin{pmatrix}
    		T(x) + 2s T(x) \II(x)T(x) + O(s^2) & O(s) \\
    		O(s) & 1+O(s^2)
    	\end{pmatrix},
    \end{equation}
    where $\II_{ij}(x) = \partial_{ij}\varphi(x)/\sqg$ is the second fundamental form of $\partial\Omega = \Graph\varphi$, and
    \begin{equation}\label{tmp:laOs2}
    	\det DG(x,s) = \sqg \left[ 1-s \Tr(T(x) \II(x)) \right]  
    	+ O(s^2),
    \end{equation}
    where the second term is the mean curvature $\Tr_g \II(x) = \Tr(T(x) \II(x))$ of $\partial\Omega$, when $x$ is sufficiently close to the origin.
	In particular as $s\to 0+$,
	\[ A(x,s) \to A(x,0) := c(x) \begin{pmatrix}
	\left( \Id_{n-1} + \left( \nabla \varphi(x) \right)^T  \nabla \varphi(x) \right)^{-1} & 0 \\
	0 & 1
\end{pmatrix}, \]
where 
\[ c(x) = \lim_{s\to 0+} \det DG(x,s) =  \frac{1+|\nabla \varphi(x)|^2 }{ \sqrt{1+|\nabla \varphi(x)|^2}} = \sqrt{1+|\nabla \varphi(x)|^2} \] 
is a strictly positive scalar function.
    By writing the symmetric matrix $A(x,s)$ as a block matrix
\[ A(x,s) = \begin{pmatrix}
	A_0(x,s) & B(x,s) \\
	B(x,s)^T & d(x,s)
\end{pmatrix}, \]
where $A_0(x,s)$ denotes an $(n-1)\times (n-1)$ matrix,
we define a symmetric matrix $\tilde{A}(x,s)$ on the entire ball $B_{2r}(0)$ as
\begin{equation}\label{def:parity}
	\tilde{A}(x,s) = \left\{\begin{array}{ll}
	A(x,s), & s > 0 \\
	A(x,0), & s=0 \\
	\begin{pmatrix}
		A_0(x,-s) & -B(x,-s) \\
		-B(x,-s)^T & d(x,-s) \\
	\end{pmatrix}, & s<0.
\end{array} \right. 
\end{equation} 
Note that $\tilde{A}(x,s)$ is $C^{k-2}$ in the tangential direction, i.e.\ in the $x$-variable. Away from $s=0$, each of the entries of $\tilde{A}(x,s)$ is $C^\infty$ in the normal direction (i.e.\ in the $s$-variable). However, by the expression \eqref{tmp:laOs}, \eqref{tmp:laOs2}, the matrix $\tilde{A}(x,s)$ is only Lipschitz in the normal direction at the boundary $s=0$. In fact, no matter how regular $\Omega$ is, as long as it is curved (i.e.\ $\II(0)$ is non-trivial), $\tilde{A}(x,s)$ is not $C^1$ across the boundary $s=0$.

    The odd reflection of $v$ is defined by
    \[ \tilde{v}(x,s) = \left\{\begin{array}{ll}
	v(x,s), & \text{ if } s\geq 0 \\
	-v(x,-s), & \text{ if } s<0.
\end{array}
 \right. \]
 By computation, $\tilde{v}$ satisfies $-\divg(\tilde{A}(x,s)\nabla \tilde{v}) = \lambda \left|\det DG(x,|s|)\right| \, \tilde{v} $ in $B_{2r}(0)$.\footnote{Since $\lim_{s\to 0\pm} \tilde{A}(x,s) = A(x,0)$ has zeros on the \textit{off-diagonal} block matrix, a computation shows that the conormal derivatives of $\tilde{v}(y,s)$ from above (i.e.\ $\bR^n_+$) and below (i.e.\ $\bR^n_-$) cancel each other out.}
 
 Now we pull back $\tilde{v}$ in the lower half space $\bR^n_-$ to $\Omega^c$. Since $G$ is a diffeomorphism in a neighborhood of the origin, we may define $\tilde{u}$ by $\tilde{v}(x,s) = \tilde{u}\circ G(x,s)$. Clearly $\tilde{u} = u$ in $\overline{\Omega}$, and this is an extension of $u$. Moreover $\tilde{u}$ satisfies
 \[ -\divg(B\nabla \tilde{u}) = \lambda \dfrac{\left|\det DG(x,|s|)\right| }{\left|\det DG(x,s)\right| } \tilde{u}, \quad \text{ in } B_{2r}(0), \]
 where the matrix $B$ is given by
 \[ B(G(x,s)) = \left|\det DG(x,s)\right|^{-1} (DG(x,s)) \tilde{A}(x,s) (DG(x,s))^T. \] 
 Clearly by \eqref{eq:Adt}, we have $B=\Id$ in $\overline{\Omega}$, as is expected. We cannot further simplify the expression for the matrix $B$ because when we map $s\mapsto -s$, due to the curvature of $\partial\Omega$ the matrix $DG(x,s)^T DG(x,s)$ does not have the same parity as in the definition \eqref{def:parity}, and $\left|\det DG(x,s)\right| $ is not necessarily even with respect to $s$. But we remark that since $\tilde{A}(x,s)$ is Lipschitz, we have that the matrix $B(y,t)$ is also Lipschitz. Indeed, the coefficient on the right hand side $\left|\det DG(x,|s|)\right|/\left|\det DG(x,s)\right|$ is also Lipschitz. Also, away from $s=0$, both $\tilde{A}$ and $B(G(x,s))$ have $C^{k-2}$ regularity.
\end{proof}

\emph{Step two.} By Lemma \ref{lm:ext}, $u$ can be extended to solve an interior elliptic problem~\eqref{eq:PDEext} with Lipschitz coefficients. We expect similar ideas as in \cite{Bers} to give the expansion of $\tilde{u}$ into a homogeneous harmonic polynomial (Proposition \ref{prop:bdexp}). However, we cannot directly appeal to \cite{Bers} since that result requires a non-divergence form operator with H\"older coefficients, whereas \eqref{eq:PDEext} is of divergence form and the leading-order coefficient matrix is only Lipschitz, not $C^{1,\alpha}$. Instead, we now prove the expansion \eqref{eq:expansion} in Proposition \ref{prop:bdexp} by hand, using ideas similar to \cite{KZexpansion}. For simplicity of the notation, when there is no confusion we still denote the extension $\tilde{u}$ as $u$. We also assume $p=0$ without loss of generality.

\begin{proof}[Proof of ~\eqref{eq:expansion}]
We may assume that $u$ does not vanish identically. Lemma \ref{lm:ext} says that $u$ solves~\eqref{eq:PDEext} with Lipschitz coefficients, so standard unique continuation tells us that $u$ vanishes to finite order at each zero. That is, there exists $k\in \mathbb{N}$, such that
\begin{equation}\label{eq:uvo}
	|u(Y)| \leq C|Y|^k, \quad \text{ and } \quad \lim_{s\to 0 } \frac{\sup_{B_s(0)} |u(Y)|}{s^{\ell}} = +\infty \text{ for any } \ell>k.
\end{equation}
By \eqref{eq:PDEext} we have
\begin{equation}
	-\Delta u = \divg((B-\Id)\nabla u) + \lambda c u = \divg(\vec{f}) + g,
\end{equation}
where we denote
\[ \vec{f} := (B - \Id) \nabla u, \qquad g := \lambda c u. \]
By the standard elliptic estimate (applied to the equation in \eqref{eq:PDEext}), for any $Y\in B_{7r/4}(0) \setminus\{0\}$ we have
\begin{equation}\label{eq:gradu}
	|\nabla u(Y)| 
	\lesssim \frac{1}{s} \left(\frac{1}{s^n} \iint_{B_{s/14}(Y)} u^2 ~dZ \right)^{1/2} \lesssim \frac{1}{s} \sup_{B_{s/14}(Y)} |u|,
\end{equation}
where we denote $s=|Y|$.
Combined with \eqref{eq:uvo} and the Lipschitz regularity of the matrix $B$, we get
\begin{equation}\label{est:vecf}
	|\vec{f}(Y)| \leq \left| B(Y) - \Id \right| \cdot |\nabla u(Y)| \lesssim \frac{|Y|}{s} \cdot \sup_{B_{s/14}(Y)} |u| \lesssim s^k.
\end{equation} 
Moreover, since $c$ is bounded, we have
\begin{equation}\label{est:g}
	|g(Y)| \lesssim |u(Y)| \lesssim s^k,
\end{equation}
with a constant depending on the value of $\lambda$.

Let $\zeta$ be a smooth cut-off function, such that $0\leq \zeta \leq 1$, $\zeta \equiv 1$ on $B_{3r/2}(0)$, and $\zeta$ is compactly supported in $B_{7r/4}(0)$. Let \[ \Gamma(\xi) = \left\{\begin{array}{ll}
    c_n |\xi|^{2-n}, & \text{ if } n \geq 3 \\
    -\frac{1}{2\pi} \log|\xi|, & \text{ if } n=2
\end{array}  \right. \] 
be the fundamental solution of the Laplacian in $\bR^n$. 
Define
\begin{equation}\label{def:w}
	w_1(Y) : = \iint \Gamma(Y-Z) (g\zeta)(Z) \, dZ, \quad \text{ and } \quad
	w_2(Y):= \iint \Gamma(Y-Z) \divg(\vec{f} \zeta)(Z) ~dZ. 
\end{equation} 
By the divergence theorem, we have 
\begin{align*}
	w_2(Y)= \iint \Gamma(Y-Z) \divg(\vec{f} \zeta)(Z) ~dZ & = -\iint \nabla_Z \left( \Gamma(Y-Z) \right) \cdot \vec{f}\zeta(Z) ~dZ \\
	& = \iint \nabla \Gamma(Y-Z) \cdot \vec{f}\zeta(Z) ~dZ.
\end{align*} 
By considering the above integral in the regions $\{|Z|< 2|Y|\}$ and $\{2|Y| \leq |Z| < 2r\}$, one can show it is well-defined, and hence $w_2(Y)$ is well-defined.
Similarly, the integral defining $w_1(Y)$ is also well-defined.
Moreover, note that
\begin{equation}
	-\Delta(w_1+w_2) = g\zeta + \divg(\vec{f}\zeta) = -\Delta u, \quad \text{ in } B_{3r/2}(0).
\end{equation}
Hence $u-w_1-w_2$ is a harmonic function in $B_{3r/2}(0)$, which can be written as the infinite sum of homogeneous harmonic polynomials.
In particular, we have
\begin{equation}\label{eq:P3exp}
	(u -w_1-w_2)(Y) = P_3(Y) + \psi_3(Y),
\end{equation} 
where $P_3$ is a harmonic polynomial of degree at most $k$, and the function $\psi_3$ satisfies $|\psi_3(Y)| \lesssim|Y|^{k+1}$.

Next, we aim to write the expansion for $w_1$ and $w_2$. Since $\nabla \Gamma(\cdot)$ is analytic away from the origin, we may consider the Taylor expansion of $\nabla \Gamma(\cdot - Z)$ at the origin for $Z \in \bR^n \setminus \{0\}$. Let $\alpha = (\alpha_1, \cdots, \alpha_n)$ denote an $n$-index. For each $j\in \{0, \cdots, k\}$, we define an $\bR^n$-valued function as follows
\[ \wGk(Y,Z) := \sum_{|\alpha|=j} D^\alpha \nabla \Gamma(-Z) \frac{Y^\alpha}{\alpha!}. \]
For fixed $Z \in \bR^n \setminus \{0\}$, the function $\wGk(\cdot, Z)$ is a harmonic homogeneous polynomial of degree $j$.
Then, since
\[ |D^\alpha \nabla \Gamma(-Z)| \lesssim |Z|^{1-n-|\alpha|}, \]
we have
\begin{equation}\label{est:wGk}
	\left| \wGk(Y,Z) \right| \leq C_j |Z|^{1-n-j}|Y|^j,
\end{equation}
where the constant $C_j$ depends on $j$ as well as the dimension $n$.

Let
\begin{equation}\label{def:P2integral}
	P_2(Y) := \iint  \sum_{j=0}^k \wGk(Y,Z) \cdot \vec{f}\zeta(Z) ~dZ. 
\end{equation} 
Since $\wGk$ is not defined at $Z=0$, we first need to justify that the integral in~\eqref{def:P2integral} is well-defined.
In fact, for any $\delta \in (0, 2r)$, let
\[ f_\delta(Y):=  \iint_{\{ \delta\leq |Z|<2r\}}  \sum_{j=0}^k \wGk(Y,Z) \cdot \vec{f}\zeta(Z) ~dZ.\]
By \eqref{est:wGk} and \eqref{est:vecf}, we have
\begin{align*}
	|f_\delta(Y)| & \leq \sum_{j=0}^k \iint_{\{ \delta\leq |Z|<2r\}}  \left| \wGk(Y,Z) \right| \cdot \left| \vec{f}(Z) \right| ~dZ \\
	& \lesssim \sum_{j=0}^k |Y|^j \int_{\delta}^{2r} \tau^{k-j} ~d\tau \\
	& \lesssim \sum_{j=0}^{k} |Y|^j r^{k-j+1}, 
\end{align*}
which is uniformly bounded as $\delta \to 0$. 
Moreover, let $\gamma = (\gamma_1, \cdots, \gamma_n)$ be an $n$-index, such that $|\gamma| = i \in \{0, 1, \cdots, k \} $. Notice that when we take the $Y$-derivative of $\wGk$, it does not affect the coefficients which just depend on $Z$. Thus, we similarly obtain
\[ |D^\gamma f_\delta(Y)| \lesssim \sum_{j=i}^{k} |Y|^{j-i} r^{k-j+1},   
\]
which is also uniformly bounded as $\delta \to 0$. 

Since $f_\delta(Y)$ is a polynomial of degree at most $k$, it is completely determined by $D^\gamma f_\delta(0)$ with indices $|\gamma| \in \{0, 1, \cdots, k\}$.
 Therefore as $\delta \to 0$ (modulo passing to a subsequence), the sequence $f_\delta(Y)$ converges to $P_2(Y)$ in $C^m_{\operatorname{loc}}(\bR^n)$, for any $m\in \mathbb{N}$, and $P_2$ is well-defined. Moreover, since $f_\delta(Y)$ is a harmonic function for any $\delta>0$, the limit function $P_2(Y)$ is a harmonic polynomial of degree less than or equal to $k$.

We will estimate the remaining term
\begin{equation}\label{eq:P2exp}
	\psi_2(Y):= w_2(Y) - P_2(Y) = \iint \left( \nabla \Gamma(Y-Z) - \sum_{j=0}^k \wGk(Y,Z) \right) \cdot \vec{f}\zeta(Z) ~dZ. 
\end{equation} 
Denoting $s=|Y|< 2r$, we split the integral in \eqref{eq:P2exp} into three parts:
\[ \I := \iint_{\{|Z|< 2s\}} \nabla \Gamma(Y-Z) \cdot \vec{f}\zeta(Z) ~dZ, \]
\[ \II := \iint_{\{ |Z|< 2s\}} \sum_{j=0}^k \wGk(Y,Z) \cdot \vec{f}\zeta(Z) ~dZ, \]
\[ \III := \iint_{\{2s \leq |Z| <2r\}} \left( \nabla \Gamma(Y-Z) - \sum_{j=0}^k \wGk(Y,Z) \right) \cdot \vec{f}\zeta(Z) ~dZ. \]
By \eqref{est:vecf} and the bound on the fundamental solution $\Gamma$, we can easily estimate
\begin{align}
	|\I| & \lesssim \iint_{\{|Z|< 2s\}} |Y-Z|^{1-n} \cdot |\vec{f}(Z)| ~dZ \nonumber \\
	& \lesssim \iint_{\{|Z|< 2s\}} |Y-Z|^{1-n} \cdot |Z|^k ~dZ \nonumber \\
	& \lesssim s^k \iint_{\{|X|< 3s\}} |X|^{1-n} ~dX \nonumber \\
	& \lesssim s^{k+1}. \label{er:Gamma1}
\end{align}
Combining \eqref{est:wGk} and \eqref{est:vecf}, we get
\begin{align}
	|\II| & \lesssim \sum_{j=0}^k s^j  \cdot  \iint_{\{ |Z| < 2s \}} |\vec{f}(Z)||Z|^{1-n-j} ~dZ \nonumber \\
	& \lesssim \sum_{j=0}^k s^j \int_0^{2s} \tau^{1-n-j} \cdot \tau^k \cdot \tau^{n-1} ~d\tau \nonumber  \\
	& \lesssim \sum_{j=0}^k s^j \cdot s^{k-j+1} \nonumber  \\
	& \lesssim_k s^{k+1}.
	\label{er:Gamma2}
\end{align}
Lastly, since $\nabla \Gamma(\cdot)$ is analytic away from the origin, on the set $\{|Z| \geq 2s\}$ we have the expansion
\[ \nabla \Gamma(Y-Z) - \sum_{j=0}^k \wGk(Y,Z) = \sum_{|\alpha|=k+1} D^\alpha \nabla \Gamma(\theta Y-Z) \frac{Y^\alpha}{\alpha!} \]
for some $\theta \in [0,1]$ depending on $Z$ and the component.
Hence by the decay of the fundamental solution, we have
\begin{equation}\label{est:wGkr}
	\left| \nabla \Gamma(Y-Z) - \sum_{j=0}^k \wGk(Y,Z) \right| \leq \sum_{|\alpha|=k+1} |\theta Y - Z|^{1-n-|\alpha|} \cdot \frac{s^{|\alpha|}}{\alpha!} \lesssim_k \frac{s^{k+1}}{|Z|^{n+k}}.
\end{equation}
Therefore
\begin{align}
	|\III| & \lesssim \iint_{\{2s \leq |Z| < 2r\}} \frac{s^{k+1}}{|Z|^{n+k}} \cdot |\vec{f}(Z)| ~dZ \nonumber \\
	& \lesssim s^{k+1} \int_{2s}^{2r} \frac{1}{\tau^{n+k}} \cdot \tau^k \cdot \tau^{n-1} ~d\tau \nonumber \\
	& \lesssim s^{k+1} \log \frac{2r}{2s} \lesssim s^{k+1} \log \frac{1}{s} \lesssim s^{k+\beta}, \label{er:Gamma3}
\end{align}
which holds for an arbitrary choice of $\beta \in (0,1)$, as long as $r$ is chosen sufficiently small (depending on the value of $\beta$).
Combining \eqref{eq:P2exp}, \eqref{er:Gamma1}, \eqref{er:Gamma2}, and \eqref{er:Gamma3}, we conclude that
\begin{equation}\label{eq:errordecay}
	\left|\psi_2(Y) \right| \lesssim s^{k+\beta}, \qquad \text{ for any } 0<\beta<1. 
\end{equation}

Given \eqref{est:g}, the analysis for $w_1(Y)$ is analogous to that for $w_2(Y)$, and the only difference is the homogeneities of $\Gamma$ and $\nabla \Gamma$. In fact, since $\Gamma(\cdot) $ is analytic away from the origin, for $Z\in \bR^n \setminus \{0\}$ fixed we have the Taylor expansion
\begin{equation}
	\Gamma(Y-Z) = \sum_{j=0}^k \sum_{|\alpha|=j} D^\alpha \Gamma(-Z) \frac{Y^\alpha}{\alpha!} + \sum_{|\alpha|=k+1} D^\alpha \Gamma(\theta Y-Z) \frac{Y^\alpha}{\alpha!},
\end{equation}
for some $\theta\in [0,1]$ depending on $Z$.
By the decay of the fundamental solution, we have
\begin{equation}\label{est:wGk1}
	\left| \sum_{|\alpha|=j} D^\alpha \Gamma(-Z) \frac{Y^\alpha}{\alpha!} \right| \leq C_j |Z|^{2-n-j}|Y|^j, 
\end{equation}
which is analogous to the estimate in \eqref{est:wGk} but has a one-degree higher decay rate in $|Z|$. The only exception is when $n=2$ and $j=0$, when we have a logarithmic control $|\Gamma(-Z)| \leq C|\log|Z||$.
Moreover, when $|Z|\geq 2|Y|$ we have
\begin{equation}\label{est:wGkr1}
	\left| \sum_{|\alpha|=k+1} D^\alpha \Gamma(\theta Y-Z) \frac{Y^\alpha}{\alpha!} \right| \lesssim \frac{|Y|^{k+1}}{|Z|^{n+k-1}},
\end{equation}
which is analogous to the estimate in \eqref{est:wGkr} but has a one-degree higher decay rate in $|Z|$.

Let
\begin{equation}
	P_1(Y) := \iint \left[ \sum_{j=0}^k \sum_{|\alpha|=j} D^\alpha \Gamma(-Z) \frac{Y^\alpha}{\alpha!} \right] (g\zeta)(Z) \, dZ, 
\end{equation}
which is a harmonic polynomial of degree less than or equal to $k$, and let
\begin{equation}\label{eq:P1exp}
	\psi_1(Y) := w_1(Y) - P_1(Y) =\iint \left[\Gamma(Y-Z) - \sum_{j=0}^k \sum_{|\alpha|=j} D^\alpha \Gamma(-Z) \frac{Y^\alpha}{\alpha!} \right] (g\zeta)(Z) \, dZ . 
\end{equation} 
We can estimate $\psi_1(Y)$ by splitting it into three integrals as above, and using the estimates \eqref{est:wGk1}, \eqref{est:wGkr1} to show that
\[ |\I|, |\II| \lesssim |Y|^{k+2}, \quad \text{ and } \quad |\III| \lesssim |Y|^{k+1}. \]
Therefore $|\psi_1(Y)| \lesssim |Y|^{k+1}$.

Finally, combining \eqref{eq:P3exp} and \eqref{eq:P2exp}, \eqref{eq:P1exp}, we have the following expansion:
\begin{equation}\label{tmp:exp}
	u(Y) = P_1(Y) + P_2(Y) + P_3(Y) +  \psi_1(Y) + \psi_2(Y) + \psi_3(Y), \quad \text{ in } B_{3r/2}(0), 
\end{equation} 
where $P_1 + P_2+P_3 $ is a harmonic polynomial of degree less than or equal to $k$, and $\left| \psi_1(Y) \right|, |\psi_3(Y)| \leq C |Y|^{k+1}$, $\left| \psi_2(Y) \right| \leq C |Y|^{k+\beta}$ for any $0<\beta<1$.
Because $k$ is the vanishing order of $u$ (see \eqref{eq:uvo}), it follows that $P_1 + P_2 + P_3$ is non-trivial and homogeneous of degree exactly $k$.
This finishes the proof of the expansion \eqref{eq:expansion}, with $P=P_1+P_2+P_3$ and $\psi = \psi_1 + \psi_2 + \psi_3$.
\end{proof}

\emph{Step three.} We aim to prove the gradient estimates \eqref{est:grad1}, \eqref{est:grad2} for the remaining term $\psi$. Since $\tilde{u}$ solves the equation \eqref{eq:PDEext} in the ball $B_{3r/2}(0)$ with only Lipschitz coefficients, we have that $\tilde{u}\in C^{1,\alpha}$; but in general $\tilde{u}$ is not necessarily $C^2$ in the entire ball (see Remark \ref{rmk:reg}). So we only expect to estimate the second-order derivatives $\nabla^2 \psi$ inside the domain (as well as in $\overline{\Omega}^c$), but there may be a jump discontinuity across the boundary. 

Note that $\psi$ satisfies the following inhomogeneous boundary value problem
\begin{equation}\label{eq:psi}
	\left\{\begin{array}{ll}
		-\Delta \psi - \lambda \psi = \lambda P, & \text{ in } \Omega \cap B_{3r/2}(0), \\
		\psi = -P, & \text{ on } \partial\Omega \cap B_{3r/2}(0).
	\end{array} \right.
\end{equation}
Let $Y\in B_{r/2}(0)\cap \Omega$ be fixed. Let $s :=|Y|< r/2$, and $d :=\dist(Y,\partial\Omega)$. Since $0\in \partial\Omega$, clearly 
\[ d = \dist(Y,\partial\Omega) \leq |Y| = s. \]
We consider two cases: either $d > r/4$ and we use interior Schauder estimate, or $d\leq r/4$ and we use boundary Schauder estimate for $C^{2,\alpha}$ domains. When $d>r/4$, we have that
\[ s\geq d>\frac{r}{4} > \frac{s}{2}. \]
Also note that
\begin{equation}\label{eq:incl}
	B_{\frac{d}{3}}(Y) \Subset B_{\frac{d}{2}}(Y) \subset B_{\frac32 s}(0) \subset B_{\frac32 r}(0). 
\end{equation} 
Defining $\psi_d(Z):= \psi(Y+dZ)$, it solves the elliptic equation
\begin{equation}\label{eq:intSch}
	-\Delta \psi_d - \lambda d^2 \psi_d = \lambda d^2 P(Y+d\cdot), \qquad \text{ in } B_{1/2}(0). 
\end{equation} 
Denote
\[ \|f\|_{\dot{C}^{0,\alpha}(D)} = \sup_{Z, Z' \in D} \frac{|f(Z) - f(Z')| }{|Z-Z'|^{\alpha}}, \quad \|f\|_{C^{j,\alpha}(D)} = \sum_{i=0}^j \|\nabla^i f\|_{L^\infty(D)} + \|\nabla^j f\|_{\dot{C}^{0,\alpha}(D)}.  \]
The interior Schauder estimate says that
\begin{align*}
	\|\psi_d\|_{C^{2,\alpha}(B_{1/3}(0)) } & \leq C\left( \|\psi_d\|_{L^\infty(B_{1/2}(0)) } + \|\lambda d^2 P(Y+d\cdot) \|_{C^{0,\alpha}(B_{1/2}(0)) } \right) \\
	& \leq C\left( \|\psi\|_{L^\infty(B_{d/2}(Y)) } + \lambda d^2 \|P \|_{L^\infty(B_{d/2}(Y)) } + \lambda d^{2+\alpha} \|P\|_{\dot{C}^{0,\alpha}(B_{d/2}(Y))} \right) \\
	& \lesssim s^{k+\beta} + d^{2+k} \\
	& \lesssim s^{k+\beta},
\end{align*} 
by \eqref{eq:incl} and $d\leq s$.
We remark that the constant $C$ above depends on the upper bound of the zero-th order coefficient in the equation \eqref{eq:intSch}, which is $\lambda d^2 \leq \lambda r^2/4$.
Since $d>s/2$, the gradient estimates \eqref{est:grad1}, \eqref{est:grad2} at $Y$ follow.

Next assume that $d\leq r/4$. Suppose $q\in \partial\Omega$ such that
\begin{equation}
	|Y-q| = \dist(Y, \partial\Omega)  = d.
\end{equation}
Recall $d\leq s< r/2$, and note that $|q| \leq d+s$. We have
\begin{equation}
	Y \in \overline{B_s(q)} \subset B_{\frac32 s}(q) \subset B_{d+\frac52 s}(0) \subset B_{\frac32 r}(0).
\end{equation}
Defining $\psi_s(Z) := \psi(q+sZ)$, it solves the inhomogeneous boundary value problem
\[ \left\{\begin{array}{ll}
	-\Delta\psi_s - \lambda s^2 \psi_s = \lambda s^2 P(q+s\cdot), & \text{ in } B_{\frac32}(0) \cap \frac{\Omega-q}{s}, \\
	\psi_s = -P(q+s\cdot), & \text{ on } B_{\frac32}(0) \cap \frac{\partial\Omega-q}{s}.
\end{array} \right. \]
Thus applying the boundary Schauder estimate (see for example \cite[Theorem 7.3]{ADN}) to the $C^{2,\alpha}$ domain $\Omega_s := (\Omega-q)/s$, we have
\begin{align}
	\|\psi_s\|_{C^{2,\alpha}(B_1(0) \cap \Omega_s)} & \leq C \left( \|\psi_s\|_{L^\infty(B_{3/2}(0) \cap \Omega_s)} +  \|\lambda s^2 P(q+s\cdot)\|_{C^{0,\alpha}(B_{3/2}(0) \cap \Omega_s)} \right. \nonumber \\
	& \qquad\qquad \left. + \|P(q+s\cdot)\|_{C^{2,\alpha}(B_{3/2}(0) \cap \partial\Omega_s ) } \right).\label{eq:bSch} 
\end{align}
Since $|q| \leq d+s \leq 2s$, the first two terms above are bounded by $s^{k+\beta}$. Naively, the final term involving the boundary data is bounded by $s^{k}$, but this is not good enough. On the other hand by the expansion \eqref{eq:expansion}, we know that (a constant multiple of) $P$ is the tangent function of $u$ at the origin. Moreover, since $u$ vanishes on $\partial\Omega$, it follows that its tangent function $P$ vanishes on the hyperplane
\[ \lim_{\tau\to 0+} \frac{\partial\Omega}{\tau} = \bR^{n-1} = \{(x,0): x\in \mathbb{R}^{n-1}\}. \]
Hence $P$ and all its tangential derivatives on $\bR^{n-1}$ vanish. We use this added information to give a better estimate of $P$ on $\partial\Omega$. 
For any $(x,\varphi(x)) \in \partial\Omega \cap B_{3s/2}(q)$, we have
\begin{equation}\label{tmp:bSch1}
	|P(x,\varphi(x))| = |P(x,\varphi(x)) - P(x,0)| = |\varphi(x)| \cdot| \partial_n P(x, t_1 \varphi(x))|
\end{equation}
for some $t_1 \in [0,1]$. Denoting $e= x/|x| \in \mathbb{R}^{n-1}$, we have
\begin{equation}\label{tmp:bSch2}
	|\varphi(x)|  = \left|\varphi(0) + \partial_e \varphi(0) |x| + \partial_e \varphi(t_2x) \frac{|x|^2}{2} \right| = |\partial_e \varphi(t_2x)| \cdot \frac{|x|^2}{2} \lesssim |x|^3 \lesssim s^3. 
\end{equation} 
And standard elliptic estimates imply that
\begin{equation}\label{tmp:bSch3}
	| \partial_n P(x, t_1 \varphi(x))| \lesssim |(x,\varphi(x))|^{k-1} \lesssim s^{k-1}. 
\end{equation} 
Combining \eqref{tmp:bSch1}, \eqref{tmp:bSch2} and \eqref{tmp:bSch3}, we obtain that
\begin{equation}\label{eq:bSch0}
	\|P\|_{L^\infty(B_{3s/2}(q) \cap \partial\Omega) } \lesssim s^{k+2}. 
\end{equation} 
Furthermore, let $\nabla_{\sslash}^\gamma$ denote the tangential derivative on $\partial\Omega_s$ or $\partial\Omega$, with index $\gamma = (\gamma_1, \cdots, \gamma_{n-1})$ and $|\gamma| = 1, 2, 3$. Then again by comparing the tangential derivatives of $P$ at $(x, \varphi(x))$ and $(x,0)$, we have
\begin{equation}\label{eq:bSch2}
	\|\nabla_{\sslash}^\gamma \left(P(q+s\cdot) \right) \|_{L^\infty(B_{3/2}(0) \cap \partial\Omega_s)} = s^{|\gamma|} \cdot \|\nabla_{\sslash}^\gamma P \|_{L^\infty(B_{3s/2}(q) \cap \partial\Omega)} \lesssim s^{|\gamma|} \cdot s^3 \cdot s^{k-1-|\gamma|} = s^{k+2}.
\end{equation}
Combining \eqref{eq:bSch0} and \eqref{eq:bSch2}, we get
\[ \|P(q+s\cdot)\|_{C^{2,\alpha}(B_{3/2}(0) \cap \partial\Omega_s )} \lesssim s^{k+2}. \]
Finally, the desired gradient estimates \eqref{est:grad1}, \eqref{est:grad2} at the point $Y$ follow from \eqref{eq:bSch}.

To estimate the derivatives of $\psi$ on $\overline{\Omega}^c$, we observe that it satisfies the following inhomogeneous boundary value problem
\begin{equation}
	\left\{\begin{array}{ll}
		-\divg(B \nabla \psi) - \lambda c \psi = \lambda c P + \divg((B-\Id) \nabla P), & \text{ in } \overline{\Omega}^c \cap B_{3r/2}(0), \\
		\psi = -P, & \text{ on } \partial\Omega \cap B_{3r/2}(0).
	\end{array} \right.
\end{equation}
This is more complicated than the equation in $\Omega$, namely \eqref{eq:psi}. However, we note that inside $\overline{\Omega}^c$, the coefficient matrix $B$ and the function $c$ have $C^{1,\alpha}$ regularity, when the domain $\Omega$ is assumed to be $C^{2,\alpha}$. Together with the quantitative $C^{2,\alpha}$ estimates of the harmonic polynomial $P$, this allows us to apply the interior and boundary Schauder estimates in a similar fashion as in $\Omega$, and obtain the desired gradient estimates for $\psi$ in $\overline{\Omega}^c$.

Lastly, since $u\in C^{1,\alpha}$ in the entire ball $B_r(0)$, the estimate of $\nabla \psi$ holds everywhere, including on the boundary $\partial\Omega$. However, $u$, and hence $\psi$, is not expected to be $C^2$ across the boundary. The estimate of $\nabla^2\psi$ only means the one-sided derivatives as we approach the boundary from $\Omega$, and they may well disagree with the one-sided derivatives as we approach the boundary from $\Omega^c$.

\section{Proof of Main Results}\label{sec:main-proof}

\begin{prop}\label{prop:main}
    Assume $\Omega$ is a bounded open $C^{2,\alpha}$ domain (resp.\ $C^k$ domain) in $\bR^n$. 
Let $u$ be a Dirichlet Laplacian eigenfunction of $\Omega$, such that
\begin{enumerate}[label=(\roman*)]
    \item\label{cond:interior-regular} Zero is a regular value of $u$ on $\Omega$,
    \item\label{cond:boundary-regular} Zero is a regular value of $\dd_{\nu} u: \dd \Omega \to \bR$, the normal derivative of $u$ on the boundary.
\end{enumerate}
Then $N_u := \overline{\{u=0\}\cap \Omega}$ is a compact, possibly disconnected, embedded hypersurface which intersects $\partial\Omega$ orthogonally, and its manifold boundary $\partial N_u$, if not empty, is a $C^{2,\alpha}$-regular (resp. $C^{k-1}$-regular) embedded submanifold of $\partial\Omega$.
\end{prop}

In Theorem \ref{thm:generic}, we have shown that any generic connected compact $C^k$ domain, $k\geq n$, satisfies the above conditions \ref{cond:interior-regular} and \ref{cond:boundary-regular}.
So Theorem \ref{thm:main} immediately follows from Proposition \ref{prop:main}. It remains to prove the proposition.

Recall that Proposition \ref{prop:bdexp} implies that for any $p\in \partial\Omega$, $u$ can be locally expanded into the sum of a homogeneous harmonic polynomial $P$ with degree $k\in \mathbb{N}$ and an error term $\psi$. Moreover, $P$ vanishes on $\Tan_p \partial\Omega$, the tangent hyperplane to $\partial\Omega$ at $p$. Modulo translations and rotations, we may assume the point $p=0$, the tangent hyperplane $\Tan_p \partial\Omega = \{x_1=0\} = \{0\} \times \bR^{n-1}$, and the unit outer normal vector $\nu$ at $p$ agrees with $-e_1$. 
Condition \ref{cond:interior-regular} implies that $\{u=0\} \cap \Omega$ is a smooth embedded hypersurface, possibly with boundary and possibly with multiple connected components.
\begin{claim}\label{cl:P2}
    We claim that condition \ref{cond:boundary-regular} implies that one of the following holds:
    \begin{itemize}
        \item  $k=1$ and $P=cx_1$ for some $\bR \ni c \neq 0$,
        \item $k=2$ and, modulo a further rotation, $P=c x_1 x_2$ for some $c>0$.
    \end{itemize}
\end{claim}

\begin{proof}[Proof of Claim]
If $\partial_\nu u(0)\neq 0$, then $\nabla u(0) \neq 0$ and thus $\nabla P(0) = \nabla u(0) \neq 0$ by the gradient estimate \eqref{est:grad1}. Hence $k=1$ and $P$ is a linear function; furthermore, since $P$ vanishes on $\Tan_0 \partial\Omega = \{x_1=0\}$, it must be equal to $cx_1$ for some nonzero real number $c$. Assume $\partial_\nu u(0) = 0$, and thus $\nabla u(0)=0$. By the gradient estimate \eqref{est:grad1} we have $\nabla P(0)=0$, and thus $k\geq 2$. Moreover, condition \ref{cond:boundary-regular} implies that $d(\partial_\nu u):\Tan_0 \partial\Omega \to \bR$ is surjective at $0$, so there exists $e \in \Tan_0 \partial\Omega \setminus\{0\}$ such that $\partial_e \partial_\nu u(0) \neq 0$. Hence it follows from \eqref{est:grad2} that 
\[ \partial_e \partial_\nu P(0) = \partial_e \partial_\nu u(0) \neq 0. \]
In particular $k=2$, i.e.\ $P$ is a homogeneous harmonic polynomial of degree $2$.

Homogeneous harmonic polynomials of degree $2$ are fully classified: They must be of the form
\begin{equation}\label{exp:P2}
	P_2(x) = \langle x, Ax \rangle,
\end{equation}
where $A = (a_{ij})$ is an $n\times n$ real symmetric square matrix with $\Tr A = 0$. 
We shall show that $P_2$ must be the simplest possible degree-$2$ homogeneous harmonic polynomial, that is, $P_2$ is two-dimensional and independent of the other $(n-2)$ variables. If this is not the case, then $\{P_2=0\}$ cannot contain a flat hyperplane, contradicting the smoothness of the boundary $\partial\Omega$. We remark that this observation is not true for homogeneous harmonic polynomials of degree bigger than $2$. For example, $P_3(x,y,z) = xyz$ is a homogeneous harmonic polynomial of degree $3$ in $\bR^3$, and its nodal set consists of the union of three planes.

Recall we have proven in step three of Section \ref{sec:bd} that $P\equiv 0$ on $\Tan_0 \partial\Omega = \{x_1=0\}$. Considering the expression \eqref{exp:P2}, the matrix $A$ must be a simple block matrix of the form
\[ A = \begin{pmatrix}
    a_{11} & a_{12} \,\, \cdots \,\,  a_{1n} \\
    a_{12} & \\
    \vdots & 0 \\
    a_{1n} & \\
\end{pmatrix}. \]
Thus $\Tr A = a_{11}=0$. Modulo a rotation in $\{x_1=0\} = \bR^{n-1}$ (to identify the unit vector in the direction of $\sum_{j=2}^n a_{1j} e_j$ as the new coordinate $e_2$), the homogeneous harmonic polynomial $P$ must be of the form $P(x) = cx_1 x_2$, with $c>0$. This finishes the proof of Claim \ref{cl:P2}.
\end{proof}

An important heuristic is that the nodal set of the homogeneous harmonic polynomial $P$ is close to that of $u$. We make this quantitative in the following lemma, which will be used to study the geometric properties of $\{u=0\}$ in combination with Claim \ref{cl:P2}. 
The proof is very similar to \cite[Lemma 4.8]{HS}.

\begin{lem}\label{lm:HS}
	Suppose $w_1, w_2$ are two functions in $B_{2r}(0)$ satisfying
	\begin{equation}\label{eq:bd}
		\|w_i\|_{L^\infty} \leq A_1 r^2, \quad \|\nabla w_i\|_{L^\infty} \leq A_1 r, \quad \|\nabla^2 w_i \|_{L^\infty} \leq A_1, \quad i\in \{1,2\},
	\end{equation}
	and moreover,
	\begin{equation}\label{eq:diff}
		\|w_1 - w_2\|_{L^\infty} \leq A_2 r^{2+\beta}, \quad \|\nabla w_1 - \nabla w_2\|_{L^\infty} \leq A_2 r^{1+\beta}, 
	\end{equation}
	for some exponent $\beta \in [3/4,1)$ and constants $A_1, A_2>0$.
	Let $\eta = r^{\frac54} \ll r$, and consider the set
	\begin{equation}\label{def:S}
		N_r = \{x\in B_r(0): w_1(x)=0, |\nabla w_1(x)| > \eta \}. 
	\end{equation} 
	If $r$ is sufficiently small, then for any $x\in N_r$ there are a hyperplane $T_x$ passing through $x$, and functions $h_1^x, h_2^x$ defined on $T_x \cap B_{\alpha\eta}(x)$ such that
	\begin{equation}
		\{w_i = 0\} \cap B_{\alpha\eta}(x)  = \{(z, g_i^x(z)): z\in T_x \cap B_{\alpha\eta}(x) \} \cap B_{\alpha\eta}(x) =: \Graph_{T_x} h_i^x
	\end{equation}
	where $\alpha>0$ is a small constant such that $\alpha A_1 < 1/12$.
	In particular, $\{w_2=0\}$ contains a point $x'$ such that $|x'-x|< \eta^{1+\epsilon}$ with $\epsilon = 1/10$.
\end{lem}
\begin{proof}
	By \eqref{eq:bd}, we know that for any $y\in B_{\alpha\eta}(x)$,
	\begin{equation}\label{eq:lbw1}
		|\nabla w_1 (y) | \geq |\nabla w_1(x)| - \|\nabla^2 w_1\|_{L^\infty} \cdot |x-y| > \eta - A_1\cdot \alpha \eta > \frac23 \eta. 
	\end{equation} 
	Moreover by \eqref{eq:diff},
	\begin{equation}\label{eq:lbw2}
		|\nabla w_2(y)| \geq |\nabla w_1(y)| - \|\nabla w_1 - \nabla w_2\|_{L^\infty} > \frac23 \eta - A_2 r^{1+\beta} > \frac{\eta}{2}, 
	\end{equation} 	
	by choosing $r$ (and thus $\eta< r$) sufficiently small depending on the value of $A_2$ (recall that $r = \eta^{\frac45}$ and $\beta \geq 3/4$).
	Let 
	\[ \nu_i(y) = \frac{\nabla w_i(y)}{|\nabla w_i(y)|}, \]
	which is the unit normal vector for the level set $\{w_i = w_i(y)\}$. Then
	\begin{align}
		|\nu_1(y) - \nu_2(y)| = \left|\frac{\nabla w_1(y)}{|\nabla w_1(y)|} - \frac{\nabla w_2(y)}{|\nabla w_2(y)|} \right| \lesssim \frac{|\nabla w_1(y) - \nabla w_2(y)|}{|\nabla w_i(y)|}  \lesssim \frac{\|\nabla w_1 - \nabla w_2\|_{L^\infty}}{\eta} \leq \eta^{\frac15}, \label{eq:diffnu}
	\end{align} 
	and 
	\begin{align}
		|\nu_i(y) - \nu_i(x)| = \left|\frac{\nabla w_i(y)}{|\nabla w_i(y)|} - \frac{\nabla w_i(x)}{|\nabla w_i(x)|} \right| & \lesssim \frac{|\nabla w_i(y) - \nabla w_i(x)|}{|\nabla w_i(x)|} \nonumber \\
		& \leq \frac{\|\nabla^2 w_i\|_{L^\infty} \cdot |y-x|}{\eta} \nonumber \\
		& \leq \alpha A_1 < \frac{1}{12}.\label{eq:oscnu}
	\end{align}
	Furthermore, letting $\epsilon = 1/10$, any $z\in B_{\eta^{1+\epsilon}}(x)$ satisfies
	\begin{align}
		|\nu_i(z) - \nu_i(x)| & \lesssim \frac{|\nabla w_i(z) - \nabla w_i(x)|}{|\nabla w_i(x)|} \leq \frac{\|\nabla^2 w_i\|_{L^\infty} \cdot |z-x|}{\eta} \lesssim \eta^{\epsilon}.\label{eq:oscnusmall}
	\end{align}
	In particular, combining \eqref{eq:diffnu} and \eqref{eq:oscnusmall} we have
	\begin{equation}\label{eq:diffnusmall}
		|\nu_2(z) - \nu_1(x)| \lesssim \eta^{\epsilon}.
	\end{equation}
	
	Let $T_x$ denote the hyperplane passing through $x$ and orthogonal to $\nu_1(x)$.
	By \eqref{eq:lbw1} and \eqref{eq:oscnu}, the nodal set $\{w_1=0\} \cap B_{\alpha\eta}(x)$ is a union of graphs over $T_x$. An elementary argument using \eqref{eq:lbw1} and the mean-value theorem for functions of one variable justifies that it is in fact the graph of a single function $h_1^x: T_x \cap B_{\alpha\eta}(x) \to \bR$, namely
	\[ \{w_1=0\} \cap B_{\alpha\eta}(x) = \Graph_{T_x} h_1^x. \]
	
	On the other hand, we claim that
	\begin{equation}\label{cl:nodalw2}
		\{w_2=0\} \cap B_{\eta^{1+\epsilon}}(x) \neq \emptyset. 
	\end{equation} 
	In fact, 
	\begin{align*}
		w_1(x+\eta^{1+\epsilon} \nu_1(x)) & = w_1(x) + \partial_{\nu_1(x)} w_1(x+t\eta^{1+\epsilon}\nu_1(x)) \cdot \eta^{1+\epsilon} \\
		& \geq 0 + \left(\partial_{\nu_1(x)}w_1(x) - \|\nabla^2 w_1\|_{L^\infty} \cdot  t\eta^{1+\epsilon} \right) \eta^{1+\epsilon} \\
		& \geq \left(|\nabla w_1(x)| - A_1\eta^{1+\epsilon} \right) \eta^{1+\epsilon} > \frac12 \eta^{2+\epsilon} > 0,
	\end{align*} 
	where in the first line we use the mean-value theorem for functions of one variable and $t\in [0,1]$. Thus
	\[ w_2(x+\eta^{1+\epsilon} \nu_1(x)) \geq w_1(x+\eta^{1+\epsilon} \nu_1(x)) - \|w_1 - w_2\|_{L^\infty} > \frac13 \eta^{2+\epsilon}. \]
	Similarly, we have
	\[ w_2(x-\eta^{1+\epsilon}\nu_1(x)) < -\frac13 \eta^{2+\epsilon} < 0. \]
	By the intermediate value theorem, $w_2$ must vanish at some point $x'$ on the line segment connecting $x\pm\eta^{1+\epsilon}\nu_1(x)$, which finishes proving the claim \eqref{cl:nodalw2}. Combining \eqref{eq:diffnusmall} (for $z=x'$) and \eqref{eq:oscnu} and using the same argument as in the proof for $w_1$, we conclude that the nodal set $\{w_2=0\} \cap B_{\alpha\eta}(x)$ is also the graph of a single function $h_2^x: T_x \cap B_{\alpha\eta}(x) \to \bR$, and moreover, the graph passes through the point $x' \in B_{\eta^{1+\epsilon}}(x)$.
\end{proof}

Next, we discuss the geometric implications of the classification in Claim \ref{cl:P2} and prove Proposition \ref{prop:main}.

Let $N_u := \overline{\{u=0\} \cap \Omega}$. Firstly we show that if $N_u \cap \partial\Omega \neq \emptyset$, then
\begin{equation}\label{cl:Nu-bd}
    N_u \cap \partial\Omega = \{ \partial_\nu u=0\}. 
\end{equation} 
Without loss of generality assume the origin is an arbitrary boundary point. If $\partial_\nu u(0) \neq 0$, by the continuity of $\partial_\nu u$ on $\partial\Omega$ we have either $\partial_\nu u\geq a>0$ or $\partial_\nu u\leq -a <0$ in a neighborhood of the origin. In either case, $u$ cannot vanish in $\Omega$ in a smaller neighborhood of the origin. Thus $0\notin N_u$. 
On the other hand, assume $0\notin N_u$. Then there exists $\epsilon_0>0$ such that $B_{\epsilon_0}(0) \cap \{u=0\} \cap \Omega = \emptyset$. In other words either $u>0$ or $u<0$ in $B_{\epsilon_0}(0) \cap \Omega$. Thus by the boundary Hopf lemma $\partial_\nu u \neq 0$. 
This finishes the proof of \eqref{cl:Nu-bd}.

Clearly we also have $N_u \cap \Omega = \{u=0\} \cap \Omega$, which is a boundaryless embedded hypersurface by the hypothesis \ref{cond:interior-regular}. But we still do not know how this interior hypersurface and the boundary part in \eqref{cl:Nu-bd} are joined together. To that end, we look for a single regular function that describes $N_u$. We focus on the case when $\Omega$ is assumed to be a $C^k$ domain, since it is more relevant to Theorem \ref{thm:main}; and the case when $\Omega$ is a $C^{2,\alpha}$ domain is analogous.

Let $\varphi$ denote the first Dirichlet eigenfunction in $\Om$. Then $\varphi>0$ in $\Omega$, and the boundary Hopf lemma implies that $\partial_\nu \varphi<0$ on $\partial\Omega$. 
We have $(\Delta + \lambda_1) \varphi = 0$ and $(\Delta+\lambda_1) u = (\lambda_1 - \lambda) u \in C^{k-1,\beta}$ in a $C^{k-1,\beta}$-regular domain $\Omega$, with vanishing boundary data. By the higher-order boundary Harnack inequality proven in \cite[Theorem 3.1]{HigherHarnack}, their ratio satisfies $\frac{u}{\varphi} \in C^{k-1,\beta}$ up to the boundary\footnote{Even though \cite[Theorem 3.1]{HigherHarnack} is stated for the interior, the uniform norm bound on $C^{k-1,\beta}(\Omega)$ indicates that we have regularity up to the boundary.}. Moreover, since both $u$ and $\varphi$ are continuously differentiable up to the boundary, its restriction to the boundary satisfies 
\[ \frac{u}{\varphi}\Big|_{\partial\Omega} = \frac{\partial_\nu u}{\partial_\nu \varphi}. 
\] 
We claim that
\[ N_u := \overline{\{u=0\} \cap \Omega} = \left\{x\in \overline{\Omega}: \frac{u}{\varphi}(x) = 0 \right\}.\]
Clearly, on $\Omega$ both sets agree with $\{u=0\}$. On the boundary, since $\partial_\nu \varphi <0$ we have
\[ \left\{x\in \partial \Omega: \frac{u}{\varphi}(x) =\frac{\partial_\nu u}{\partial_\nu \varphi}(x) = 0 \right\} = \{x\in \partial\Omega: \partial_\nu u(x) = 0\} = N_u \cap \partial\Omega \]
by \eqref{cl:Nu-bd}.
Moreover, for any $p=0 \in \{\partial_\nu u=0\}$, the tangential derivative of $\frac{u}{\varphi}$ on $\partial\Omega$ satisfies
\[ d\left(\frac{\partial_\nu u}{\partial_\nu \varphi} \right)(0) = \frac{d(\partial_\nu u)(0)}{\partial_\nu \varphi(p)}, \]
which is non-trivial by the hypothesis \ref{cond:boundary-regular}.
In fact, by Claim \ref{cl:P2} the leading-order homogeneous harmonic polynomial at the origin is, modulo a change of coordinates, $P=cx_1x_2$ for some $c>0$, and its only non-trivial double derivative is the mixed derivative $\partial_{x_2}\partial_{x_1} P = c$. So Proposition \ref{prop:bdexp} implies that 
\[ \partial_{x_2} \partial_{\nu} u(0) = -\partial_{x_2} \partial_{x_1} u(0) = -\partial_{x_2}\partial_{x_1} P \neq 0. \]
Applying the implicit function theorem to $\frac{u}{\varphi} \in C^{k-1,\beta}(\overline{\Omega})$, we obtain that locally near the origin, $N_u$ is the intersection of a $C^{k-1,\beta}$-graph over $\{x_2=0\}$ and $\overline{\Omega}$. Since this holds for any point in $N_u \cap \partial\Omega$, it follows that $N_u$ is a $C^{k-1,\beta}$-embedded hypersurface with boundary, and its manifold boundary satisfies
\[ \partial N_u = N_u \cap \partial\Omega = \{\partial_\nu u = 0\} = \left\{ \frac{u}{\varphi}\Big|_{\partial\Omega} = 0 \right\}, \]
which is a $C^{k-1,\beta}$-embedded submanifold of $\partial\Omega$.

Next, we prove that $N_u$ intersects $\partial\Omega$ orthogonally at the origin. In fact we prove a stronger statement: Assuming the leading-order homogeneous harmonic polynomial at the origin is $P=cx_1 x_2$, we show that $N_u\cap \Omega$ is quantitatively close to $\{P=0\} \setminus \Tan_0 \partial\Omega = \{x_2=0\} \cap \{x_1 \neq 0\}$. 
Let $X = (z, y) \in \bR^2_x \times \bR^{n-2}_y$ be a point in $\{u=0\} \cap \Omega$ that is sufficiently close to the origin. In addition, assume that $X$ has non-trivial projection onto the positive $x_1$-axis, i.e. there exists $\theta_0\in [0,\pi/2)$ depending on $X$ such that
\begin{equation}\label{tmp:Xidist}
    \left\langle \frac{X}{|X|}, e_1 \right\rangle = \cos \theta_0 >0.
\end{equation}
Let $L = \{x_1 = x_2 = 0\}$, and let $s := \dist(X, L) = |z|$. Since $ \langle X ,e_1\rangle = \cos\theta_0 |X|$, it follows that 
\[ s = \sqrt{\langle X ,e_1\rangle^2 + \langle X ,e_2\rangle^2 } \geq |\langle X ,e_1\rangle| \gtrsim |X| > 0. \]
Note that
\[ |\nabla P(X)| = cs, \qquad |\nabla \psi(X)| \lesssim |X|^{1+\beta} \lesssim s^{1+\beta}, \]
and thus
\[ |\nabla u(X)| \geq \frac{cs}{2} =: \eta. \]
Applying Lemma \ref{lm:HS} to the pair of functions $w_1 := u, w_2 := P$ with the parameter $\eta$, we obtain that
\begin{equation}\label{tmp:XidistUB}
	\dist(X, \{P = 0\}) \leq \eta^{1+\epsilon} \approx s^{1+\epsilon}. 
\end{equation} 
On the other hand, the assumption \eqref{tmp:Xidist} implies that
\begin{equation}\label{tmp:XidistLB}
	\dist(X, \{x_1 = 0\}) = |\langle X, e_1 \rangle | \gtrsim |X| \geq s. 
\end{equation} 
Recalling $\{P=0\} = \{x_1=0\} \cup \{x_2=0\}$, \eqref{tmp:XidistUB} and \eqref{tmp:XidistLB} indicate that 
\[ \dist(X, \{x_2 = 0\}) \lesssim s^{1+\epsilon} \approx |X|^{1+\epsilon}. \]
In particular, not only is $N_u$ graphical over the hyperplane $\{x_2=0\}$, we also have $\Tan_0 N_u \subset \{x_2=0\}$; and the tangent space of $N_u \cap \partial\Omega$ at the origin is $L \approx \mathbb{R}^{n-2}_y$.

\section{The two-dimensional case}\label{sec:2d}

We prove Proposition \ref{prop:2D} in this section. Firstly, by the continuous dependence of the Dirichlet eigenvalues $\lambda(\Omega)$ on the domain, it is well-known that the simplicity of eigenvalues is an open condition. Hence small perturbations of $\Omega$ still have simple eigenspaces for eigenvalues in a small neighborhood of $\lambda$. It remains to prove the regularity of the nodal curves. We will focus on the boundary case, since it is more relevant to the goal of the paper and the proof is slightly more involved. An analogous result for the interior has been proven (see \cite[Section 4]{Albert}). For completeness, towards the end of the section, we mention the necessary changes to prove the interior case.

\textbf{Boundary asymptotics near degenerate points in dimension two}.
For any planar domain $\Omega$ and a Dirichlet eigenfunction $u$ on $\Omega$, we say $z \in \partial\Omega$ is a \emph{degenerate boundary point} if 
\begin{equation}
    \partial_\nu u(z) = 0, \qquad \partial_\tau \partial_\nu u(z) = 0,
\end{equation}
where $\nu$ denotes the unit outer normal vector, and $\tau$ denotes the unit tangent vector in $\Tan_z \partial\Omega$. Recalling Proposition \ref{prop:bdexp}, let $P$ be the leading-order homogeneous harmonic polynomial in the expansion of $u$ near $z$, that is,
\[ u(\zeta) = P(\zeta-z) + \psi(\zeta - z), \qquad \text{ in } B_r(z) \cap \Omega. \]
Modulo a rotation, we have $P(\omega) = c \Im \omega^k$ for some $k\in \mathbb{N}$, where the complex variable $\omega = x+iy$ with the positive $x$-axis pointing in the $\tau_z$-direction, and positive $y$-axis pointing in the $(-\nu_z)$-direction. By the error estimates in Proposition \ref{prop:bdexp}, we have
\[ \partial_\nu u(z) = 0 \iff \nabla u(z) = 0 \iff \nabla P(0)= 0 \iff k>1, \]
and additionally,
\[ \partial_\tau \partial_\nu u(z) = 0 \iff \partial_x\partial_y P(0)=0 \iff k>2, \]
where we use $P(\omega) = c\Im \omega^k$, and thus $\partial_y P(\omega) = ck x^{k-1} + y^2 (\cdots)$.
In particular, if $\Om$ has no degenerate boundary point, either $\partial_\nu u$ has a sign, in which case the nodal curve of $u$ does not intersect the boundary; or whenever $\partial_\nu u(z) = 0$, we have $\partial_\tau \partial_\nu u(z) \neq 0$ and thus the leading-order homogeneous harmonic polynomial is quadratic. In the latter case, the nodal curve of $u$ intersects $\partial\Om$ at $z$ orthogonally.

Moreover, the mixed derivative of $u$ can be bounded from below near $z$:

\begin{claim}\label{cl:nondeg-mix-der}
    We have $|\partial_\tau \partial_\nu u(\zeta)| \geq A|\zeta-z|^{k-2}$ for any $\zeta \in B_r(z) \cap \partial\Omega$, where the constant $A$ only depends on $P$ and the constants of the error bounds for $\psi$.
\end{claim}

In fact, by the $C^2$ regularity of $\partial\Omega$ and the bound of $|\nabla^2 u(\zeta)|$ we have
\begin{align*}
    \left| \langle \nabla^2 u(\zeta) \tau_\zeta, -\nu_\zeta \rangle - \langle \nabla^2 u(\zeta) \tau_z, -\nu_z \rangle \right| \lesssim |\nabla^2 u(\zeta)| \cdot |\zeta - z| \lesssim |\zeta-z|^{k-2+1}.
\end{align*}
Since $\tau_z$ is parallel to the positive $x$-axis and $-\nu_z$ is parallel to the positive $y$-axis, we can rewrite 
\[ \langle \nabla^2 u(\zeta) \tau_z, -\nu_z \rangle = \partial_x \partial_y u(\zeta), \] 
and thus
\[ \left| \langle \nabla^2 u(\zeta) \tau_z, -\nu_z \rangle - \partial_x \partial_y P(\zeta-z) \right| = |\partial_x \partial_y \psi(\zeta-z)| \lesssim |\zeta-z|^{k+\beta-2}. \]
By the expression of $P$ and the $C^2$ regularity of $\partial\Omega$, we further have
\[ |\partial_x \partial_y P(\omega)| \gtrsim |\omega|^{k-2}. \]
The combination of the above estimates proves the claim, for sufficiently small radius $r$.

In particular, Claim \ref{cl:nondeg-mix-der} also implies that a degenerate boundary point is isolated. Assume that $z_1, \cdots, z_D$ are all the degenerate boundary points of $\Omega$. 

\textbf{Admissible directions of perturbations}.
Let $s\in C^{2,\alpha}(\overline{\Om}, \mathbb{R}^2)$ denote the admissible direction of perturbations, and let $\rho_s:= s\cdot \nu \in C^{1,\alpha}(\partial\Omega)$. We remark that even though the H\"older space $C^{2,\alpha}(\overline{\Om}, \mathbb{R}^2)$ is not separable, this is not an issue since we do not appeal to Sard-Smale or a transversality theorem.
Assume the eigenfunction $u$ has unit $L^2$ norm. Recall the notation from Section \ref{sec:generic}. In particular, rewriting \eqref{eqn:TQ} at $h=\id$ gives
\begin{equation}
    T_{(u, \la, \id)} Q = \{(v, \eta, s) \in T_u S \times \bR \times C^2(\overline{\Omega}, \bR^n) \,|\, (\Delta + \la) v + \eta u  - (\Delta + \la)(s \cdot \nabla u) =0\}. 
\end{equation}
\begin{lem}\label{lm:unique-triple}
    Assume $\lambda$ is a simple eigenvalue and $u$ is the corresponding eigenfunction. There exists a unique triplet $(v_s, \eta_s, s)\in T_{(u,\lambda, \id)} Q$. Moreover, at the degenerate boundary point $z_i$ we have
    \begin{equation}\label{eq:normder-vs}
        \partial_\nu v_s(z_i) =  - \int_{\partial\Om} \rho_s \,\partial_\nu u \, \partial_\nu P_{z_i},
    \end{equation}
    where $P_{z_i}$ denotes the reduced Poisson kernel for the operator $\Delta+\lambda$ in $\Omega$ anchored at $z_i$.
\end{lem}
\begin{rmk}
    We remark that the above lemma holds in any dimension, not specific to $\bR^2$.
\end{rmk}
\begin{proof}
    The existence of $v_s$ and $\eta_s$ such that $(v_s, \eta_s, s) \in T_{(u,\lambda,\id)} Q$ has been proven in Remark \ref{rmk:simple}, where 
    \begin{equation}\label{tmp:etas}
        \eta_s = -\int_{\partial\Om} s \cdot \nu \, |\partial_\nu u|^2 = -\int_{\partial\Om} \rho_s \, |\partial_\nu u|^2, 
    \end{equation} 
    and $v_s$ is the unique $H^2 \cap H_0^1(\Omega)$ solution that is $L^2(\Om)$-orthogonal to $u$ and satisfies
    \begin{equation}\label{tmp:vs}
        (\Delta + \lambda) v_s = (\Delta+\lambda) (s\cdot \nabla u) - \eta_s u. 
    \end{equation} 
    Moreover, letting $p>2$, by the elliptic estimates and the Sobolev embedding $W^{2,p}(\Om) \hookrightarrow C^1(\Om)$, we have
    \begin{equation}\label{eq:vs-C1}
        \|v_s\|_{C^1(\Om)} \lesssim \|v_s\|_{W^{2,p}(\Om)} \lesssim \|(\Delta+ \lambda) (s\cdot \nabla u) - \eta_s u \|_{L^p(\Om) } \lesssim \|s\|_{C^2}. 
    \end{equation}

    As in Subsection \ref{sec:boundary-reg}, by multiplying \eqref{tmp:vs} by $P_{z_i}$ and integrating, we have
    \begin{equation}\label{tmpeq:zi-1}
        \int_\Om P_{z_i} (\Delta + \lambda)v_s = \int_\Om P_{z_i} (\Delta+ \lambda) (s\cdot \nabla u) - \eta_s \int_{\Om} uP_{z_i}. 
    \end{equation} 
    Recall by Lemma \ref{lem:reduced-poisson-kernel} that $P_{z_i}$ is orthogonal to the eigenfunction $u$ (or the kernel of the operator $L:=\Delta+\lambda$), and thus the last term disappears.    
    Recalling that $\partial_\nu u(z_i) = 0$, the identity for the reduced Poisson kernel implies
    \begin{equation}\label{tmpeq:zi-2}
        \int_\Om P_{z_i} (\Delta + \lambda)v_s = \partial_\nu v_s(z_i);  
    \end{equation} 
    moreover, assuming that $s$ is supported away from $z_i$, we have
    \begin{equation}\label{tmpeq:zi-3}
        \int_\Om P_{z_i} (\Delta+ \lambda) (s\cdot \nabla u) = - \int_{\partial\Om} s\cdot \nabla u \, \partial_\nu P_{z_i} = -\int_{\partial\Om} \rho_s \, \partial_\nu u \, \partial_\nu P_{z_i}.
    \end{equation}
    We claim that \eqref{tmpeq:zi-3} in fact holds without assuming $s$ is supported away from $z_i$, by a density argument and the assumption that $z_i$ is a degenerate boundary point. Therefore the identity \eqref{eq:normder-vs} follows by combining \eqref{tmpeq:zi-1}, \eqref{tmpeq:zi-2} and \eqref{tmpeq:zi-3}.

    To prove the claim, let $s\in C^{2,\alpha}(\overline{\Om}, \mathbb{R}^2)$ be arbitrary, let $\epsilon>0$ be arbitrarily small, and let $\chi_\epsilon$ be a smooth function such that $\chi_\epsilon = 0$ on $B_\epsilon(z_i)$, $\chi_\epsilon = 1$ outside of $B_{2\epsilon}(z_i)$, $|\nabla \chi_\epsilon| \lesssim 1/\epsilon$ and $|\nabla^2 \chi_\epsilon| \lesssim 1/\epsilon^2$. Since $s_\epsilon:= s\chi_\epsilon$ is supported away from $z_i$, \eqref{tmpeq:zi-3} holds, i.e.
    \begin{equation}\label{tmp:Pk-cutoff-nearpole}
        \int_\Om P_{z_i} (\Delta+ \lambda) (s_\epsilon \cdot \nabla u) = - \int_{\partial\Om} \rho_{s_\epsilon} \, \partial_\nu u \, \partial_\nu P_{z_i}. 
    \end{equation} 
    
    To pass to the limit, we first show that $\partial_\nu u \, \partial_\nu P_{z_i}$ is integrable near the pole $z_i$.
    Let $\zeta \in \partial\Om \cap B_{2\epsilon}(z_i) \setminus \{z_i\}$ be fixed, and let $\tau:= |\zeta-z_i|$. Since $P_{z_i}$ satisfies $(\Delta+\lambda) P_{z_i} = 0$ on $\Om \cap B_{\tau/2}(\zeta)$ and $P_{z_i} = 0$ on $\partial\Om \cap B_{\tau/2}(\zeta)$, we may apply the Moser estimate and boundary Schauder estimate to show that
    \begin{equation}\label{eq:Pkbd-near-pole}
        |\partial_\nu P_{z_i}(\zeta) | \lesssim \tau^{-1} \sup_{\Om \cap B_{\tau/4 }(\zeta)} |P_{z_i}| \lesssim \tau^{-1-\frac{n}{q}} \|P_{z_i}\|_{L^q(\Om)} 
    \end{equation} 
    where $1<q<\frac{n}{n-1}$. Recall by the proof of Lemma \ref{lem:reduced-poisson-kernel} that $\|P_{z_i}\|_{L^q(\Om)}$ is uniformly bounded independent of $z_i$. 
    In particular, by choosing $q$ as close to $\frac{n}{n-1}$ as possible, \eqref{eq:Pkbd-near-pole} says that $\partial_\nu P_{z_i}(\zeta)$ decays at the rate $\tau^{-n-\delta }$ for any $\delta>0$ -- recall that the boundary normal derivative of the classical Poisson kernel in the upper half-space decays at the rate $\tau^{-n}$, so this estimate is close to being optimal.    
    On the other hand, letting $k_i \geq 3$ denote the degree of the leading-order homogeneous harmonic polynomial of $u$ at $z_i$, we have 
    \begin{equation}\label{tmp:der-near-pole}
        |\partial_\nu u(\zeta)| = |\nabla u(\zeta)| \lesssim \tau^{k_i-1} 
    \end{equation} 
    as long as $\tau<2\epsilon$ is sufficiently small.
    Combining \eqref{eq:Pkbd-near-pole} and \eqref{tmp:der-near-pole}, we get
    \begin{equation}\label{eq:integ-Pk-near-pole}
        \int_{\partial\Om \cap B_{2\epsilon}(z_i)} |\partial_\nu u \, \partial_\nu P_{z_i} | \lesssim \int_0^{2\epsilon} \tau^{k_i-1} \tau^{-1-\frac{n}{q}} \tau^{n-2} \, d\tau < +\infty, 
    \end{equation} 
    since the exponent 
    \[ k_i-\frac{n}{q} +n - 4 > k_i - 4 \geq -1.\]
    We remark that \eqref{eq:integ-Pk-near-pole} fails when $k_i=2$ at $z_i$, even if we choose $q$ to be close to the optimal value $\frac{n}{n-1}$.
    Hence 
    \begin{equation}\label{tmp:limit-rhs}
        \int_{\partial\Om} \rho_{s_\epsilon} \, \partial_\nu u \partial_\nu P_{z_i} \to \int_{\partial\Om} \rho_{s} \, \partial_\nu u \partial_\nu P_{z_i}, \qquad \text{ as } \epsilon \to 0. 
    \end{equation} 

    On the other hand,
    \begin{align*}
        |(\Delta+\lambda)(s_\epsilon \cdot \nabla u) - \chi_\epsilon (\Delta+\lambda) (s\cdot \nabla u)| & \leq 2|\nabla (s\cdot \nabla u)| \cdot |\nabla \chi_\epsilon| + |s\cdot \nabla u| \cdot |\Delta \chi_\epsilon| \\
        & \lesssim \epsilon^{k_i-2} \cdot \epsilon^{-1} + \epsilon^{k_i-1} \cdot \epsilon^{-2} \\
        & \lesssim \epsilon^{k_i-3},
    \end{align*}
    where we use the bound on $|\nabla u|, |\nabla^2 u|$ on the support of the derivatives of $\chi_\epsilon$, which is $B_{2\epsilon}(z_i) \setminus B_{\epsilon}(z_i)$. Together with the $L^q$-integrability of $P_{z_i}$ and $k_i\geq 3$, this implies
    \begin{equation}\label{tmp:limit-lhs}
        \int_\Om P_{z_i} (\Delta+ \lambda) (s_\epsilon \cdot \nabla u) \to \int_\Om P_{z_i} (\Delta+ \lambda) (s \cdot \nabla u), \qquad \text{ as } \epsilon \to 0. 
    \end{equation} 
    Letting $\epsilon \to 0$ in \eqref{tmp:Pk-cutoff-nearpole} and using \eqref{tmp:limit-rhs} and \eqref{tmp:limit-lhs}, we obtain \eqref{tmpeq:zi-3} for all $s$, finishing the proof of the claim.
\end{proof}

From the above discussions, the following map
\[ \ell_{z_i}: s \in C^{2,\alpha}(\overline{\Om}, \bR^2) \longmapsto \partial_\nu v_s(z_i) =  - \int_{\partial\Om} \rho_s \,\partial_\nu u \, \partial_\nu P_{z_i} \in \bR \]
is a continuous linear functional. Furthermore, recall by the unique continuation and the discussion above that $\partial_\nu u \, \partial_\nu P_{z_i}$ is a non-trivial integrable function on $\partial\Om$. Thus $\ell_{z_i}$ is not identically zero, and 
\[ \ker \ell_{z_i} = \{ s\in C^{2,\alpha}(\overline\Om, \bR^2): \text{ its restriction on the boundary } s\cdot \nu \perp \partial_\nu u \, \partial_\nu P_{z_i} \}. \]
From now on, we let 
\begin{equation}\label{asmp:s}
    s\in C^{2,\alpha}(\overline{\Om}, \bR^2) \setminus  \bigcup_{i=1}^D \ker \ell_{z_i} . 
\end{equation} 

\textbf{Properties of eigenfunctions on the perturbed domains}. Consider a one-parameter family of $C^{2,\alpha}$ diffeomorphisms $h_t = \id + ts$ on $\overline{\Om}$ and let $ \Om_t:= h_t(\Om)$. Since the Dirichlet eigenvalues depend continuously on smooth perturbations of the domains, in a small enough neighborhood of the simple eigenvalue $\lambda$ (depending on the spectral gap near $\lambda$), $\Omega_t$ has a unique simple eigenvalue, denoted by $\lambda_t$, as long as $|t|$ is sufficiently small. Let $\tilde u_t \in H^2 \cap H_0^1(\Om_t)$ denote the Dirichlet eigenfunction of $\Om_t$ with eigenvalue $\lambda_t$, and let $u_t = \tilde u_t \circ h_t\in H^2 \cap H_0^1(\Om)$ be its pullback function on $\Om$. We normalize so that each $u_t$ has unit $L^2$ norm, i.e.\ $u_t \in S(\Om)$.
Then 
\begin{equation}\label{eq:ef-pertb}
    \varphi(u_t, \lambda_t, h_t) := (L_{h_t} + \lambda_t) u_t = (\Delta_{\Om_t} + \lambda_t) \tilde u_t \circ h_t = 0. 
\end{equation} 
Let 
\begin{equation}\label{eq:cv-tang}
    w_t = \frac{u_t-u}{t}, \qquad \delta_t = \frac{\lambda_t - \lambda}{t}. 
\end{equation} 
Moreover, recall by Lemma \ref{lm:unique-triple} that the simplicity assumption implies the tangent space $T_{(u, \la, \id)} Q = \ker d\varphi|_{(u, \la, \id)}$ has a one-to-one correspondence with $s \in C^{2,\alpha}(\overline{\Om}, \bR^2)$. That is, for each $s$ fixed, there is a unique element $(v,\eta,s) \in T_{(u, \la, \id)} Q$. Together with \eqref{eq:ef-pertb} and $h_t=\id+ts$, this implies
\begin{equation}\label{CV:tangent}
    (w_t, \delta_t, s) \longrightarrow (v,\eta, s), \qquad \text{ as } t\to 0. 
\end{equation} 
A priori the above convergence is in $L^2(\Om) \times \bR$. However, since the map $h_t \in C^{2,\alpha} \mapsto L_{h_t} \in \mathcal{L}(C^{2,\alpha}, C^{0,\alpha})$ is continuously differentiable, together with the simplicity of the eigenspace, the convergence in \eqref{CV:tangent} can be upgraded to $C^{2,\alpha}(\Omega) \times \mathbb{R}$.

Our goal is to show that for $|t|>0$ sufficiently small, $\tilde u_t$ has no degenerate boundary points. Firstly, we claim that if $p_t \in \partial\Om_t$ is a degenerate boundary point for $\tilde u_t$, then the corresponding point $q_t := h_t^{-1}(p_t)$ on $ \partial\Om$ satisfies
\begin{equation}\label{asmp:degenerate}
    \partial_\nu u_t(q_t) = 0, \qquad \partial_\tau \partial_\nu u_t(q_t) = 0. 
\end{equation} 
The proof follows from \eqref{eq:change-normalder} and \eqref{eq:change-normalder-2}. Since $h_t$ is a diffeomorphism from $\overline{\Om}$ to $\overline{\Om_t}$, there are real numbers $c(q_t), \bar c(q_t) > 0$ such that
\[ \left(Dh_t(q_t)\right)^T \nu_t(p_t) = c(q_t) \nu(q_t), 
\qquad Dh_t(q_t) \tau(q_t) = \bar c(q_t) \tau_t(p_t), \]
where $\nu, \nu_t$ denote the unit outer normal vectors of $\partial\Om, \partial\Om_t$, and $\tau, \tau_t$ denote  unit tangent vectors for $\partial\Om, \partial\Om_t$ ($\tau_t$ is chosen so that $\bar c(q_t)$ above is positive). Then
\[ \partial_\nu u_t(q_t) = c(q_t) \partial_{\nu_t} \tilde u_t(p_t), \]
\[ \partial_\tau \partial_\nu u_t(q_t) = c(q_t) \bar{c}(q_t) \cdot  \partial_{\tau_t} \partial_{\nu_t} \tilde u_t(p_t) = \det Dh_t(q_t) \cdot \partial_{\tau_t} \partial_{\nu_t} \tilde u_t(p_t). \]
Hence the assumption $\partial_{\nu_t} \tilde u_t(p_t)=0$ and $\partial_{\tau_t} \partial_{\nu_t} \tilde u_t(p_t) = 0$ implies \eqref{asmp:degenerate}.

Next, we show that \eqref{asmp:degenerate} is impossible near each $z_i$. For each $z_i$, let $k_i>2$ denote the degree of the leading-order homogeneous harmonic polynomial of $u$ at $z_i$, and let $\mathcal{O}_i$ denote the neighborhood so that the mixed derivative estimate in Claim \ref{cl:nondeg-mix-der} holds for $u$. Recall by the assumption \eqref{asmp:s}, we know $\partial_\nu v(z_i) \neq 0$ for all $i=1, \cdots, D$. By the $C^1$ continuity of $v$ and by making $\mathcal{O}_i$ smaller if necessary, we have
\begin{equation}\label{asmp:nondeg}
    J := \min_{i=1, \cdots, D} \inf_{\mathcal{O}_i} |\partial_\nu v| >0. 
\end{equation} 
Assume for contradiction that $q_t \in \mathcal{O}_i$.
Then
\[ 0 = \partial_\tau \partial_\nu u_t(q_t) = \partial_\tau \partial_\nu u(q_t) + t \partial_\tau \partial_\nu w_t(q_t). \]
It follows from \eqref{CV:tangent} and Claim \ref{cl:nondeg-mix-der} that
\begin{equation}\label{tmp:mix-der}
    A|q_t - z_i|^{k_i-2} \leq |\partial_\tau \partial_\nu u(q_t)| \leq |t| \cdot |\partial_\tau \partial_\nu w_t(q_t)| \leq 2|t|\cdot \|v\|_{C^2(\Om)} \leq C|t| \cdot \|s\|_{C^{2, \alpha}}. 
\end{equation} 
On the other hand, we also have
\[ 0 = \partial_\nu u_t(q_t) = \partial_\nu u(q_t) + t \partial_\nu w_t(q_t). \]
By  \eqref{CV:tangent} and \eqref{asmp:nondeg}, we have $|\partial_\nu w_t(q_t)| \geq J/2>0$. By the expansion of $u$ and the $C^1$ error estimate in Proposition \ref{prop:bdexp}, we have
\[ |\partial_\nu u(q_t)| = |\nabla u(q_t)| \leq C |q_t - z_i|^{k_i -1}. \]
Combined, we have
\begin{equation}\label{tmp:1st-der}
    \frac{J}{2} |t| \leq |t \partial_\nu w_t(q_t)| = |\partial_\nu u(q_t)| \leq C |q_t - z_i|^{k_i -1}. 
\end{equation} 
The estimates \eqref{tmp:mix-der} and \eqref{tmp:1st-der} together imply
\[ |t| \lesssim |q_t-z_i|^{k_i-1} \lesssim |t|^{\frac{k_i-1}{k_i-2}}, \]
with constants depending on $s$ and $u$. This is impossible if $|t|$ is chosen sufficiently small. Therefore \eqref{asmp:degenerate} is impossible in $\mathcal{O}_i$.

On $\partial\Om \setminus \cup_{i=1}^D \mathcal{O}_i$, either $|\partial_\nu u|$ has a positive lower bound, which implies $|\partial_\nu u_t|$ has a (smaller) positive lower bound if $|t|$ is chosen small; or $\partial_\nu u(z) = 0$ but $\partial_\tau \partial_\nu u(z) \neq 0$. In the latter case, the mixed derivative $|\partial_\tau \partial_\nu u_t|$ has a positive lower bound, if $|t|$ is chosen small. Therefore \eqref{asmp:degenerate} is impossible on all of $\partial\Om$, and thus $\tilde u_t$ does not have degenerate boundary points on $\Om_t$.

\textbf{Brief sketch of the interior case}. For the interior case, the aim is to show that $u_t$ has no interior singular points, such that
\begin{equation}\label{eq:t-int-sing}
	u_t(z_t) = 0, \qquad \nabla u_t(z_t) = 0. 
\end{equation} 

Suppose $x \in \Om$ is an interior singular point of $u$, that is $u(x)=0$ and $\nabla u(x) = 0$. Then by the asymptotic expansion of $u$ near $x$ by \cite{Bers}, we can show that the degree of the leading-order homogeneous harmonic polynomial satisfies $m\geq 2$, and 
\[ |\nabla u(\zeta)| \gtrsim |\zeta - x|^{m-1}, \qquad \text{ in a neighborhood of } x. \]
The analysis of the boundary expansion above shows that the interior singular points cannot accumulate at $\partial\Om$. Together, they imply that each fixed eigenfunction $u$ has finitely many interior singular points $x_1, \cdots, x_S \in \Om$. 

Assume that $s$ is an admissible perturbation in \eqref{asmp:s}, such that the unique element $(v_s, \eta_s, s) \in T_{(u,\lambda, \id)} Q$ satisfies the additional conditions 
\begin{equation}\label{cond:interior}
	v_s(x_j) \neq 0, \qquad \text{ for all } j=1, \cdots, S.
\end{equation}
 Then by an analogous argument as the proof of \eqref{asmp:degenerate}, we show that \eqref{eq:t-int-sing} is impossible near each $x_j$. Near every $x\in \overline\Omega$ such that either $|u(x)|\geq c_0$ or $|\nabla u(x)| \geq c_0 $, for some fixed value of $c_0>0$, we can choose a uniform constant $t_0 = t_0(c_0, s)>0$ such that \eqref{eq:t-int-sing} is impossible for all $0<|t|<t_0$.
It remains to rule out \eqref{eq:t-int-sing} near each boundary point $q\in \partial\Omega$ such that $u(q)=|\nabla u(q)|=0$. We remark that these boundary singular points are also isolated, by Proposition \ref{prop:bdexp}. Assume without loss of generality that $q$ is the origin. Then the leading-order homogeneous harmonic polynomial of $u$ at the origin is $P(z) = c \Im z^k$ for some $k\geq 2$ and $c\in \R$. 

\setcounter{footnote}{0}
\renewcommand{\thefootnote}{\fnsymbol{footnote}}
For any $z\in \Omega$ near the origin, let $d_z:= \dist(z, \partial\Omega)$ denote its distance to the boundary. Then by Proposition \ref{prop:bdexp} and the assumption that $u$ vanishes on the boundary of the $C^2$ domain $\Omega$\footnote{The latter assumption is used to get the desired estimates when $d_z \ll |z|$.}, we have
\begin{equation}\label{tmp:bdonu}
	\frac{|u(z)|}{d_z} = \frac{|P(z) + \psi(z)|}{d_z} \lesssim |z|^{k-1}, \qquad \frac{|\nabla u(z)|}{d_z} = \frac{|\nabla P(z) + \nabla \psi(z)|}{d_z} \gtrsim |z|^{k-2}. 
\end{equation} 
Assuming \eqref{eq:t-int-sing} holds for some $z_t$ near the origin, then
\begin{equation}\label{tmp:int-sing-near-bd}
	0 = \frac{u_t(z_t)}{d_{z_t}} = \frac{u(z_t) + t w_t(z_t) }{d_{z_t}}, \qquad 0 = \frac{\nabla u_t(z_t)}{d_{z_t}} = \frac{\nabla u(z_t) + t \nabla w_t(z_t) }{d_{z_t}}. 
\end{equation}

Recall \eqref{CV:tangent} that $w_t \to v$ in $C^{2,\alpha}$ up to the boundary; moreover, they all vanish on the boundary. Let $p\in \partial\Omega$ be such that $|z-p|= d_z$, and let $\nu_p$ denote the unit normal vector at $p$ pointing inwards, so that the line segment $[p,z]$ in the $\nu_p$-direction lies in $\Omega$.  We have
\begin{equation}\label{tmp:bd-comp} 
	\frac{w_t(z)}{d_z} - \frac{v(z)}{d_z} = \frac{w_t(z)-w_t(p)}{d_z} - \frac{v(z)-v(p)}{d_z} = \int_0^1 \left[ \partial_{\nu_p} w_t(p+\rho d_z) - \partial_{\nu_p} v(p+\rho d_z) \right] d\rho,  
\end{equation}
which tends to $0$ uniformly in $z$ as $t\to 0$.
Moreover, by perturbing $z$ in the tangential direction (i.e. in the direction of $\tau_p$ orthogonal to $\nu_p$) \eqref{tmp:bd-comp} gives 
\[ \left|\frac{\partial_{\tau_p} w_t(z)}{d_z} - \frac{\partial_{\tau_p} v(z)}{d_z} \right| \lesssim \|w_t - v\|_{C^2}, \]
and by perturbing $z$ in the normal $\nu_p$ direction, \eqref{tmp:bd-comp} gives
\[ \left|\frac{\partial_{\nu_p} w_t(z)}{d_z} - \frac{\partial_{\nu_p} v(z)}{d_z} \right| \lesssim \|w_t - v\|_{C^2}. \]
Combined, we obtain
\begin{equation}\label{tmp:bd-comp-grad}
	\left|\frac{\nabla w_t(z)}{d_z} - \frac{\nabla v(z)}{d_z} \right| \lesssim \|w_t - v\|_{C^2} \to 0, \qquad \text{ as } t\to 0. 
\end{equation} 

If the origin is one of the degenerate boundary points $z_i$, then the assumption \eqref{asmp:s} on $s$, together with the $C^{1,\alpha}$ continuity of $v$ up to the boundary, guarantees that $\frac{|v(z)|}{d_z} \geq \frac{J}{2} >0$ nearby, where $J$ is defined in \eqref{asmp:nondeg}. Therefore, combining \eqref{tmp:bd-comp}, \eqref{tmp:int-sing-near-bd} and  \eqref{tmp:bdonu}, we get
\begin{equation}\label{tmp:int-1}
	\frac{J}{3} |t| \leq |t| \frac{|w_t(z_t)|}{d_{z_t}} = \frac{|u(z_t)|}{d_{z_t}} \lesssim |z_t|^{k-1}; 
\end{equation}  
combining \eqref{tmp:bdonu}, \eqref{tmp:int-sing-near-bd} and \eqref{tmp:bd-comp-grad}, we get
\begin{equation}\label{tmp:int-2}
	|z_t|^{k-2} \lesssim \frac{|\nabla u(z_t)|}{d_{z_t}} = |t| \frac{|\nabla w_t(z_t)|}{d_{z_t}} \lesssim |t| \left(\|v\|_{C^2} + \|w_t - v\|_{C^2} \right) \lesssim |t| \|v\|_{C^2}.  
\end{equation} 
These two inequalities \eqref{tmp:int-1}, \eqref{tmp:int-2} can not hold at the same time, if $|t|$ is chosen sufficiently small (depending on $J$ and $\|v\|_{C^2} \lesssim \|s\|_{C^{2,\alpha}}$).
If the origin is not a degenerate boundary point, then necessarily its leading-order degree $k=2$. Thus \eqref{tmp:bdonu} gives a definite lower bound of $\frac{|\nabla u(z)|}{d_z}$. Combining \eqref{tmp:int-sing-near-bd} and \eqref{tmp:bd-comp-grad}, we can guarantee $\frac{\nabla u_t(z_t)}{d_{z_t}} \neq 0$ near the origin by choosing $|t|$ sufficiently small (depending on the perturbation $s$ and $c$, the coefficient for $P$). To summarize, we have identified $t_0>0$ such that for all $0<|t|<t_0$, the eigenfunction $u_t$ has no interior singular points near every boundary and interior singular points of $u$. Away from these points, either $|u|$ or $|\nabla u|$ has a uniform lower bound $c_0>0$ in $\overline{\Omega}$, so by choosing $t_0$ smaller, we guarantee that $u_t$ has no interior singular points in all of $\Omega$.

It remains to find perturbations $s$ that satisfy the conditions \eqref{cond:interior}.
As in Subsection \ref{sec:interior-reg}, we multiply the equation \eqref{tmp:vs} for $v_s$ by $G_{x_j}$, the reduced Green's function with pole at $x_j$, and integrate. Using the assumption $u(x_j)=|\nabla u(x_j)|=0$, and the properties of the reduced Green's function, we obtain
\begin{equation}
    v_s(x_j)  = \eta_s \int_\Om u G_{x_j} + \int_{\partial\Om} \rho_s \, \partial_\nu u \, \partial_\nu G_{x_j} = \int_{\partial\Om} \rho_s \, \partial_\nu u \, \partial_\nu G_{x_j}. 
\end{equation}
This implies that the continuous linear map $\mu_{x_j}: s\in C^{2,\alpha}(\overline{\Om}, \bR^2) \mapsto v_s(x_j) \in \bR$ is not identically zero. Finally, we let
\begin{equation}\label{def:good-perturbation}
    \mathcal{A}_\lambda := C^{2,\alpha}(\overline{\Om}, \bR^2) \mathbin{\big\backslash} \left( \bigcup_{j=1}^S \ker \mu_{x_j} \cup \bigcup_{i=1}^D \ker \ell_{z_i}\right). 
\end{equation} 
Since we only remove finitely many hyperplanes, this is clearly open and dense in $C^{2,\alpha}(\overline{\Om}, \bR^2)$. This finishes the proof of Proposition \ref{prop:2D}.

The simplest example to illustrate the above perturbation is the following. Assume the origin is a degenerate boundary point on $\partial\Om$, $\partial\Om$ is tangent to the $x$-axis at the origin, and $\partial_\nu u|_{\partial\Om} = x^2$ near the origin (this is the case when $u(z)=-\frac13 \Im z^3 + o(|z|^3)$ near the origin). Assume also that $\partial_\nu v_s(0)=-1 \neq 0$. Then
\[ \partial_\nu u_t|_{\partial\Om} \sim  x^2 - t \qquad \text{ near the origin}. \]
In particular, when $t<0$, we always have $\partial_\nu u_t > 0$ and geometrically, the nodal curve of $u_t \circ h_t^{-1}$ is closed; when $t>0$, we have $\partial_\nu u_t=0$ at two points $x \sim \pm \sqrt{t}$ and geometrically, the nodal curve of $u_t \circ h_t^{-1}$ touches the boundary orthogonally at two points near the origin. The simple eigenfunctions $u_t$ of the one-parameter family $\Om_t$ have a bifurcation behavior at $t=0$.

\bibliographystyle{alpha}
\bibliography{bibliography}

@article{Bers,
    author = {Bers, L.},
    title = {Local behavior of solutions of general linear elliptic equations},
    journal = {Comm. Pure Appl. Math.},
    year = {1955},
    volume = {8},
    pages = {473-496},
    url = {https://doi.org/10.1002/cpa.3160080404}
}

@article{HartmanWintner,
  title={On the local behavior of solutions of non-parabolic partial differential equations: {III}. {A}pproximations by spherical harmonics},
  author={Hartman, P. and Wintner, A.},
  journal={Am. J. Math.},
  volume={77},
  number={3},
  pages={453--474},
  year={1955},
  publisher={JSTOR}
}

@article{ADN,
    author = {Agmon, S. and Douglis, A. and Nirenberg, L.},
    title = {Estimates near the boundary for solutions of elliptic partial differential equations satisfying general boundary conditions. {I}},
    journal = {Comm. Pure Appl. Math.},
    year = {1959},
    volume = {12},
    page = {623-727},
    url = {https://doi.org/10.1002/cpa.3160120405}
}

@article{HS,
    author = {Hardt, R. and Simon, L.},
    title = {Nodal sets for solutions of elliptic equations},
    journal = {J. Differential Geom.},
    year = {1989},
    volume = {30},
    pages = {505-522},
    url = {https://doi.org/10.4310/jdg/1214443599}
}

@article{KZexpansion,
    author = {C. Kenig and Z. Zhao},
    title = {Expansion of harmonic functions near the boundary of {D}ini domains},
    journal = {Rev. Mat. Iberoam.},
    year = {2022},
    volume = {38},
    pages = {2117-2152}
}

@article{UhlenbeckGeneric,
  author = {Uhlenbeck, K.},
  title={Generic properties of eigenfunctions},
  journal={Am. J. Math.},
  volume={98},
  number={4},
  pages={1059--1078},
  year={1976},
  publisher={JSTOR}
}

@article{PereiraHenry22,
  title={Dan {H}enry’s work on perturbation of the boundary problems},
  author={Pereira, A. L.},
  journal={S{\~a}o Paulo J. Math Sci.},
  volume={16},
  number={1},
  pages={157--170},
  year={2022},
  publisher={Springer}
}

@book{HenryPerturbationbook05,
  title={Perturbation of the boundary in boundary-value problems of partial differential equations},
  author={Henry, D.},
  series = {London Mathematical Society Lecture Note Series},
  volume={318},
  year={2005},
  publisher={Cambridge University Press}
}

@article{KZcritical,
  title={A note on the critical set of harmonic functions near the boundary},
  author={Kenig, C. and Zhao, Z.},
  journal={Trans. Am. Math. Soc.},
  volume={378},
  number={05},
  pages={3721--3754},
  year={2025}
}

@article {Micheletti-1,
    AUTHOR = {Micheletti, A. M.},
     TITLE = {Perturbazione dello spettro dell'operatore di {L}aplace, in
              relazione ad una variazione del campo},
   JOURNAL = {Ann. Scuola Norm. Sup. Pisa Cl. Sci. (3)},
  FJOURNAL = {Annali della Scuola Normale Superiore di Pisa. Classe di
              Scienze. Serie III},
    VOLUME = {26},
      YEAR = {1972},
     PAGES = {151--169},
      ISSN = {0391-173X},
   MRCLASS = {35P05},
  MRNUMBER = {367480},
MRREVIEWER = {W.\ M.\ Greenlee},
}

@article {Micheletti-2,
    AUTHOR = {Micheletti, A. M.},
     TITLE = {Metrica per famiglie di domini limitati e propriet\`a{}
              generiche degli autovalori},
   JOURNAL = {Ann. Scuola Norm. Sup. Pisa Cl. Sci. (3)},
  FJOURNAL = {Annali della Scuola Normale Superiore di Pisa. Classe di
              Scienze. Serie III},
    VOLUME = {26},
      YEAR = {1972},
     PAGES = {683--694},
      ISSN = {0391-173X},
   MRCLASS = {35P99 (35J05)},
  MRNUMBER = {377306},
MRREVIEWER = {W.\ M.\ Greenlee},
}

@article{Albert,
 ISSN = {00029947},
 URL = {http://www.jstor.org/stable/1997812},
 author = {Jeffrey H. Albert},
 journal = {Transactions of the American Mathematical Society},
 pages = {341--354},
 publisher = {American Mathematical Society},
 title = {Generic Properties of Eigenfunctions of Elliptic Partial Differential Operators},
 urldate = {2026-07-13},
 volume = {238},
 year = {1978}
}

@article{SpectralPartition,
    author = {Helffer, B. and Hoffmann-Ostenhof, T. and Terracini, S.},
    title = {Nodal Domains and Spectral Minimal Partitions},
    journal = {Annales de l'Institut Henri Poincaré. C, Analyse non linéaire},
    pages = {101-138},
    volume = {26},
    url = {10.1016/j.anihpc.2007.07.004},
    year = {2009}
}

@article{Alessandrini,
    author = {Alessandrini, G.},
    title = {Nodal lines of eigenfunctions of the fixed membrane problem in general convex domains},
    journal = {Commentarii Mathematici Helvetici},
    volume = {69}, 
    pages = {142-154},
    year = {1994}
}

@article{HigherHarnack,
    author = {De Silva, D. and Savin, O.},
    title = {A note on higher regularity boundary {H}arnack inequality},
    journal = {Discrete Contin. Dyn. Syst.},
    volume = {35},
    pages = {6155–6163},
    year = {2015}
}

@book{GT,
	Author = {Gilbarg, D. and Trudinger, N. S.},
	Edition = {Reprint of the 1998},
	Isbn = {3-540-41160-7/pbk},

	Msc2010 = {35-02 35J65 35B45 35J25 35B50 35B05 47H10},
	Pages = {xiii + 517},
	Publisher = {Berlin: Springer},
	Title = {{Elliptic partial differential equations of second order}},
	Year = {2001},
	Zbl = {1042.35002}}

@book{gilbarg_trudinger_1977,
  author    = {Gilbarg, D. and Trudinger, N.},
  title     = {Elliptic Partial Differential Equations of Second Order},
  series    = {Grundlehren der mathematischen Wissenschaften},
  volume    = {224},
  publisher = {Springer-Verlag},
  address   = {Berlin, Heidelberg},
  year      = {1977},
  edition   = {1st},
  pages     = {x+401},
  doi       = {10.1007/978-3-642-96379-7}
}

\end{document}